\documentclass[12pt,a4paper]{article}

\usepackage{amssymb,amsmath,amsfonts,amstext}
\usepackage{bbm}
\usepackage{geometry}	
\usepackage{graphicx}
\usepackage[colorlinks=true,breaklinks=true]{hyperref}
\usepackage{ifthen}
\usepackage[utf8]{inputenc}
\usepackage[bbgreekl]{mathbbol}
\usepackage{multirow}
\usepackage[l2tabu,orthodox]{nag}
\usepackage[numbers,square]{natbib}
\usepackage{tikz}
\usepackage{tikz-cd}
\usepackage{url}

\usetikzlibrary{arrows,matrix,shapes,positioning}
\usetikzlibrary{calc,patterns,shapes.geometric}
\usetikzlibrary{decorations.markings,decorations.pathmorphing,decorations.pathreplacing}
\usetikzlibrary{decorations.text,decorations.shapes}

\newcommand*{\doi}[1]{\href{http://dx.doi.org/\detokenize{#1}}{doi: \detokenize{#1}}}

\numberwithin{equation}{section}

\let\div\undefined

\DeclareMathOperator{\div}{div}
\DeclareMathOperator{\curl}{curl}
\DeclareMathOperator{\grad}{grad}

\DeclareMathOperator{\Ima}{Im}
\DeclareMathOperator{\Ker}{Ker}

\DeclareMathOperator{\spn}{span}
\DeclareMathAlphabet{\mathpzc}{OT1}{pzc}{m}{it}

\newcommand{\dist}{\ensuremath{\mathpzc{f}}}
\newcommand{\ham}{\ensuremath{\mathcal{H}}}
\newcommand{\jac}{\ensuremath{\mathbb{J}}}
\newcommand{\cas}{\ensuremath{\mathcal{C}}}
\newcommand{\mom}{\ensuremath{\mathcal{P}}}
\newcommand{\ff}{\ensuremath{\mathcal{F}}}
\newcommand{\fg}{\ensuremath{\mathcal{G}}}
\newcommand{\fh}{\ensuremath{\mathcal{H}}}
\newcommand{\dff}{\ensuremath{F}}
\newcommand{\dfg}{\ensuremath{G}}
\newcommand{\dham}{\ensuremath{H}}
\newcommand{\dcas}{\ensuremath{C}}

\newcommand{\dd}{\ensuremath{\mathrm{d}}}
\newcommand{\DD}{\ensuremath{\mathrm{D}}}
\newcommand{\dt}{\ensuremath{\,\dd t}}
\newcommand{\dx}[1][\empty]{\ensuremath{\,\dd \ifthenelse{\equal{#1}{\empty}}{\xx}{\ifthenelse{\equal{#1}{1}}{x}{x_{#1}}}}}
\newcommand{\dv}[1][\empty]{\ensuremath{\,\dd \ifthenelse{\equal{#1}{\empty}}{\vv}{v_{#1}}}}
\newcommand{\ext}{\ensuremath{\mathsf{d}}}
\newcommand{\pa}{\ensuremath{\partial}}
\newcommand{\Dt}{\ensuremath{\Delta t}}

\renewcommand{\vec}[1]{\ensuremath{\mathbf{#1}}}
\newcommand{\coeff}[1]{\ensuremath{\lowercase{#1}}}
\newcommand{\abs}[1]{\ensuremath{\left \vert #1 \right \vert}}
\newcommand{\cbra}[2]{\ensuremath{\left \{ #1 , #2 \right \}}}
\newcommand{\sbra}[2]{\ensuremath{\left  [ #1 , #2 \right  ]}}
\newcommand{\pair}[2]{\ensuremath{\left < #1 \, , \, #2 \right >}}
\newcommand{\mcal}[1]{\ensuremath{\mathcal{#1}}}
\newcommand{\order}[1]{\ensuremath{\mcal{O} (#1)}}
\newcommand{\fracd}[2]{\ensuremath{\frac{\delta   #1}{\delta   #2}}}
\newcommand{\fracp}[2]{\ensuremath{\frac{\partial #1}{\partial #2}}}

\newcommand{\uu}{\ensuremath{\vec{u}}}
\newcommand{\vv}{\ensuremath{\vec{v}}}

\newcommand{\xx}{\ensuremath{\vec{x}}}

\newcommand{\zz}{\ensuremath{\vec{z}}}

\newcommand{\cu}{\ensuremath{\mathrm{u}}}
\newcommand{\cv}{\ensuremath{\mathrm{v}}}

\newcommand{\cx}{\ensuremath{\mathrm{x}}}

\newcommand{\xa}{\ensuremath{\xx_{a}}}
\newcommand{\va}{\ensuremath{\vv_{a}}}

\newcommand{\f}{\ensuremath{\dist}}
\newcommand{\A}{\ensuremath{\vec{A}}}
\newcommand{\B}{\ensuremath{\vec{B}}}
\newcommand{\E}{\ensuremath{\vec{E}}}
\newcommand{\J}{\ensuremath{\vec{J}}}
\newcommand{\X}{\ensuremath{\vec{X}}}
\newcommand{\V}{\ensuremath{\vec{V}}}

\newcommand{\vara}{\ensuremath{\mathrm{a}}}

\newcommand{\varae}{\ensuremath{\mathrm{f}}}

\newcommand{\ca}{\ensuremath{\vec{\coeff{A}}}}
\newcommand{\cb}{\ensuremath{\vec{\coeff{B}}}}
\newcommand{\cc}{\ensuremath{\vec{\coeff{C}}}}
\newcommand{\cd}{\ensuremath{\vec{\coeff{D}}}}
\newcommand{\ce}{\ensuremath{\vec{\coeff{E}}}}

\newcommand{\cj}{\ensuremath{\vec{\coeff{J}}}}
\newcommand{\cae}{\ensuremath{\boldsymbol{f}}}
\newcommand{\cphi}{\ensuremath{\boldsymbol{\varphi}}}
\newcommand{\crho}{\ensuremath{\boldsymbol{\varrho}}}

\newcommand{\MA}{\ensuremath{\mathbb{A}}}
\newcommand{\BB}{\ensuremath{\mathbb{B}}}
\newcommand{\MM}{\ensuremath{\mathbb{M}}}

\newcommand{\C}{\ensuremath{\mathbb{C}}}
\newcommand{\D}{\ensuremath{\mathbb{D}}}
\newcommand{\G}{\ensuremath{\mathbb{G}}}

\newcommand{\R}{\ensuremath{\mathbb{R}}}

\newcommand{\Lab}{\ensuremath{\boldsymbol{\Lambda}}}
\newcommand{\LaB}{\ensuremath{\mathbb{\Lambda}}}
\newcommand{\sib}{\ensuremath{\boldsymbol{\sigma}}}
\newcommand{\Sib}{\ensuremath{\boldsymbol{\Sigma}}}
\newcommand{\omb}{\ensuremath{\boldsymbol{\omega}}}

\newcommand{\rsp}{\ensuremath{\mathbb{R}}}

\newcommand{\de}{\ensuremath{\delta}}
\newcommand{\De}{\ensuremath{\Delta}}

\newcommand{\ve}{\ensuremath{\varepsilon}}

\newcommand{\la}{\ensuremath{\la}}
\newcommand{\La}{\ensuremath{\Lambda}}

\newcommand{\phy}{\ensuremath{\varphi}}

\newcommand{\si}{\ensuremath{\sigma}}
\newcommand{\Si}{\ensuremath{\Sigma}}

\newcommand{\vb}{\mathbf{v}}
\newcommand{\eb}{\mathbf{e}}
\newcommand{\jb}{\mathbf{j}}
\newcommand{\Eb}{\mathbf{E}}

\newcommand{\bb}{\mathbf{b}}
\newcommand{\db}{\mathbf{d}}

\title{GEMPIC: Geometric ElectroMagnetic\\Particle-In-Cell Methods} 

\author{
\vspace{.5em}\\
\large{Michael Kraus$^{1,2}$},
\large{Katharina Kormann$^{1,2}$}, \\
\large{Philip~J.~Morrison$^{3}$},
\large{Eric Sonnendr\"ucker$^{1,2}$}
\vspace{.5em}\\
$^1$\normalsize{Max-Planck-Institut f\"ur Plasmaphysik}\\
\normalsize{Boltzmannstra\ss{}e 2, 85748 Garching, Deutschland}%
\vspace{.5em}\\
$^2$\normalsize{Technische Universit\"at M\"unchen, Zentrum Mathematik}\\
\normalsize{Boltzmannstra\ss{}e 3, 85748 Garching, Deutschland}%
\vspace{.5em}\\
$^3$\normalsize{Department of Physics and Institute for Fusion Studies}\\
\normalsize{The University of Texas at Austin, Austin, TX, 78712, USA}%
\vspace{1em}\\
}

\date{1. June 2017}

\begin{document}

\maketitle

\begin{abstract}

We present a novel framework for Finite Element Particle-in-Cell methods based on the discretization of the underlying Hamiltonian structure of the Vlasov--Maxwell system. We derive a semi-discrete Poisson bracket, which retains the defining properties of a bracket, anti-symmetry and the Jacobi identity, as well as conservation of its Casimir invariants, implying that the semi-discrete system is still a Hamiltonian system.
In order to obtain a fully discrete Poisson integrator, the semi-discrete bracket is used in conjunction with Hamiltonian splitting methods for integration in time.
Techniques from Finite Element Exterior Calculus ensure conservation of the divergence of the magnetic field and Gauss' law as well as stability of the field solver.
The resulting methods are gauge invariant, feature exact charge conservation and show excellent long-time energy and momentum behaviour.
Due to the generality of our framework, these conservation properties are guaranteed independently of a particular choice of the Finite Element basis, as long as the corresponding Finite Element spaces satisfy certain compatibility conditions.

\end{abstract}

\newpage

\tableofcontents

\newpage

\section{Introduction}
\label{sec:intro}

We consider a structure-preserving numerical implementation of the Vlasov--Maxwell system, which is a system of kinetic equations, describing the dynamics of charged particles in a plasma, coupled to Maxwell's equations describing electrodynamic phenomena arising from the motion of the particles as well as from externally applied fields.
While the design of numerical methods for the Vlasov--Maxwell (and Vlasov--Poisson) system has attracted considerable attention since the early 1960s (see~\cite{Sonnendruecker:2017} and references therein), the systematic development of structure-preserving or geometric numerical methods started only recently.

The Vlasov--Maxwell system exhibits a rather large set of structural properties, which should be considered in the discretization. Most prominently, the Vlasov--Maxwell system features a variational~\cite{Low:1958, Ye:1992, Cendra:1998} as well as a Hamiltonian~\cite{Morrison:1980, Weinstein:1981, Morrison:1982, Marsden:1982} structure. This implies a range of conserved quantities, which by Noether's theorem are related to symmetries of the Lagrangian and the Hamiltonian, respectively. In addition, the degeneracy of the Poisson brackets in the Hamiltonian formulation implies the conservation of several families of so-called Casimir functionals (see e.g.~\cite{Morrison:1998} for a review).

Maxwell's equations have a rich structure themselves. The various fields and potentials appearing in these equations are most naturally described as differential forms~\cite{Bossavit:1990, Warnick:1997, Warnick:2006, BaezMuniain:1994} (see also~\cite{Dray:2014, Morita:2001, Darling:1994}). The spaces of these differential forms build what is called a deRham complex. This implies certain compatibility conditions between the spaces, essentially boiling down to the identities from vector calculus, $\curl \grad = 0$ and $\div \curl = 0$. It has been realized that it is of utmost importance to preserve this complex structure in the discretization in order to obtain stable numerical methods. This goes hand in hand with preserving  two more structural properties provided by the constraints on the electromagnetic fields, namely that the divergence of the magnetic field $\B$ vanishes, $\div \B = 0$, and Gauss' law, $\div \E = \rho$, stating that the divergence of the electromagnetic field $\E$ equals the charge density $\rho$.

The compatibility problems of discrete Vlasov--Maxwell solvers has been widely discussed in the Particle-In-Cell (PIC) literature~\cite{Villasenor.Buneman.1992.cpc,umeda2003,barthelme2005,esirkepov2001exact, eastwood1991virtual,yu2013high} for exact charge conservation. An alternative is to modify Maxwell's equations by adding Lagrange multipliers to relax the constraint~\cite{munz1999, Boris:1970, marder1987,langdon1992,munz2000divergence}.
For a more geometric perspective on charge conservation based on Whitney forms one can refer to~\cite{moon2015exact}.
Even though it has attracted less interest the problem also exists for grid-based discretizations of the Vlasov equations and the same recipes apply there as discussed in~\cite{crouseilles2014,sircombe2009split}.
Note also that the infinite-dimensional kernel of the curl operator has made it particularly hard to find good discretizations for Maxwell's equations, especially for the eigenvalue problem~\cite{Boffi.2006.csd, boffi2010finite, Buffa.Perugia.2006.sinum, Caorsi.Fernandes.Raffetto.2000.sinum, Hesthaven.Warburton.2004.phil}. 

Geometric Eulerian (grid-based) discretizations for the Vlasov--Poisson system have been proposed based on spline differential forms~\cite{Back:2014} as well as variational integrators~\cite{KrausMajSonnendruecker:2015, Kraus:2013:thesis}. While the former guarantees exact local conservation of important quantities like mass, momentum, energy and the $L^{2}$ norm of the distribution function after a semi-discretization in space, the latter retains these properties even after the discretization in time.
Recently, also various discretizations based on discontinuous Galerkin methods have been proposed for both, the Vlasov--Poisson~\cite{deDios:2011, deDios:2012b, deDios:2012a, HeathGamba:2012, ChengGamba:2013, Madaule:2014} and the Vlasov--Maxwell system~\cite{Cheng:2014b, Cheng:2014c, Cheng:2014a}. Even though these are usually not based on geometric principles, they tend to show good long-time conservation properties with respect to momentum and/or energy.

First attempts to obtain geometric semi-discretizations for Particle-in-Cell methods for the Vlasov--Maxwell system have been made by~\citet{Lewis:1970, Lewis:1972}. In his works, Lewis presents a fairly general framework for discretizing Low's Lagrangian~\cite{Low:1958} in space. After fixing the Coulomb gauge and applying a simple Finite Difference approximation to the fields, he obtains semi-discrete, energy and charge-conserving Euler--Lagrange equations. For integration in time the leapfrog method is used.
In a similar way, Evstatiev, Shadwick and Stamm performed a variational semi-discretization of Low's Lagrangian in space, using standard Finite Difference and Finite Element discretizations of the fields and an explicit symplectic integrator in time~\cite{Evstatiev:2013, Shadwick:2014, Stamm:2014}. On the semi-discrete level, energy is conserved exactly but momentum and charge are only conserved in an average sense.

The first semi-discretization of the noncanonical Poisson bracket formulation of the Vlasov--Maxwell system~\cite{Morrison:1980, Weinstein:1981, Morrison:1982, Marsden:1982} can be found in the work of~\citet{Holloway:1996}. Spatial discretizations based on Fourier--Galerkin, Fourier collocation and Legendre-–Gauss-–Lobatto collocation methods are considered. The semi-discrete system is automatically guaranteed to be gauge invariant as it is formulated in terms of the electromagnetic fields instead of the potentials. The different discretization approaches are shown to have varying properties regarding the conservation of momentum maps and Casimir invariants but none preserves the Jacobi identity. It was already noted by~\citet{Morrison:1981} and~\citet{Scovel:1994}, though, that grid-based discretizations of noncanonical Poisson brackets do not appear to inherit a Poisson structure from the continuous problem and \citeauthor{Scovel:1994} suggested that one should turn to particle-based discretizations instead. In fact, for the vorticity equation it was shown by~\citet{Morrison:1981b} that using discrete vortices leads to a semi-discretization that retains the Hamiltonian structure.
Such an integrator for the Vlasov--Amp\`ere Poisson bracket was first presented by~\citet{Evstatiev:2013}, based on a mixed semi-discretization in space, using particles for the distribution function and a grid-based discretization for the electromagnetic fields. However, this work lacks a proof of the Jacobi identity for the semi-discrete bracket, which is crucial for a Hamiltonian integrator.

The first fully discrete geometric Particle-in-Cell method for the Vlasov--Maxwell system has been proposed by~\citet{Squire:2012}, applying a fully discrete action principle to Low's Lagrangian and discretizing the electromagnetic fields via discrete exterior calculus (DEC)~\cite{Hirani:2003, Desbrun:2008, Stern:2014}. This leads to gauge invariant variational integrators that satisfy exact charge conservation in addition to approximate energy conservation.
\citet{Xiao:2015} suggest a Hamiltonian discretization using Whitney form interpolants for the fields. Their integrator is obtained from a variational principle, so that the Jacobi identity is satisfied automatically. Moreover, the Whitney form interpolants preserve the deRham complex structure of the involved spaces, so that the algorithm is also charge-conserving.
\citet{Qin:2016} use the same interpolants to directly discretize the canonical Vlasov--Maxwell bracket~\cite{Marsden:1982} and integrate the resulting finite dimensional system with the symplectic Euler method.
\citet{he2016} introduce a discretization of the noncanonical Vlasov--Maxwell bracket, based on first order Finite Elements, which is a special case of our framework. The system of ordinary differential equations obtained from the semi-discrete bracket is integrated in time using the splitting method developed by~\citet{Crouseilles:2015} with a correction provided by~\citet{He:2015} (see also~\cite{Qin:2015}). The authors prove the Jacobi identity of the semi-discrete bracket but skip over the Casimir invariants, which also need to be conserved for the semi-discrete system to be Hamiltonian.

In this work, we unify many of the preceding ideas in a general, flexible and rigorous framework based on Finite Element Exterior Calculus (FEEC)~\cite{Arnold.Falk.Winther.2006.anum,Arnold.Falk.Winther.2010.bams,Christiansen:2011,Monk:2003}. We provide a semi-discretization of the noncanonical Vlasov--Maxwell Poisson structure, which preserves the defining properties of the bracket, anti-symmetry and the Jacobi identity, as well as its Casimir invariants, implying that the semi-discrete system is still a Hamiltonian system. Due to the generality of our framework, the aforementioned conservation properties are guaranteed independently of a particular choice of the Finite Element basis, as long as the corresponding Finite Element spaces satisfy certain compatibility conditions. In particular, this includes the spline spaces presented in Section~\ref{sec:feec_splines}. In order to ensure that these properties are also conserved by the fully discrete numerical scheme, the semi-discrete bracket is used in conjunction with Poisson time integrators provided by the previously mentioned splitting method~\cite{Crouseilles:2015, Qin:2015, He:2015} and higher-order compositions thereof.
A semi-discretization of the noncanonical Hamiltonian structure of the relativistic Vlasov--Maxwell system with spin and that for the gyrokinetic Vlasov--Maxwell system have recently been described by~\citet{Burby:2017}.

It is worth emphasizing that the aim and use of preserving the Hamiltonian structure in the course of discretization is not limited to good energy and momentum conservation properties. These are merely by-products but not the goal of the effort.
Furthermore, from a practical point of view, the significance of global energy or momentum conservation by some numerical scheme for some Hamiltonian partial differential equation should not be overestimated. Of course, these are important properties of any Hamiltonian system and should be preserved within suitable error bounds in any numerical simulation.
However, when performing a semi-discretization in space, the resulting finite-dimensional system of ordinary differential equations usually has millions or billions degrees of freedom. 
Conserving only a very small number of invariants hardly restricts the numerical solution of such a large system. It is not difficult to perceive that one can conserve the total energy of a system in a simulation and still obtain false or even unphysical results.
It is much more useful to preserve local conservation laws like the local energy and momentum balance or multi-symplecticity~\cite{Reich:2000, Moore:2003}, thus posing much more severe restrictions on the numerical solution than just conserving the total energy of the system.
A symplectic or Poisson integrator, on the other hand, preserves the whole hierarchy of Poincar{\'e} integral invariants of the finite-dimensional system~\cite{Channell:1990, SanzSerna:1993}. For a Hamiltonian system of ordinary differential equations with $n$ degrees of freedom, e.g., obtained from a semi-discrete Poisson bracket, these are $n$ invariants.
In addition, such integrators often preserve Noether symmetries and the associated local conservation laws as well as Casimir invariants.

We proceed as follows.
In Section~\ref{sec:vlasov-maxwell}, we provide a short review of the Vlasov--Maxwell system and its Poisson bracket formulation, including a discussion of the Jacobi identity, Casimir invariants and invariants commuting with the specific Vlasov--Maxwell Hamiltonian.
In Section~\ref{sec:feec}, we introduce the Finite Element Exterior Calculus framework using the example of Maxwell's equation, we introduce the deRham complex and Finite Element spaces of differential forms.
The actual discretization of the Poisson bracket is performed in Section~\ref{sec:bracket}.
We prove the discrete Jacobi identity and the conservation of discrete Casimir invariants, including the discrete Gauss' law.
In Section~\ref{sec:splitting}, we introduce a splitting for the Vlasov--Maxwell Hamiltonian, which leads to an explicit time stepping scheme. 
Various compositions are used in order to obtain higher order methods.
Backward error analysis is used in order to study the long-time energy behaviour.
In Section~\ref{sec:example}, we apply the method to the Vlasov--Maxwell system in 1d2v using splines for the discretization of the fields. 
Section~\ref{sec:numerics} concludes the paper with numerical experiments, using nonlinear Landau damping and the Weibel instability to verify the favourable properties of our scheme.

\section{The Vlasov--Maxwell System}\label{sec:vlasov-maxwell}

The non-relativistic Vlasov equation for a particle species $s$ of charge $q_s$ and mass $m_s$ reads
\begin{equation}\label{eq:vlasov}
\fracp{\f_s}{t} + \vv \cdot \nabla_\xx \f_s + \frac{q_s}{m_s} (\E + \vv \times \B) \cdot \nabla_\vv \f_s = 0,
\end{equation}
and couples nonlinearly to the Maxwell equations,
\begingroup
\allowdisplaybreaks
\begin{align}
\fracp{\E}{t} - \curl \B &= - \J, \label{eq:ampere} \\
\vphantom{\fracp{}{}}
\fracp{\B}{t}  + \curl \E &= 0, \label{eq:faraday} \\
\vphantom{\fracp{}{}}
\div \E &= \rho, \label{eq:gauss} \\
\vphantom{\fracp{}{}}
\div \B &= 0. \label{eq:gaussm}
\end{align}
\endgroup
These equations are to be solved with suitable initial and boundary conditions.
Here, $\f_{s}$ is the phase space distribution function of particle species $s$, $\E$ and $\B$ are the electric and magnetic fields, respectively, and we have scaled the variables, but retained the mass $m_s$ and the signed charge $q_s$ to distinguish species.
Observe that we use $\grad$, $\curl$, $\div$ to denote $\nabla_\xx$, $\nabla_\xx\times$, $\nabla_\xx\cdot$, respectively, when they act on variables depending only on $\xx$.
The sources for the Maxwell equations, the charge density $\rho$ and the current density $\J$, are obtained from the distribution functions $\f_s$ by
\begin{equation}\label{eq:sources}
\rho = \sum_s q_s \int \f_s \dv, ~~~~ \J = \sum_s q_s \int \f_s \vv \dv.
\end{equation}
Taking the divergence of  Amp\`ere's equation \eqref{eq:ampere} and using Gauss' law \eqref{eq:gauss} gives the continuity equation for charge conservation
\begin{equation}\label{eq:continuity}
\fracp{\rho}{t} + \div \J = 0.
\end{equation}
Equation~\eqref{eq:continuity} serves as a compatibility condition for Maxwell's equations, which are ill posed when~\eqref{eq:continuity} is not satisfied.
Moreover it can be shown that if the divergence constraints \eqref{eq:gauss}-\eqref{eq:gaussm} are satisfied at the initial time, they remain satisfied for all times by the solution of Amp\`ere's equation~\eqref{eq:ampere} and Faraday's law~\eqref{eq:faraday}, which have a unique solution by themselves provided adequate initial and boundary conditions are imposed.
This follows directly from the fact that the divergence of the curl vanishes and Eq.~\eqref{eq:continuity}.
The continuity equation follows from the Vlasov equation by integration over velocity space and using the definitions of charge and current densities.
However this does not necessarily remain true when the charge and current densities are approximated numerically.
The problem for numerical methods is then to find a way to have discrete sources, which satisfy a discrete continuity equation compatible with the discrete divergence and curl operators. 
Another option is to modify the Maxwell equations, so that they are well posed independently of the sources, by introducing two additional scalar unknowns that can be seen as Lagrange multipliers for the divergence constraints.
These should become arbitrarily small when the continuity equation is close to being satisfied.

\subsection{Noncanonical Hamiltonian Structure}

The Vlasov--Maxwell system possesses a noncanonical Hamiltonian structure. The system of equations~\eqref{eq:vlasov}-\eqref{eq:faraday} can be obtained from the following Poisson bracket, a bilinear, anti-symmetric bracket that satisfies Leibniz' rule and the Jacobi identity: 
\begin{align}\label{eq:poisson_bracket_maxwell}
\cbra{\ff}{\fg} [\f_s, \E, \B]
\nonumber
&= \sum \limits_{s} \int \dfrac{\f_s}{m_s} \sbra{ \fracd{\ff}{\f_s} }{ \fracd{\fg}{\f_s} } \dx \dv \\
\nonumber
&+ \sum \limits_{s} \dfrac{q_s}{m_s} \int \f_s \left(
        \nabla_\vv \fracd{\ff}{\f_s} \cdot \fracd{\fg}{\E}
      - \nabla_\vv \fracd{\fg}{\f_s} \cdot \fracd{\ff}{\E}
   \right) \dx \dv \\
\nonumber
&+ \sum \limits_{s} \dfrac{q_s}{m_s^2} \int \f_s \, \B \cdot \left( \nabla_\vv \fracd{\ff}{\f_s} \times \nabla_\vv \fracd{\fg}{\f_s} \right) \dx \dv \\
&+ \int \left( \curl\fracd{\ff}{\E} \cdot  \fracd{\fg}{\B} - \curl\fracd{\fg}{\E} \cdot \fracd{\ff}{\B} \right) \dx ,
\end{align}
where $[f,g]=\nabla_\xx f \cdot  \nabla_\vv g - \nabla_\xx g \cdot  \nabla_\vv f$. 
This bracket was introduced in~\cite{Morrison:1980}, with a term corrected in~\cite{Marsden:1982} (see also~\cite{Weinstein:1981, Morrison:1982}), and its limitation to divergence-free magnetic fields first pointed out in~\cite{Morrison:1982}. See also~\cite{Chandre:2013} and~\cite{Morrison:2013}, where the latter contains the details of the direct proof of the Jacobi identity
\begin{align}
\{ \ff , \{ \fg , \fh \} \} + \{ \fg , \{ \fh , \ff \} \} + \{ \fh , \{ \ff , \fg \} \} = 0 .
\end{align}

The time evolution of any functional $\ff [\f_s, \E, \B]$ is given by
\begin{align}
\dfrac{\dd}{\dd t} \ff [\f_s, \E, \B] = \cbra{\ff}{\ham} ,
\label{eq:eqs_of_motion_time_evolution}
\end{align}
with the Hamiltonian $\ham$ given as the sum of the kinetic energy of the particles and the electric and magnetic field energies,
\begin{align}
\ham = \sum \limits_{s} \dfrac{m_s}{2} \int \abs{\vv}^{2} \, \f_s (t,\xx,\vv) \dx \dv + \dfrac{1}{2} \int \Big( \abs{\E (t, \xx)}^{2} + \abs{\B (t, \xx)}^{2} \Big) \dx .
\label{eq:hamiltonian_vlasov_maxwell}
\end{align}
In order to obtain the Vlasov equations, we consider the functional
\begin{align}
\delta_{\xx\vv} [\f_s]
= \int \f_s (t, \xx', \vv') \, \de (\xx-\xx') \, \de (\vv-\vv') \dx' \dv'
= \f_s (t,\xx,\vv) ,
\label{eq:vlasov_maxwell_dist_functional}
\end{align}
for which the equations of motion \eqref{eq:eqs_of_motion_time_evolution} are computed as
\begin{align}
\fracp{\f_s}{t} (t,\xx,\vv)
\nonumber
&= \int \de (\xx-\xx') \, \de (\vv-\vv') \sbra{ \tfrac{1}{2} \abs{\vv'}^{2} }{ \f_s (t, \xx', \vv') } \dx' \dv' \\
\nonumber
&- \dfrac{q_s}{m_s} \int \de (\xx-\xx') \, \de (\vv-\vv') \left( \nabla_\vv \f_s (t, \xx', \vv') \right) \cdot \E (t, \xx') \dx' \dv' \\
\nonumber
&- \dfrac{q_s}{m_s} \int \de (\xx-\xx') \, \de (\vv-\vv') \left( \nabla_\vv \f_s (t,\xx', \vv') \right) \cdot \left( \B (t, \xx') \times \vv' \right) \dx' \dv' \\
&= - \vv \cdot \nabla_\xx \f_s (t,\xx,\vv) - \dfrac{q_s}{m_s} \big( \E (t, \xx) + \vv \times \B (t, \xx) \big) \cdot \nabla_\vv \f_s (t, \xx,\vv) .
\end{align}
For the electric field, we consider
\begin{align}
\delta_{\xx} [\E] = \int \E (t, \xx') \, \de (\xx-\xx') \dx' = \E (t, \xx) ,
\label{eq:electric_field_functional}
\end{align}
so that from \eqref{eq:eqs_of_motion_time_evolution} we obtain  Amp\`ere's equation,
\begin{align}
\fracp{\E}{t} (t, \xx)
\nonumber
&= \int \bigg( \curl \B (t, \xx') - \sum \limits_{s} q_s \f_s (t, \xx',\vv') \, \vv' \bigg) \, \de (\xx-\xx') \dx' \dv' \\
&= \curl \B (t, \xx) - \J (t, \xx) ,
\end{align}
where the current density $\J$ is given by
\begin{align}
\J (t, \xx) = \sum \limits_{s} q_s \int \f_s (t,\xx,\vv) \, \vv \dv .
\end{align}
And for the magnetic field, we consider
\begin{align}
\delta_{\xx} [\B] = \int \B (t, \xx') \, \de (\xx-\xx') \dx' = \B (t, \xx) ,
\label{eq:magnetic_field_functional}
\end{align}
and obtain the Faraday equation,
\begin{align}
\fracp{\B}{t} (t, \xx)
&= - \int \left( \curl \E (t, \xx') \right) \, \de (\xx-\xx') \dx
 = - \curl \E (t, \xx) .
\end{align}
Our aim is to preserve this noncanonical Hamiltonian structure and its features at the discrete level. This can be done by taking only a finite number of initial positions for the particles instead of a continuum and by taking the electromagnetic fields in finite-dimensional subspaces of the original function spaces. A good candidate for such a discretization is the Finite Element Particle-In-Cell framework. In order to satisfy the continuity equation as well as the identities from vector calculus and thereby preserve Gauss' law and the divergence of the magnetic field, the Finite Element spaces for the different fields cannot be chosen independently. The right framework is given by Finite Element Exterior Calculus (FEEC).

Before describing this framework in more detail, we shortly want to discuss some conservation laws of the Vlasov--Maxwell system. In Hamiltonian systems, there are two kinds of conserved quantities, Casimir invariants and momentum maps.

\subsection{Invariants} \label{sec:casimirs}

A family of conserved quantities are Casimir invariants (Casimirs), which originate from the degeneracy of the Poisson bracket. Casimirs are functionals $\cas (\f_s, \E, \B)$ which Poisson commute with every other functional $\mathcal{G} (\f_s, \E, \B)$, i.e., $\{ \cas , \mathcal{G} \} = 0$.
For the Vlasov--Maxwell bracket, there are several such Casimirs~\cite{Morrison:1987, MorrisonPfirsch:1989, Chandre:2013}. First, the integral of any real function $h_{s}$ of each distribution function $\dist_{s}$ is preserved, i.e.,
\begin{align}
\cas_{s} = \int h_{s} ( \dist_{s} ) \dx \dv .
\end{align}
This family of Casimirs is a manifestation of Liouville's theorem and corresponds to conservation of phase space volume.
Further we have two Casimirs related to Gauss' law~\eqref{eq:gauss} and the divergence-free property of the magnetic field~\eqref{eq:gaussm},
\begin{align}
\cas_{E} &= \int h_{E} (\xx) \big( \div \E - \rho \big) \dx , \\
\cas_{B} &= \int h_{B} (\xx) \, \div \B \dx ,
\end{align}
where $h_{E}$ and $h_{B}$ are arbitrary real functions of $\xx$.
The latter functional, $\cas_{B}$, is not a true Casimir but should rather be referred to as pseudo-Casimir. It acts like a Casimir in that it Poisson commutes with any other functional, but the Jacobi identity is only satisfied when $\div B = 0$ (see \cite{Morrison:1982,Morrison:2013}).

A second family of conserved quantities are momentum maps $\Phi$, which arise from symmetries that preserve the particular Hamiltonian $\ham$, and therefore also the equations of motion. This means that the Hamiltonian is constant along the flow of $\Phi$, i.e.,
\begin{align}\label{eq:momentum_maps_symmetry_condition}
\{ \ham , \Phi \} = 0 .
\end{align}
From Noether's theorem it follows that the generators $\Phi$ of the symmetry are preserved by the time evolution, i.e.,
\begin{align}
\dfrac{\dd \Phi}{\dd t} = 0 .
\end{align}
If the symmetry condition~\eqref{eq:momentum_maps_symmetry_condition} holds, this is obvious by the antisymmetry of the Poisson bracket as
\begin{align}
\dfrac{\dd \Phi}{\dd t} = \{ \Phi , \ham \} = - \{ \ham , \Phi \} = 0 .
\end{align}
Therefore $\Phi$ is a constant of motion if and only if $\{ \Phi , \ham \} = 0$.

The complete set of constants of motion, the algebra of invariants, will be discussed elsewhere. However, as an example of a momentum map we shall consider here the total momentum
\begin{align}\label{eq:momentum}
\mom = \sum_s \int m_s \vv \f_s \dx \dv + \int \E \times \B \dx.
\end{align}
By direct computations, assuming periodic boundary conditions, it can be shown that
\begin{align}\label{dotP}
\dfrac{\dd \mom}{\dd t} = \cbra{\mom}{\ham} = \int \E (\rho - \div \E) \dx = \int \E \, Q(\xx)\dx. 
\end{align}
defining  $Q(\xx):=\rho - \div \E$, which is a local version of the Casimir $\cas_{E}$.
 Therefore, if at $t=0$ the Casimir  $Q \equiv 0$, then  momentum is conserved.   If at $t=0$ the Casimir  $Q\not\equiv 0$, then  momentum is not conserved and it changes in accordance with \eqref{dotP}.  For a multi-species plasma $Q\equiv 0$ is equivalent to the physical requirement that  Poisson's equation be  satisfied.  If for some reason it is not exactly satisfied, then we have violation of momentum conservation. 
 
For a single species plasma, say electrons, with a neutralizing positive background charge $\rho_B(\xx)$, say ions, Poisson's equation is
\begin{align}
\div \E = \rho_B - \rho_e .
\end{align}
The Poisson bracket for this case has the local Casimir
\begin{align}
Q_e = \div \E + \rho_e ,
\end{align}
and it does not recognize the background charge.  Because the background is stationary, the total momentum is
\begin{align}
\mom = \int m_e \vv \, \f_e \dx \dv + \int \E \times \B \dx,
\end{align}
and it satisfies 
\begin{align}\label{eq:dotPe}
\frac{\dd \mom_e}{\dt} = \{ \mom_e, \ham \} = - \int \E \, \rho_B(\xx) \dx .
\end{align}
We will verify this relation in the numerical experiments of Sec.~\ref{sec:numerics_momentum}.

\section{Finite Element Exterior Calculus}\label{sec:feec}

Finite Element Exterior Calculus (FEEC) is a mathematical framework for mixed Finite Element methods, which uses geometrical and topological ideas for systematically analysing the stability and convergence of Finite Element discretizations of partial differential equations.
This proved to be a particularly difficult problem for Maxwell's equation, which we will use in the following as an example for reviewing this framework.

\subsection{Maxwell's Equations}

When Maxwell's equations are used in some material medium, they are best understood by introducing two additional fields. The electromagnetic properties are then defined by the electric and magnetic fields, usually denoted by $\E$ and $\vec{H}$, the displacement field $\vec{D}$ and the magnetic induction $\B$. For simple materials, the electric field is related to the displacement field and the magnetic field to the magnetic induction by
$$\vec{D} = \ve\E, ~~~~   \B = \mu\mathbf{H},$$
where $\ve$ and $\mu$ are the permittivity and permeability tensors reflecting the material properties. In vacuum they become the scalars $\ve_0$ and $\mu_0$, which are unity in our scaled variables, while for more complicated media such as plasmas they can be nonlinear operators \cite{Morrison:2013}.
 The Maxwell equations with the four fields read
\begin{align}
\fracp{\vec{D}}{t} -  \curl \mathbf{H}  &= -\J, \label{eq:ampere_dheb}
\\
\vphantom{\fracp{}{}}
\fracp{\B}{t} + \curl \E &= 0 , \label{eq:faraday_dheb}
\\
\vphantom{\fracp{}{}}
\div \vec{D} &= \rho, \label{eq:gauss_dheb}
\\
\vphantom{\fracp{}{}}
\div \B &= 0. \label{eq:gaussm_dheb}
\end{align}
The mathematical interpretation of these fields become clearer when interpreting them as differential forms:
$\E$ and $\vec{H}$ are 1-forms, $\vec{D}$ and $\B$ are 2-forms. The charge density $\rho$ is a 3-form and the current density $\J$ a 2-form. Moreover, the electrostatic potential $\phi$ is a 0-form and the vector potential $\A$ a 1-form. The $\grad$, $\curl$, $\div$ operators represent the exterior derivative applied respectively to 0-forms, 1-forms and 2-forms.
To be more precise, there are two kinds of differential forms, depending on the orientation. Straight differential forms have an intrinsic (or inner) orientation, whereas twisted differential forms have an outer orientation, defined by the ambient space.
Faraday's equation and $\div \B=0$ are naturally inner oriented, whereas Amp\`ere's equation and Gauss' law are outer oriented. This knowledge can be used to define a natural discretization for Maxwell's equations.
For Finite Difference approximations a dual mesh is needed for the discretization of twisted forms. This can already be found in Yee's scheme~\cite{Yee:1966}.
In the Finite Element context, only one mesh is used, but dual operators are used for the twisted forms.
As an implication, the charge density $\rho$ will be treated as a 0-form and the current density $J$ as a 1-form, instead of a (twisted) 3-form and a (twisted) 2-form, respectively. Another consequence is that Amp\`ere's equation and Gauss' law are being treated weakly while Faraday's equation and $\div \B=0$ are treated strongly.
A detailed description of this formalism can be found, e.g., in Bossavit's lecture notes~\cite{Bossavit:2005}.

\subsection{Finite Element Spaces of Differential Forms}

The full mathematical theory for the Finite Element discretization of differential forms is due to Arnold, Falk and Winther~\cite{Arnold.Falk.Winther.2006.anum,Arnold.Falk.Winther.2010.bams} and is called Finite Element Exterior Calculus (FEEC) (see also~\cite{Monk:2003,Christiansen:2011}). Most Finite Element spaces appearing in this theory were known before, but their connection in the context of differential forms was not made clear.
The first building block of FEEC is the following commuting diagram:
\begin{equation} \label{CGL_df}
    \begin{tikzpicture}[baseline=(current  bounding  box.center)]
  \matrix (m) [matrix of math nodes,row sep=3em,column sep=4em,minimum width=2em] {
           H\Lambda^0(\Omega)
                & H\Lambda^1(\Omega) &H\Lambda^2(\Omega) &H\Lambda^3(\Omega)
    \\
             V_0    & V_1 & V_2 & V_3
    \\
    };
  \path[-stealth]
    (m-1-1) edge node [left] {$\Pi_0$} (m-2-1)
     (m-1-1)       edge node [above] {$\ext$} (m-1-2)
    (m-1-2) edge node [left] {$\Pi_1$} (m-2-2)
    (m-2-1) edge node [above] {$\ext$} (m-2-2)
    (m-1-3) edge node [left] {$\Pi_2$} (m-2-3)
    (m-1-4) edge node [left] {$\Pi_3$} (m-2-4)
     (m-1-2) edge node [above] {$\ext$} (m-1-3)
     (m-2-2) edge node [above] {$\ext$} (m-2-3)
      (m-1-3) edge node [above] {$\ext$} (m-1-4)
     (m-2-3) edge node [above] {$\ext$} (m-2-4)
    ;
\end{tikzpicture}
\end{equation}
where $\Omega \subset \R^3$, $\Lambda^k(\Omega)$ is the space of $k$-forms on $\Omega$ that we endow with the inner product $\langle\alpha,\beta\rangle=\int \alpha \wedge \star \beta$, $\star$ is the Hodge operator and $\ext$ is the exterior derivative that generalizes the gradient, curl and divergence. Then we define 
$$L^2\Lambda^k(\Omega) =\{ \omega \in \Lambda^k(\Omega) \;|\; \langle\omega,\omega \rangle < +\infty \}$$
 and the Sobolev spaces of differential forms 
$$H\Lambda^k(\Omega) = \{ \omega \in L^2\Lambda^k(\Omega) \;|\; \ext\omega \in L^2\Lambda^{k+1}(\Omega) \}.$$
Obviously in a three-dimensional manifold the exterior derivative of a 3-form vanishes so that $H\Lambda^3(\Omega)=L^2(\Omega)$. 
This diagram can also be expressed using the standard vector calculus formalism:
\begin{equation}\label{CGL}
\begin{tikzpicture}[baseline=(current  bounding  box.center)]
  \matrix (m) [matrix of math nodes,row sep=3em,column sep=4em,minimum width=2em] {
           H^1(\Omega) & H(\curl,\Omega) & H(\div,\Omega) & L^2(\Omega)
    \\
           V_0 & V_1 & V_2 & V_3
    \\
    };
  \path[-stealth]
    (m-1-1) edge node [left]  {$\Pi_0$} (m-2-1)
    (m-1-1) edge node [above] {$\grad$} (m-1-2)
    (m-1-2) edge node [left]  {$\Pi_1$} (m-2-2)
    (m-2-1) edge node [above] {$\grad$} (m-2-2)
    (m-1-3) edge node [left]  {$\Pi_2$} (m-2-3)
    (m-1-4) edge node [left]  {$\Pi_3$} (m-2-4)
    (m-1-2) edge node [above] {$\curl$} (m-1-3)
    (m-2-2) edge node [above] {$\curl$} (m-2-3)
    (m-1-3) edge node [above] {$\div$}  (m-1-4)
    (m-2-3) edge node [above] {$\div$}  (m-2-4)
    ;
\end{tikzpicture}
\end{equation}
The first row of~\eqref{CGL} represents the sequence of function spaces involved in Maxwell's equations. Such a sequence is called a complex if at each node, the image of the previous  operator is in the kernel of the next operator, i.e., $\Ima(\grad) \subseteq \Ker(\curl)$ and $\Ima(\curl) \subseteq \Ker(\div)$. The power of the conforming Finite Element framework is that this complex can be reproduced at the discrete level by choosing the appropriate finite-dimensional subspaces $V_0$, $V_1$, $V_2$, $V_3$. The order of the approximation is dictated by the choice made for $V_0$ and the requirement of having a complex at the discrete level. The projection operators $\Pi_i$ are the Finite Element interpolants, which have the property that the diagram is commuting. This means for example, that the $\grad$ of the projection on $V_0$ is identical to the projection of the $\grad$ on $V_1$.
As proven by Arnold, Falk and Winther, their choice of Finite Elements naturally leads to stable discretizations. 

There are many known sequences of Finite Element spaces that fit this diagram. The sequences proposed by Arnold, Falk and Winther are based on well-known Finite Element spaces. On tetrahedra these are $H^1$ conforming $ \mathbb{P}_k$ Lagrange Finite Elements for $V_0$, the $H(\curl)$ conforming N\'ed\'elec elements for $V_1$, the $H(\div)$ conforming Raviart--Thomas elements for $V_2$ and Discontinuous Galerkin elements for $V_3$. A similar sequence can be defined on hexahedra based on the $H^1$ conforming $\mathbb{Q}_k$ Lagrange Finite Elements for $V_0$.
 
Other sequences that satisfy the complex property are available. Let us in particular cite the mimetic spectral elements~\cite{kreeft2011, Gerritsma:2012, Palha:2014} and the spline Finite Elements~\cite{buffa2010isogeometric, Buffa:2011, Ratnani:2012} that we shall use in this work, as splines are generally favoured in PIC codes due to their smoothness properties that enable noise reduction.

\subsection{Finite Element Discretization of Maxwell's Equations}

This framework is enough to express discrete relations between all the straight (or primal forms), i.e.,
$\E$, $\B$, $\A$ and $\phi$.
 The commuting diagram yields a direct expression of the discrete Faraday equation. Indeed projecting all the components of the equation onto $V_2$ yields
\begin{align*}
\fracp{\Pi_2 \B}{t} + \Pi_2 \curl \E =0,
\end{align*}
which is equivalent, due to the commuting diagram property, to 
\begin{align*}
\fracp{\Pi_2 \B}{t} + \curl \Pi_1\E =0 .
\end{align*}
Denoting with an $h$ index the discrete fields, $\B_h=\Pi_2 \B$, $\E_h=\Pi_1 \E$, this yields the discrete Faraday equation,
\begin{equation}\label{eq:fardisc}
\fracp{\B_h}{t} + \curl \E_h=0.
\end{equation}
In the same way, the discrete electric and magnetic fields are defined exactly as in the continuous case from the discrete potentials, thanks to the compatible Finite Element spaces,
\begin{align}\label{eq:potentialdisc}
\E_h
&= \Pi_{1} \E
 = - \Pi_{1} \grad \phi - \Pi_{1} \fracp{\A}{t}
 = - \grad \Pi_{0} \phi - \fracp{\Pi_{1} \A}{t}
 = - \grad \phi_h - \fracp{\A_h}{t} ,
\\
\B_h
&= \Pi_{2} \B
 = \Pi_{2} \curl \A
 = \curl \Pi_{1} \A
 = \curl \A_h ,
\end{align}
so that automatically we get
\begin{equation}\label{eq:gaussmdisc}
\div \B_h = 0 .
\end{equation}

On the other hand, Amp\`ere's equation and Gauss' law relate expressions involving twisted differential forms. In the Finite Element framework, these should be expressed on the dual complex to~\eqref{CGL}.
But due to the property that the dual of an operator in $L^2(\Omega)$ can be identified with its $L^2$ adjoint via an inner product, the discrete dual spaces are such that $V_0^*=V_3$, $V_1^*=V_2$, $V_2^*=V_1$ and $V_3^*=V_0$, so that the dual operators and spaces are not explicitly needed.
They are most naturally used seamlessly by keeping the weak formulation of the corresponding equations.
The weak form of Amp\`ere's equation is found by taking the dot product of \eqref{eq:ampere} with a test function $\bar{\E} \in H(\curl,\Omega)$ and applying a Green identity.
Assuming periodic  boundary conditions, the weak solution of Amp\`ere's equation $ (\E, \B) \in H(\curl,\Omega) \times H(\div,\Omega)$
is characterized by
\begin{equation}\label{eq:varampere}
\frac{\dd}{\dd t} \int_{\Omega} \E \cdot \bar{\E} \dx - \int_{\Omega} \B \cdot
\curl \bar{\E} \dx = - \int_{\Omega} \J \cdot \bar{\E} \dx
\quad \forall \bar{\E} \in H(\curl,\Omega).
\end{equation}
The discrete version is obtained by replacing the continuous spaces by their finite-dimensional subspaces.  The approximate solution $ (\E_h,\B_h) \in V_1 \times V_2$ is characterized by
\begin{equation}\label{eq:varamperedisc}
\frac{\dd}{\dd t} \int_{\Omega} \E_h \cdot \bar{\E}_h \dx - \int_{\Omega} \B_h \cdot \curl \bar{\E}_h \dx = - \int_{\Omega} \J_h \cdot \bar{\E}_h \dx
\quad \forall \bar{\E}_h \in V_1.
\end{equation}
In the same way the weak solution of Gauss' law with $\E \in H(\curl,\Omega)$ is characterized by
\begin{align}\label{eq:vargauss}
\int_{\Omega} \E \cdot \nabla \bar{\phi} \dx = -  \int_{\Omega} \rho \bar{\phi} \dx \quad \forall \bar{\phi} \in H^1(\Omega),
\end{align}
its discrete version for $\E_h \in V_1$ being characterized by
\begin{equation}\label{eq:vargaussdisc}
\int_{\Omega} \E_h \cdot \nabla \bar{\phi}_h \dx = - \int_{\Omega}
\rho_h \bar{\phi}_h \dx \quad \forall \bar{\phi}_h \in V_0.
\end{equation}
The last step for the Finite Element discretization is to define a basis for each of the finite-dimensional spaces $V_0,V_1,V_2,V_3$, with $\dim V_k=N_k$ and to find equations relating the coefficients on these bases.
Let us denote by $\{\Lambda_i^0\}_{i=1 \dots N_0}$ and $\{\Lambda_i^3\}_{i=1 \dots N_3}$ a basis of $V_0$ and $V_3$, respectively, which are spaces of scalar functions, and $\{\Lab_{i,\mu}^1\}_{i=1 \dots N_1, \mu=1 \dots 3}$ a basis of $V_1\subset H(\curl,\Omega)$ and $\{\Lab_{i,\mu}^2\}_{i=1 \dots N_2, \mu=1 \dots 3}$ a basis of $V_2\subset H(\div,\Omega)$, which are vector valued functions,
\begin{align}
\Lab_{i,1}^k &= \begin{pmatrix} \Lambda_i^{k,1} \\ 0 \\ 0 \end{pmatrix} , &
\Lab_{i,2}^k &= \begin{pmatrix} 0 \\ \Lambda_i^{k,2} \\ 0 \end{pmatrix} , &
\Lab_{i,1}^k &= \begin{pmatrix} 0 \\ 0 \\ \Lambda_i^{k,3} \end{pmatrix} , &
k &= 1,2 .
\end{align}
Let us note that the restriction to a basis of this form is not strictly necessary and the generalization to more general bases is straightforward. However, for didactical reasons we stick to this form of the basis as it simplifies some of the computations and thus helps to clarify the following derivations.

In order to keep a concise notation, and by slight abuse of the same, we introduce vectors of basis functions $\Lab^k = (\Lab_{1,1}^{k}, \Lab_{1,2}^{k}, \dots, \Lab_{N_k,3}^{k})^\top$ for $k=1,2$, which are indexed by $I = 3(i-1)+\mu = 1...3N_{k}$ with $i = 1 \dots N_{k}$ and $\mu = 1 \dots 3$, and $\Lab^k = (\La_1^{k}, \La_2^{k}, \dots, \La_{N_k}^{k})^\top$ for $k=0,3$, which are indexed by $i = 1 \dots N_{k}$.

We shall also need for each basis the dual basis, which in Finite Element terminology corresponds to the degrees of freedom.
For each basis $\La_i^k$ for $k=0,3$ and $\Lab_I^k$ for $k=1,2$, the dual basis is denoted by $\Si_i^k$ and $\Sib_I^k$, respectively, and defined by
\begin{align}
\pair{ \Si^{k}_{i} }{ \La^{k}_{j} }
&= \int \Si^{k}_{i} (\xx) \, \La^{k}_{j} (\xx) \dx
 = \delta_{ij} , &
k &= 0,3 ,
\end{align}
for the scalar valued bases $\La_i^k$, and
\begin{align}
\pair{ \Sib^{k}_{I} }{ \Lab^{k}_{J} }
&= \int \Sib^{k}_{I} (\xx) \cdot \Lab^{k}_{J} (\xx) \dx
 = \delta_{IJ} , &
k &= 1,2 ,
\end{align}
for the vector valued bases $\Lab_{I}^k$, where $\pair{\cdot}{\cdot}$ denotes the $L^{2}$ inner product in the appropriate space and $\delta_{IJ}$ is the Kronecker symbol, whose value is unity for $I=J$ and zero otherwise.
We introduce the linear functionals $L^2 \Lambda^k (\Omega) \rightarrow \rsp$, which are denoted by $\sigma_i^k$ for $k=0,3$ and by $\sib_I^k$ for $k=1,2$, respectively.
On the finite element space they are represented by the dual basis functions $\Sib_I^{k}$ and defined by
\begin{align}\label{eq:dual_basis_03}
\si_i^k (\omega) &= \pair{ \Si^{k}_{i} }{ \omega }
\quad \forall \omega \in L^2 \Lambda^k (\Omega) , &
k &= 0,3 ,
\end{align}
and
\begin{align}\label{eq:dual_basis_12}
\sib_I^k (\omb) &= \pair{ \Sib^{k}_{I} }{ \omb }
\quad \forall \omb \in L^2 \Lambda^k (\Omega) , &
k &= 1,2 ,
\end{align}
so that $\sigma_i^k(\Lambda_j^k)=\delta_{ij}$ and $\sib_I^k(\Lab_J^k)=\delta_{IJ}$ for the appropriate $k$.
Elements of the finite-dimensional spaces can be expanded on their respective bases, e.g, elements of $V_1$ and $V_2$, respectively, as
\begin{align}\label{eq:fields_discrete}
\E_h (t,\xx) &= \sum_{i=1}^{N_1} \sum_{\mu=1}^{3} e_{i,\mu} (t) \, \Lab_{i,\mu}^1 (\xx), &
\B_h (t,\xx) &= \sum_{i=1}^{N_2} \sum_{\mu=1}^{3} b_{i,\mu} (t) \, \Lab_{i,\mu}^2 (\xx),
\end{align}
denoting by $\ce = (e_{1,1}, \, e_{1,2}, \, \dots, \, e_{N_1,3})^\top \in \rsp^{3N_1}$ and $\cb = (b_{1,1}, \, b_{1,2}, \, \dots, \, b_{N_2,3})^\top \in \rsp^{3N_2}$ the corresponding degrees of freedom with $e_{i,\mu} = \sib^1_{i,\mu} (\E_h)$ and $b_{i,\mu} = \sib^2_{i,\mu} (\B_h)$, respectively.
Denoting the elements of $\ce$ by $\ce_{I}$ and the elements of $\cb$ by $\cb_{I}$, we have that $\ce_{I} = \sib^1_{I} (\E_h)$ for $I = 1 \hdots 3N_{1}$ and $\cb_{I} = \sib^2_{I} (\B_h)$ for $I = 1 \hdots 3N_{2}$, respectively, and can re-express~\eqref{eq:fields_discrete} as
\begin{align}\label{eq:fields_discrete2}
\E_h (t,\xx) &= \sum_{I=1}^{3N_1} \ce_{I} (t) \, \Lab_{I}^1 (\xx) , &
\B_h (t,\xx) &= \sum_{I=1}^{3N_2} \cb_{I} (t) \, \Lab_{I}^2 (\xx) .
\end{align}
Henceforth we will use both notations in parallel, choosing whichever is more practical at any given time.

Due to the complex property we have that $\curl \E_h \in V_2$ for all  $\E_h \in V_1$, 
so that $\curl \E_h$ can be expressed in the basis of $V_2$ by
\begin{align*}
\curl \E_h = \sum_{i=1}^{N_2} \sum_{\mu=1}^{3} c_{i,\mu} \Lab_{i,\mu}^2 .
\end{align*}
Let us also denote by $\cc=(c_{1,1}, \, c_{1,2}, \, \dots, \, c_{N_2,3})^\top$, so that $\curl \E_h$ can also be written as
\begin{align*}
\curl \E_h = \sum_{I=1}^{3 N_2} \cc_{I} \Lab_{I}^2 .
\end{align*}
On the other hand
\begin{align*}
\curl \E_h
&= \curl \left( \sum_{I=1}^{3N_1} \ce_{I} \, \Lab_{I}^1 \right)
 = \sum_{I=1}^{3N_1} \ce_{I} \, \curl \Lab_{I}^1 , \\
\sib_{I}^2 (\curl \E_h)
&= \sum_{J=1}^{3N_1} e_{J} \, \sib_{I}^2 (\curl \Lab_{J}^1) .
\end{align*}
Denoting by $\C$ the discrete curl matrix, 
\begin{align}\label{eq:rotation_matrix}
\C = (\sib_I^2 (\curl\Lab_J^1))_{1 \leq I \leq 3N_2, \, 1 \leq J \leq 3N_1} ,
\end{align}
the degrees of freedom of $\curl \E_h$ in $V_2$ are related to the degrees of freedom of $\E_h$ in $V_1$ by $\cc = \C \ce$. In the same way we can define the discrete
gradient matrix $\G$ and the discrete divergence matrix $\D$, given by
\begin{align}
\G &= (\sib_I^1 (\grad\Lab_J^0))_{1 \leq I \leq 3N_1, \, 1 \leq J \leq N_0} &
& \text{and} &
\D &= (\sib_I^3 (\div\Lab_j^2))_{1 \leq I \leq N_3, \, 1 \leq J \leq 3N_2} ,
\end{align}
respectively.
Denoting by $\cphi = (\phy_1, \dots, \phy_{N_0})^\top$ and $\ca = (\vara_{1,1}, \, \vara_{1,2}, \, \dots ,\, \vara_{N_1,3})^\top$ the degrees of freedom of the potentials $\phi_{h}$ and $\A_{h}$, with $\varphi_i = \sigma^0_i (\phi_h)$ for $1 \leq i \leq N_0$ and $\ca_I = \sib^1_I (\A_h)$ for $1 \leq I \leq 3N_1$, the relation~\eqref{eq:potentialdisc} between the discrete fields~\eqref{eq:fields_discrete} and the potentials can be written using only the degrees of freedom as
\begin{align}
\ce = - \G \cphi - \frac{\dd \ca}{\dd t}, ~~~~~~
\cb =   \C \ca.
\end{align}
Finally, we need to define the so-called mass matrices in each of the discrete spaces $V_i$, which define the discrete Hodge operator linking the primal complex with the dual complex.
We denote by $(\MM_{0})_{ij}= \int_\Omega \Lambda_i^0 (\xx) \, \Lambda_j^0 (\xx) \dx$ with $1 \leq i,j \leq N_{0}$ and $(\MM_{1})_{IJ}= \int_\Omega \Lab_I^1 (\xx) \cdot \Lab_J^1 (\xx) \dx$ with $1 \leq I, J \leq 3N_1$ the mass matrices in $V_0$ and $V_1$, respectively, and similarly $\MM_{2}$ and $\MM_{3}$ the mass matrices in $V_2$ and $V_3$.
Using these definitions as well as $\crho = (\varrho_1, \, \dots, \, \varrho_{N_0})^\top$ and $\cj = (j_{1,1}, \, j_{1,2}, \, \dots, \, j_{N_1,3})^\top$ with $\varrho_i = \sigma^0_i(\rho_h)$ for $1 \leq i \leq N_0$ and $\cj_I = \sib^1_I (\J_h)$ for $1 \leq I \leq 3N_1$ (recall that the charge density $\rho$ is treated as a 0-form and the current density $J$ as a 1-form), we obtain a system of ordinary differential equations for each of the continuous equations, namely
\begin{align}
\MM_{1} \frac{\dd \ce}{\dd t} - \C^\top \MM_{2} \cb &= -\cj \label{systfe2damp} , \\
\frac{\dd \cb}{\dd t} +\C \ce &= 0 \label{systfe2dfar} , \\
\G^\top \MM_{1} \ce &= \crho \label{systfe2dgauss} , \\
\D \cb &= 0 \label{systfe2dgaussm} .
\end{align}
It is worth emphasizing that $\div B = 0$ is satisfied in strong form, which is important for the Jacobi identity of the discretized Poisson bracket (cf.~Section~\ref{sec:jacobi_identity}).
The complex properties can also be expressed at the matrix level. The primal sequence being
\begin{equation}\label{eq:exseqmath_primal}
\begin{tikzcd}[column sep = 3em] 
\R^{N_0} \arrow{r}{\G} &
\R^{N_1} \arrow{r}{\C} &
\R^{N_2} \arrow{r}{\D} &
\R^{N_3} ,
\end{tikzcd}
\end{equation}
with $\Ima \G \subseteq \Ker \C$, $\Ima \C \subseteq \Ker \D$, and the dual sequence being
\begin{equation}\label{eq:exseqmath_dual}
\begin{tikzcd}[column sep = 3em] 
\R^{N_3} \arrow{r}{\D^\top} &
\R^{N_2} \arrow{r}{\C^\top} &
\R^{N_1} \arrow{r}{\G^\top} &
\R^{N_0} ,
\end{tikzcd}
\end{equation}
with $\Ima \D^\top \subseteq \Ker \C^\top$, $\Ima \C^\top \subseteq \Ker \G^\top$.

\subsection{Example: B-Spline Finite Elements}
\label{sec:feec_splines}

In the following, we will use so-called basic splines or B-splines as bases for the Finite Element function spaces.
B-splines are piecewise polynomials. The points where two polynomials connect are called knots. 
The $j$-th basic spline (B-spline) of degree $p$ can be defined recursively by
\begin{align}\label{eq:spline_definition_recursive}
N_{j}^{p} (x) = w_{j}^{p} (x) \, N_{j}^{p-1} (x) + (1 - w_{j+1}^{p} (x) ) \, N_{j+1}^{p-1} (x) ,
\end{align}
where
\begin{align}\label{eq:spline_definition_weights}
w_{j}^{p} (x) = \dfrac{x - x_{j}}{x_{j+p} - x_{j}} ,
\end{align}
and
\begin{align}\label{eq:spline_definition_zero}
N_{j}^{0} (x) =
\begin{cases}
1 & x \in [ x_{j} , x_{j+1} ) , \\
0 & \text{else} ,
\end{cases}
\end{align}
with the knot vector $\Xi = \{ x_{i} \}_{1 \leq i \leq N+k}$ being a non-decreasing sequence of points.
The knot vector can also contain repeated knots. If a knot $x_{i}$ has multiplicity $m$, then the B-spline is $C^{p-m}$ continuous at $x_{i}$.
The derivative of a B-spline of degree $p$ can easily be computed as the difference of two B-splines of degree $p-1$,
\begin{align}\label{eq:spline_derivative_recursive}
\frac{\dd}{\dd x} N_{j}^{p} (x) = p \, \bigg( \dfrac{N_{j}^{p-1} (x)}{x_{j+p} - x_{j}} - \dfrac{N_{j+1}^{p-1} (x)}{x_{j+p+1} - x_{j+1}} \bigg) .
\end{align}
For convenience, we introduce the following shorthand notation for differentials,
\begin{align}\label{eq:spline_derivative_shorthand}
D_{j}^{p} (x) = p \, \dfrac{N_{j}^{p-1} (x)}{x_{j+p} - x_{j}} .
\end{align}
In the case of an equidistant grid with grid step size $\Delta x = x_{j+1} - x_{j}$, this simplifies to
\begin{align}
D_{j}^{p} (x) = \dfrac{N_{j}^{p-1} (x)}{\Delta x} .
\end{align}
Using $D_{j}^{p}$ the recursion formula~\eqref{eq:spline_definition_recursive} becomes
\begin{align}
\frac{\dd}{\dd x} N_{j}^{p} (x) = D_{j}^{p} (x) - D_{j+1}^{p} (x) .
\end{align}
In more than one dimension, we can define tensor-product B-spline basis functions, e.g., for three dimensions as
\begin{align}
N^{p}_{ijk} = N^{p}_{i} (x_{1}) \otimes N^{p}_{j} (x_{2}) \otimes N^{p}_{k} (x_{3}) .
\end{align}
The bases of the differential form spaces will be tensor products of the basis functions $N_{i}^{p}$ and the differentials $D_{j}^{p}$.
The discrete function spaces are given by
\begin{subequations}\label{eq:spline_basis}
\begin{align}
\Lambda_{h}^{0} (\Omega) = \spn \bigg\{ &
N_{i}^{p} (x_{1}) \, N_{j}^{p} (x_{2}) \, N_{k}^{p} (x_{3}) \;
\Big\vert \; 1 \leq i \leq N_{1} , \, 1 \leq j \leq N_{2} , \, 1 \leq k \leq N_{3}
\bigg\} , \\
\nonumber
\Lambda_{h}^{1} (\Omega) = \spn \Bigg\{ &
\left( \begin{matrix}
D_{i}^{p} (x_{1}) \, N_{j}^{p} (x_{2}) \, N_{k}^{p} (x_{3}) \\
0 \\
0 
\end{matrix} \right) , \,
\left( \begin{matrix}
0 \\
N_{i}^{p} (x_{1}) \, D_{j}^{p} (x_{2}) \, N_{k}^{p} (x_{3}) \\
0
\end{matrix} \right) , \\
&
\left( \begin{matrix}
0 \\
0 \\
N_{i}^{p} (x_{1}) \, N_{j}^{p} (x_{2}) \, D_{k}^{p} (x_{3})
\end{matrix} \right)
\bigg\vert \; 1 \leq i \leq N_{1} , \, 1 \leq j \leq N_{2} , \, 1 \leq k \leq N_{3}
\Bigg\} , \\
\nonumber
\Lambda_{h}^{2} (\Omega) = \spn \Bigg\{ &
\left( \begin{matrix}
N_{i}^{p} (x_{1}) \, D_{j}^{p} (x_{2}) \, D_{k}^{p} (x_{3}) \\
0 \\
0 
\end{matrix} \right) , \,
\left( \begin{matrix}
0 \\
D_{i}^{p} (x_{1}) \, N_{j}^{p} (x_{2}) \, D_{k}^{p} (x_{3}) \\
0 
\end{matrix} \right) , \\
&
\left( \begin{matrix}
0 \\
0 \\
D_{i}^{p} (x_{1}) \, D_{j}^{p} (x_{2}) \, N_{k}^{p} (x_{3})
\end{matrix} \right)
\bigg\vert \; 1 \leq i \leq N_{1} , \, 1 \leq j \leq N_{2} , \, 1 \leq k \leq N_{3}
\Bigg\} , \\
\Lambda_{h}^{3} (\Omega) = \spn \bigg\{ &
D_{i}^{p} (x_{1}) \, D_{j}^{p} (x_{2}) \, D_{k}^{p} (x_{3}) \;
\Big\vert \; 1 \leq i \leq N_{1} , \, 1 \leq j \leq N_{2} , \, 1 \leq k \leq N_{3}
\bigg\} .
\end{align}
\end{subequations}
These choices appear quite natural when one considers the action of the gradient on 0-forms, the curl on 1-forms and the divergence on 2-forms.
In the following, we will exemplify this using the potentials and fields of Maxwell's equations.
The semi-discrete potentials are written in the respective spline basis (\ref{eq:spline_basis}) as
\begin{align}
\phi_{h} (t,\xx) &= \sum \limits_{i,j,k} \phy_{ijk} (t) \, N_{i}^{p} (x_{1}) \, N_{j}^{p} (x_{2}) \, N_{k}^{p} (x_{3}) ,
\end{align}
and
\begin{align}
\A_{h}    (t,\xx) &= \sum \limits_{i,j,k} \left( \begin{matrix}
a^{1}_{ijk} (t) \, D_{i}^{p} (x_{1}) \, N_{j}^{p} (x_{2}) \, N_{k}^{p} (x_{3}) \\
a^{2}_{ijk} (t) \, N_{i}^{p} (x_{1}) \, D_{j}^{p} (x_{2}) \, N_{k}^{p} (x_{3}) \\
a^{3}_{ijk} (t) \, N_{i}^{p} (x_{1}) \, N_{j}^{p} (x_{2}) \, D_{k}^{p} (x_{3}) 
\end{matrix} \right) .
\end{align}
Computing the gradient of the semi-discrete scalar potential $\phi_{h}$, we find
\begin{align}
\grad \phi_{h}
\nonumber
&= \sum \limits_{i,j,k} \phy_{ijk} (t) \, \grad \big[ N_{i}^{p} (x_{1}) \, N_{j}^{p} (x_{2}) \, N_{k}^{p} (x_{3}) \big] , \\
&= \sum \limits_{i,j,k} \phy_{ijk} (t) \, \begin{pmatrix}
\big[ D_{i}^{p} (x_{1}) - D_{i+1}^{p} (x_{1}) \big] \, N_{j}^{p} (x_{2}) \, N_{k}^{p} (x_{3}) \\
N_{i}^{p} (x_{1}) \, \big[ D_{j}^{p} (x_{2}) - D_{j+1}^{p} (x_{2}) \big] \, N_{k}^{p} (x_{3}) \\
N_{i}^{p} (x_{1}) \, N_{j}^{p} (x_{2}) \, \big[ D_{k}^{p} (x_{3}) - D_{k+1}^{p} (x_{3}) \big] \\
\end{pmatrix} .
\end{align}
Assuming periodic boundary conditions, the sums can be re-arranged to give
\begin{align}
\grad \phi_{h}
&= \sum \limits_{i,j,k} \begin{pmatrix}
\big[ \phy_{ijk} (t) - \phy_{i-1jk} (t) \big] \, D_{i}^{p} (x_{1}) \, N_{j}^{p} (x_{2}) \, N_{k}^{p} (x_{3}) \\
\big[ \phy_{ijk} (t) - \phy_{ij-1k} (t) \big] \, N_{i}^{p} (x_{1}) \, D_{j}^{p} (x_{2}) \, N_{k}^{p} (x_{3}) \\
\big[ \phy_{ijk} (t) - \phy_{ijk-1} (t) \big] \, N_{i}^{p} (x_{1}) \, N_{j}^{p} (x_{2}) \, D_{k}^{p} (x_{3}) \\
\end{pmatrix} .
\end{align}
Similarly the curl of the semi-discrete vector potential $\A_{h}$ is computed as
\begin{align}
\curl \A_{h}
\nonumber
&= \sum \limits_{i,j,k} \begin{pmatrix}
a^{3}_{ijk} (t) \, N_{i}^{p} (x_{1}) \, \big[ D_{j}^{p} (x_{2}) - D_{j+1}^{p} (x_{2}) \big] \, D_{k}^{p} (x_{3}) \\
a^{1}_{ijk} (t) \, D_{i}^{p} (x_{1}) \, N_{j}^{p} (x_{2}) \, \big[ D_{k}^{p} (x_{3}) - D_{k+1}^{p} (x_{3}) \big] \\
a^{2}_{ijk} (t) \, \big[ D_{i}^{p} (x_{1}) - D_{i+1}^{p} (x_{1}) \big] \, D_{j}^{p} (x_{2}) \, N_{k}^{p} (x_{3}) \\
\end{pmatrix} \\
&- \sum \limits_{i,j,k} \begin{pmatrix}
a^{2}_{ijk} (t) \, N_{i}^{p} (x_{1}) \, D_{j}^{p} (x_{2}) \, \big[ D_{k}^{p} (x_{3}) - D_{k+1}^{p} (x_{3}) \big] \\
a^{3}_{ijk} (t) \, \big[ D_{i}^{p} (x_{1}) - D_{i+1}^{p} (x_{1}) \big] \, N_{j}^{p} (x_{2}) \, D_{k}^{p} (x_{3}) \\
a^{1}_{ijk} (t) \, D_{i}^{p} (x_{1}) \, \big[ D_{j}^{p} (x_{2}) - D_{j+1}^{p} (x_{2}) \big] \, N_{k}^{p} (x_{3}) \\
\end{pmatrix} .
\end{align}
Again, assuming periodic boundary conditions, the the sums can be re-arranged to give
\begin{align}
\curl \A_{h} &= \sum \limits_{i,j,k} \begin{pmatrix}
\Big( \big[ a^{3}_{ijk} (t) - a^{3}_{ij-1k} (t) \big] - \big[ a^{2}_{ijk} (t) - a^{2}_{ijk-1} (t) \big] \Big) \, N_{i}^{p} (x_{1}) \, D_{j}^{p} (x_{2}) \, D_{k}^{p} (x_{3}) \\
\Big( \big[ a^{1}_{ijk} (t) - a^{1}_{ijk-1} (t) \big] - \big[ a^{3}_{ijk} (t) - a^{3}_{i-1jk} (t) \big] \Big) \, D_{i}^{p} (x_{1}) \, N_{j}^{p} (x_{2}) \, D_{k}^{p} (x_{3}) \\
\Big( \big[ a^{2}_{ijk} (t) - a^{2}_{i-1jk} (t) \big] - \big[ a^{1}_{ijk} (t) - a^{1}_{ij-1k} (t) \big] \Big) \, D_{i}^{p} (x_{1}) \, D_{j}^{p} (x_{2}) \, N_{k}^{p} (x_{3}) \\
\end{pmatrix} .
\end{align}
Given the above, we determine the semi-discrete electromagnetic fields $\E_{h}$ and $\B_{h}$.
The electric field $\E_{h} = - \grad \phi_{h} - \partial \A_{h} /\partial t$ is computed as
\begin{align}
\E_{h} (t, \xx)
&= \sum \limits_{i,j,k} \begin{pmatrix}
\big[ \phy_{i-1jk} (t) - \phy_{ijk} (t) - \dot{a}_{ijk}^{1} (t) \big] \, D_{i}^{p} (x_{1}) \, N_{j}^{p} (x_{2}) \, N_{k}^{p} (x_{3}) \\
\big[ \phy_{ij-1k} (t) - \phy_{ijk} (t) - \dot{a}_{ijk}^{2} (t) \big] \, N_{i}^{p} (x_{1}) \, D_{j}^{p} (x_{2}) \, N_{k}^{p} (x_{3}) \\
\big[ \phy_{ijk-1} (t) - \phy_{ijk} (t) - \dot{a}_{ijk}^{3} (t) \big] \, N_{i}^{p} (x_{1}) \, N_{j}^{p} (x_{2}) \, D_{k}^{p} (x_{3}) 
\end{pmatrix} ,
\end{align}
and the magnetic field $B_{h} = \curl A_{h}$ as
\begin{align}
\B_{h} (t, \xx)
&= \sum \limits_{i,j,k} \begin{pmatrix}
\big[ \big( a^{3}_{ijk} (t) - a^{3}_{ij-1k} (t) \big) - \big( a^{2}_{ijk} (t) - a^{2}_{ijk-1} (t) \big) \big] \, N_{i}^{p} (x_{1}) \, D_{j}^{p} (x_{2}) \, D_{k}^{p} (x_{3}) \\
\big[ \big( a^{1}_{ijk} (t) - a^{1}_{ijk-1} (t) \big) - \big( a^{3}_{ijk} (t) - a^{3}_{i-1jk} (t) \big) \big] \, D_{i}^{p} (x_{1}) \, N_{j}^{p} (x_{2}) \, D_{k}^{p} (x_{3}) \\
\big[ \big( a^{2}_{ijk} (t) - a^{2}_{i-1jk} (t) \big) - \big( a^{1}_{ijk} (t) - a^{1}_{ij-1k} (t) \big) \big] \, D_{i}^{p} (x_{1}) \, D_{j}^{p} (x_{2}) \, N_{k}^{p} (x_{3})
\end{pmatrix} .
\end{align}
Now it becomes clear why definitions~\eqref{eq:spline_basis} are the natural choice for the spline bases and it is straightforward to verify that $\div \B_{h} = 0$.

\section{Discretization of the Hamiltonian Structure}\label{sec:bracket}

The continuous bracket~\eqref{eq:poisson_bracket_maxwell} is for the Eulerian (as opposed to Lagrangian) formulation of the Vlasov equation, and operates on functionals of the distribution function and the electric and magnetic fields. 
We incorporate a discretization that uses a finite number of characteristics instead of the continuum particle distribution function. This is done by localizing the distribution function on particles, which amounts to a Monte Carlo discretization of the first three integrals in~\eqref{eq:poisson_bracket_maxwell} if the initial phase space positions are randomly drawn.
Moreover instead of allowing the fields $\E$ and $\B$ to vary in $H(\curl,\Omega)$ and $H(\div,\Omega)$, respectively, we keep them in the discrete subspaces $V_1$ and $V_2$.
This procedure yields a discrete Poisson bracket, from which one obtains the dynamics of a large but finite number of degrees of freedom: the particle phase space positions $\zz_a = (\xa,\va)$, where $a = 1, ..., N_{p}$, with $N_{p}$ the number of particles, and the coefficients of the fields in the Finite Element basis, where we denote by $\ce_I$ the degrees of freedom for $\E_h$ and by $\cb_I$ the degrees of freedom for $\B_h$, as introduced in Sec.~\ref{sec:feec}. 
The FEEC framework of Sec.~\ref{sec:feec} automatically provides the following discretization spaces for the potentials, the fields and the densities:
\begin{align*}
\phi_h, \rho_h \in V_0 , ~~~~
\A_h, \E_h, \J_h \in V_1 , ~~~~
\B_h \in V_2 .
\end{align*}
Recall that the coefficient vectors of the fields are denoted $\ce$ and $\cb$.
In order to also get a vector expression for the particle quantities, we denote by
\begin{align}
\X = ( \xx_1, \dots, \xx_{N_p} )^\top ,  ~~~~
\V = ( \vv_1, \dots, \vv_{N_p} )^\top .
\end{align}
We use this setting to transform~\eqref{eq:poisson_bracket_maxwell} into a discrete Poisson bracket for the dynamics of the coefficients $\ce$, $\cb$, $\X$ and $\V$.

\subsection{Discretization of the Functional Field Derivatives}

Upon inserting \eqref{eq:fields_discrete}, any functional $\ff [\E_{h}]$ can be considered as a function $\dff (\ce)$ of the Finite Element coefficients,
\begin{align}
\ff [\E_{h}] = \dff (\ce) .
\end{align}
Therefore, we can write the first variation of $\ff[\E]$,
\begin{align}
\delta \ff [\E] = \pair{ \fracd{\ff [\E]}{\E} }{ \de \E } ,
\end{align}
as
\begin{align}\label{eq:functional_derivative_fields_1}
  \pair{ \fracd{\ff [\E_{h}]}{\E} }{ \bar{\E}_{h} }_{L^{2}}
= \pair{ \fracp{\dff (\ce)}{\ce} }{ \bar{\ce} }_{\R^{N_{1}}} ,
\end{align}
with
\begin{align}
\bar{\E}_{h} (t, \xx) &= \sum_{I=1}^{3N_1} \bar{\ce}_{I} (t) \, \Lab_{I}^1 (\xx) , &
\bar{\ce} &= (\bar{e}_{1,1}, \, \bar{e}_{1,2}, \, \dots, \, \bar{e}_{N_1,3})^\top \in \rsp^{3N_1} .
\end{align}
Let $\Sib^{1}_{I} (\xx)$ denote the dual basis of $ \Lab^{1}_{I} (\xx)$ with respect to the $L^{2}$ inner product~\eqref{eq:dual_basis_12} and let us express the functional derivative on this dual basis, i.e.,
\begin{align}
\fracd{\ff [\E_{h}]}{\E} &= \sum_{I=1}^{3N_1} \cae_{I} (t) \, \Sib^{1}_{I} (\xx) , &
\cae &= (\varae_{1,1}, \, \varae_{1,2}, \, \dots, \, \varae_{N_1,3})^\top \in \rsp^{3N_1} .
\end{align}
Using \eqref{eq:functional_derivative_fields_1} for $\bar{\ce} = (0, \, \dots, \, 0, \, 1, \, 0, \, \dots, \, 0)^\top$ with $1$ at the $I$-th position and $0$ everywhere else, so that $\bar{\E}_{h} = \Lab^{1}_{I}$, we find that
\begin{align}
\cae_{I} = \fracp{\dff (\ce)}{\ce_{I}} ,
\end{align}
and can therefore write
\begin{equation}
\fracd{\ff [\E_{h}]}{\E} = \sum_{I=1}^{3N_1} \fracp{\dff (\ce)}{\ce_{I}} \Sib^{1}_{I} (\xx) .
\end{equation}
On the other hand, expanding the dual basis in the original basis,
\begin{align}
\Sib^{1}_{I} (\xx) = \sum_{J=1}^{3N_1} \ca_{IJ} \Lab^{1}_{J} (\xx) ,
\end{align}
and taking the $L^2$ inner product with $\Lab^{1}_{J} (\xx)$, we find that the matrix $\MA=(\ca_{IJ})$ verifies $\MA\MM_{1}= \mathbb{I}_{3N_1}$, where $\mathbb{I}_{3N_1}$ denotes the $3N_1 \times 3N_1$ identity matrix, so that $\MA$ is the inverse of the mass matrix $\MM_{1}$ and
\begin{align}
\fracd{\ff [\E_{h}]}{\E} = \sum \limits_{I,J=1}^{3N_{1}} \fracp{\dff (\ce)}{\ce_{I}} \, (\MM_{1}^{-1})_{IJ} \, \Lab^{1}_{J} (\xx) .
\end{align}
In full analogy we find
\begin{align}
\fracd{\ff [\B_{h}]}{\B} = \sum \limits_{I,J=1}^{3N_{2}} \fracp{\dff (\cb)}{\cb_{I}} \, (\MM_{2}^{-1})_{IJ} \, \Lab^{2}_{J} (\xx) .
\end{align}
Next, using~\eqref{eq:rotation_matrix}, we find
\begin{align}
\curl \fracd{\ff [\E_{h}]}{\E} = \sum \limits_{I,J=1}^{3N_{1}} \sum \limits_{K=1}^{3N_{2}} \fracp{\dff (\ce)}{\ce_{I}} \, (\MM_{1}^{-1})_{IJ} \mathbb{C}_{JK} \Lab^{2}_{K} (\xx) .
\end{align}
Finally, we can re-express the following term in the Poisson bracket
\begin{align}
\nonumber
\int \curl &\, \fracd{\ff [\E_{h}]}{\E} \cdot  \fracd{\fg [\B_{h}]}{\B} \dx
= \\
\nonumber
&=\sum \limits_{I,J=1}^{3N_{1}} \sum_{K,L,M=1}^{3N_2}
\fracp{\dff (\ce)}{\ce_{I}} \, (\MM_{1}^{-1})_{IJ} \mathbb{C}_{JK}
\fracp{\dfg (\cb)}{\cb_{L}} \, (\MM_{2}^{-1})_{LM}
\int  \Lab^{2}_{K} (\xx) \cdot \Lab^{2}_{M} (\xx) \dx \\
&=\sum \limits_{I,J=1}^{3N_{1}} \sum_{K=1}^{3N_2}  
\fracp{\dff (\ce)}{\ce_{I}} \, (\MM_{1}^{-1})_{IJ} \mathbb{C}_{JK}
\fracp{\dfg (\cb)}{\cb_{K}}
 =\sum \limits_{I=1}^{3N_{1}} \sum_{K=1}^{3N_2}  
\fracp{\dff (\ce)}{\ce_{I}} \, (\MM_{1}^{-1} \mathbb{C})_{IK}
\fracp{\dfg (\cb)}{\cb_{K}} .
\end{align}
The other terms in the bracket involving functional derivatives with respect to the fields are handled similarly.
In the next step we need to discretize the distribution function and the corresponding functional derivatives.

\subsection{Discretization of the Functional Particle Derivatives}

We proceed by assuming a particle-like distribution function for $N_{p}$ particles labelled by~$a$,
\begin{align}
\f_{h} (\xx,\vv,t) = \sum \limits_{a=1}^{N_{p}} \, w_a \, \de \big( \xx - \xa (t) \big) \, \de \big( \vv - \va (t) \big) ,
\label{eq:vlasov_maxwell_distribution_discrete}
\end{align}
with mass $m_{a}$, charge $q_{a}$, weights $w_{a}$, particle positions $\xa$ and particle velocities $\va$. Here, $N_{p}$ denotes the total number of particles of all species, with each particle carrying its own mass and signed charge.
Functionals of the distribution function, $\ff[\f]$, can be considered as functions of the particle phase space trajectories, $F(\X,\V)$, upon inserting \eqref{eq:vlasov_maxwell_distribution_discrete},
\begin{align}\label{eq:vlasov_maxwell_functionals_of_f}
\ff [\f_{h}] = \dff (\X,\V).
\end{align}
Variation of the left-hand side and integration by parts gives,
\begingroup
\allowdisplaybreaks
\begin{align}
\nonumber
\de \ff[\f_{h}]
&= \int \fracd{\ff}{\f} \, \de \f_{h} \dx \dv \\
\nonumber
&= - \sum \limits_{a=1}^{N_{p}} w_a \int \dfrac{\de \ff}{\de \f} \Bigg(
     \de \big( \vv - \va (t) \big) \, \nabla_\xx \de \big( \xx - \xa (t) \big) \cdot \de \xa \\
\nonumber
&  \hspace{10em}
   + \de \big( \xx - \xa (t) \big) \, \nabla_\vv \de \big( \vv - \va (t) \big) \cdot \de \va \Bigg) \dx \dv \\
&= \sum \limits_{a=1}^{N_{p}} w_a \left(
     \nabla_\xx \fracd{\ff}{\f} \bigg\vert_{(\xa, \va)} \cdot \de \xa
   + \nabla_\vv \fracd{\ff}{\f} \bigg\vert_{(\xa, \va)} \cdot \de \va
   \right) .
\end{align}
\endgroup
Upon equating this expression with the variation of the right-hand side of~\eqref{eq:vlasov_maxwell_functionals_of_f},
\begin{align}
\de \dff (\X,\V)
&= \sum \limits_{a=1}^{N_{p}} \left( \fracp{\dff}{\xa} \, \de \xa + \fracp{\dff}{\va} \, \de \va \right) ,
\end{align}
we obtain
\begin{align}\label{eq:vlasov_maxwell_discrete_derivatives_f}
\fracp{\dff}{\xa} = w_a \nabla_\xx \fracd{\ff}{\f} \bigg\vert_{(\xa, \va)}
\hspace{3em}
\text{and}
\hspace{3em}
\fracp{\dff}{\va} = w_a \nabla_\vv \fracd{\ff}{\f} \bigg\vert_{(\xa, \va)} .
\end{align}
Considering the kinetic part of the Poisson bracket~\eqref{eq:poisson_bracket_maxwell}, the discretization proceeds in two steps.
First, replace $\f_s$ with $\f_{h}$ to get
\begin{multline}
\sum \limits_{a=1}^{N_{p}} \, \int \dfrac{w_a}{m_a} \, \de \big( \xx - \xa (t) \big) \, \de \big( \vv - \va (t) \big) \sbra{ \fracd{\ff}{\f} }{ \fracd{\fg}{\f} } \dx \dv =
\\
= \sum \limits_{a=1}^{N_{p}} \, \dfrac{w_a}{m_a} \bigg( \nabla_\xx \fracd{\ff}{\f} \cdot \nabla_\vv \fracd{\fg}{\f} - \nabla_\vv \fracd{\ff}{\f} \cdot \nabla_\xx \fracd{\fg}{\f} \bigg) \bigg\vert_{(\xa, \va)} .
\end{multline}
Then insert~\eqref{eq:vlasov_maxwell_discrete_derivatives_f} in order to obtain the discrete kinetic bracket.

\subsection{Discrete Poisson Bracket}

Replacing all functional derivatives in~\eqref{eq:poisson_bracket_maxwell} as outlined in the previous two sections, we obtain the semi-discrete Poisson bracket
\begingroup
\allowdisplaybreaks
\begin{align}
\{\dff,\dfg\} & [\X, \V, \ce, \cb]
\nonumber
 = \sum \limits_{a=1}^{N_{p}} \frac{1}{m_a w_a} \bigg(
       \fracp{\dff}{\xa} \cdot \fracp{\dfg}{\va}
     - \fracp{\dfg}{\xa} \cdot \fracp{\dff}{\va}
   \bigg) \\
\nonumber
&+ \sum \limits_{a=1}^{N_{p}} \sum \limits_{I,J=1}^{3N_{1}} \dfrac{q_a}{m_a} \left( 
        \fracp{\dfg}{\ce_{I}} \, (\MM_{1}^{-1})_{IJ} \, \Lab^{1}_{J} (\xa) \cdot \fracp{\dff}{\va}
      - \fracp{\dff}{\ce_{I}} \, (\MM_{1}^{-1})_{IJ} \, \Lab^{1}_{J} (\xa) \cdot \fracp{\dfg}{\va}
   \right) \\
\nonumber
&+ \sum \limits_{a=1}^{N_{p}} \sum \limits_{I=1}^{3N_{2}} \dfrac{q_a}{m_a^2 w_a} \, \cb_{I} (t) \, \Lab^{2}_{I} (\xa) \cdot \left( \fracp{\dff}{\va} \times \fracp{\dfg}{\va} \right) \\
&+ \sum \limits_{I,J=1}^{3N_{1}} \sum \limits_{K,L=1}^{3N_{2}} \left(
        \fracp{\dff}{\ce_{I}} \, (\MM_{1}^{-1})_{IJ} \, \C^{\top}_{JK} \, \fracp{\dfg}{\cb_{L}}
      - \fracp{\dfg}{\ce_{I}} \, (\MM_{1}^{-1})_{IJ} \, \C^{\top}_{JK} \,  \fracp{\dff}{\cb_{L}}
    \right) ,
\label{eq:poisson_bracket_maxwell_discrete}
\end{align}
\endgroup
with the curl matrix $\C$ as given in~\eqref{eq:rotation_matrix}.
In order to express the semi-discrete Poisson bracket~\eqref{eq:poisson_bracket_maxwell_discrete} in matrix form, we denote by $\LaB^1(\X)$ the $3N_p \times 3N_1$ matrix with generic term $\Lab^1_I(\xa)$, where $1 \leq a \leq N_p$ and $1 \leq I \leq 3N_1$, and by $\BB(\X, \cb)$ the $3N_p \times 3N_p$ block diagonal matrix with generic block
\begin{align}\label{eq:Bfield_twoform}
\widehat{\B}_{h} (\xa, t) = \sum_{i=1}^{N_2} 
\begin{pmatrix}
\hphantom{-} 0 & \hphantom{-} b_{i,3}(t) \, \La^{2,3}_{i} (\xa) & - b_{i,2}(t) \, \La^{2,2}_{i} (\xa) \\
- b_{i,3}(t) \, \La^{2,3}_{i} (\xa) & \hphantom{-} 0 & \hphantom{-} b_{i,1}(t) \, \La^{2,1}_{i} (\xa) \\
\hphantom{-} b_{i,2}(t) \, \La^{2,2}_{i} (\xa) & - b_{i,1}(t) \, \La^{2,1}_{i} (\xa) & \hphantom{-} 0 \\
\end{pmatrix} .
\end{align}
Further, let us introduce a mass matrix $M_{p}$ and a charge matrix $M_{q}$ for the particles.
Both are diagonal $N_{p} \times N_{p}$ matrices with elements $(M_{p})_{aa} = m_a w_a$ and $(M_{q})_{aa} = q_a w_a$, respectively. 
Additionally, we will need the $3N_{p} \times 3N_{p}$ matrices
\begin{align}
\MM_{p} = M_{p} \otimes \mathbb{I}_{3} ,
\hspace{3em}
\MM_{q} = M_{q} \otimes \mathbb{I}_{3} ,
\end{align}
where $\mathbb{I}_{3}$ denotes the $3 \times 3$ identity matrix.
This allows us to rewrite
\begin{align}
\nonumber
\sum \limits_{a=1}^{N_{p}} \sum \limits_{I=1}^{3N_{2}} \dfrac{q_a}{m_a^2 w_a} &\, \cb_{I} (t) \, \Lab^{2}_{I} (\xa) \cdot \left( \fracp{\dff}{\va} \times \fracp{\dfg}{\va} \right) = \\
\nonumber
&= \sum \limits_{a=1}^{N_{p}} \sum \limits_{i=1}^{N_{2}} \sum_{\mu=1}^{3} \dfrac{q_a}{m_a^2 w_a} \, b_{i,\mu} (t) \, \Lab^{2}_{i,\mu} (\xa) \cdot \left( \fracp{\dff}{\va} \times \fracp{\dfg}{\va} \right) \\
\nonumber
&= - \sum \limits_{a=1}^{N_{p}} \fracp{\dff}{\va} \dfrac{q_a}{m_a} \cdot \sum \limits_{i=1}^{N_{2}} \sum_{\mu=1}^{3} b_{i,\mu} (t) \, \Lab^{2}_{i,\mu} (\xa) \times \dfrac{1}{m_a w_a} \fracp{\dfg}{\va} \\
\nonumber
&= - \sum \limits_{a=1}^{N_{p}} \fracp{\dff}{\va} \dfrac{q_a w_a}{m_a w_a} \cdot \widehat{\B}_{h} (\xa, t) \cdot \dfrac{1}{m_a w_a} \fracp{\dfg}{\va} \\
&= - \bigg( \fracp{\dff}{\V} \bigg)^{\top} \MM_{p}^{-1} \MM_{q} \BB (\X, \cb) \, \MM_{p}^{-1} \bigg( \fracp{\dfg}{\V} \bigg) .
\end{align}
Here, the derivatives are represented by the $3 N_{p}$ vector
\begin{align}
\fracp{\dff}{\V}
= \left(
\fracp{\dff}{\vv_1} , \,
\hdots , \,
\fracp{\dff}{\vv_{N_{p}}} 
\right)^{\top}
= \left(
\fracp{\dff}{\cv_1^{1}} , \,
\fracp{\dff}{\cv_1^{2}} , \,
\fracp{\dff}{\cv_1^{3}} , \,
\hdots , \,
\fracp{\dff}{\cv_{N_{p}}^{1}} , \,
\fracp{\dff}{\cv_{N_{p}}^{2}} , \,
\fracp{\dff}{\cv_{N_{p}}^{3}} 
\right)^{\top} ,
\end{align}
and correspondingly for $\pa \dfg / \pa \V$, $\pa \dff / \pa \ce$, $\pa \dff / \pa \cb$, etc.
Thus, the discrete Poisson bracket (\ref{eq:poisson_bracket_maxwell_discrete}) becomes
\begin{align}
\{\dff,\dfg\} [\X, \V, \ce, \cb]
\nonumber
&= \fracp{\dff}{\X} \MM_{p}^{-1} \fracp{\dfg}{\V}
 - \fracp{\dfg}{\X} \MM_{p}^{-1} \fracp{\dff}{\V}
\\
\nonumber
&+ \bigg( \fracp{\dff}{\V}  \bigg)^{\top} \MM_{p}^{-1} \MM_{q} \LaB^1(\X)^{\top} \MM_{1}^{-1} \bigg( \fracp{\dfg}{\ce} \bigg)
\\
\nonumber
&- \bigg( \fracp{\dff}{\ce} \bigg)^{\top} \MM_{1}^{-1} \LaB^1(\X) \, \MM_{q} \MM_{p}^{-1}     \bigg( \fracp{\dfg}{\V}  \bigg)
\\
\nonumber
&+ \bigg( \fracp{\dff}{\V}  \bigg)^{\top} \MM_{p}^{-1} \MM_{q} \BB (\X, \cb) \, \MM_{p}^{-1} \bigg( \fracp{\dfg}{\V} \bigg)
\\
&+ \bigg( \fracp{\dff}{\ce} \bigg)^{\top} \MM_{1}^{-1} \C^{\top}  \bigg( \fracp{\dfg}{\cb} \bigg)
 - \bigg( \fracp{\dff}{\cb} \bigg)^{\top}  \C        \MM_{1}^{-1} \bigg( \fracp{\dfg}{\ce} \bigg)
.
\label{eq:poisson_bracket_maxwell_matrix}
\end{align}
The action of this bracket on two functions $\dff$ and $\dfg$ can also be expressed as
\begin{align*}
\{\dff,\dfg\} = \DD\dff^\top \jac (\uu) \, \DD\dfg ,
\end{align*}
denoting by $\DD$ the derivative with respect to the dynamical variables 
\begin{align}
\uu = (\X, \V, \ce, \cb)^{\top} ,
\end{align}
and by $\jac$ the Poisson matrix, given by
\begin{align}\label{eq:poisson_matrix_discrete}
\jac (\uu) =
\begin{pmatrix}
0 & \MM_p^{-1} & 0 & 0 \\
- \MM_p^{-1} & \MM_p^{-1} \MM_q \BB (\X,\cb) \, \MM_p^{-1} & \MM_p^{-1} \MM_q \LaB^1(\X) \MM_{1}^{-1} & 0 \\
0 & - \MM_{1}^{-1} \LaB^1 (\X)^\top \MM_q \MM_p^{-1} & 0 &  \MM_{1}^{-1} \C^\top\\
0 & 0 & - \C \MM_{1}^{-1} & 0 \\
\end{pmatrix} .
\end{align}
We immediately see that $\jac (\uu)$ is anti-symmetric, but it remains to be shown that it satisfies the Jacobi identity.

\subsection{Jacobi Identity}
\label{sec:jacobi_identity}

The discrete Poisson bracket~\eqref{eq:poisson_bracket_maxwell_matrix} satisfies the Jacobi identity if and only if the following condition holds (see e.g. \cite[Section IV]{Morrison:1998} or \cite[Section VII.2, Lemma 2.3]{HairerLubichWanner:2006}):
\begin{align}\label{eq:poisson_condition}
\sum \limits_{l} \left( 
     \fracp{\jac_{ij} (\uu)}{\cu_{l}} \, \jac_{lk} (\uu)
   + \fracp{\jac_{jk} (\uu)}{\cu_{l}} \, \jac_{li} (\uu)
   + \fracp{\jac_{ki} (\uu)}{\cu_{l}} \, \jac_{lj} (\uu)
\right) = 0
\quad \text{for all} \quad  i, j, k .
\end{align}
Here, all indices $i,j,k,l$ run from $1$ to $6 N_p + 3 N_1 + 3 N_2$.
To simplify the verification of~\eqref{eq:poisson_condition}, we start by identifying blocks of the Poisson matrix $\jac$ whose elements contribute to the above condition.
Therefore, we write
\begin{align}
\jac (\uu) = \begin{pmatrix}
J_{11} (\uu) & J_{12} (\uu) & J_{13} (\uu) & J_{14} (\uu) \\
J_{21} (\uu) & J_{22} (\uu) & J_{23} (\uu) & J_{24} (\uu) \\
J_{31} (\uu) & J_{32} (\uu) & J_{33} (\uu) & J_{34} (\uu) \\
J_{41} (\uu) & J_{42} (\uu) & J_{43} (\uu) & J_{44} (\uu) \\
\end{pmatrix}
 = \begin{pmatrix}
0 & J_{12} & 0 & 0 \\
J_{21} & J_{22} (\X, \cb) & J_{23} (\X) & 0 \\
0 & J_{32} (\X) & 0 & J_{34} \\
0 & 0 & J_{43} & 0 \\
\end{pmatrix} .
\end{align}
The Poisson matrix $\jac$ only depends on $\X$ and $\cb$, so in~\eqref{eq:poisson_condition} we only need to sum $l$ over the corresponding indices, $1 \leq l \leq 3 N_{p}$ and $6 N_{p} + 3 N_{1} < l \leq 6 N_{p} + 3 N_{1} + 3 N_{2}$, respectively. Considering the terms $\jac_{li} (\uu)$, $\jac_{lj} (\uu)$ and $\jac_{lk} (\uu)$, we see that in the aforementioned index ranges for $l$, only $J_{12} = \MM_{p}^{-1}$ and $J_{43} = - \MM_{2}^{-1} \C \MM_{1}^{-1}$ are non-vanishing, so that we have to account only for those two blocks, i.e., for $1 \leq l \leq 3 N_{p}$ we only need to consider $3 N_p < i,j,k \leq 6 N_p$ and for $6 N_{p} + 3 N_{1} < l \leq 6 N_{p} + 3 N_{1} + 3 N_{2}$ we only need to consider $6 N_p < i,j,k \leq 6 N_p + 3 N_1$. Note that $J_{12}$ is a diagonal matrix, therefore $(J_{12})_{ab} = (J_{12})_{aa} \delta_{ab}$ with $1 \leq a,b \leq N_p$. Further, only $J_{22}$, $J_{23}$ and $J_{32}$ depend on $\cb$ and/or $\X$, so only those blocks have to be considered when computing derivatives with respect to $\uu$.
In summary, we obtain two conditions.
The contributions involving $J_{22}$ and $J_{12}$ are
\begin{multline}\label{eq:poisson_condition1}
\sum \limits_{a=1}^{3 N_{p}} \Bigg( 
  \fracp{\big( J_{22} (\X, \cb) \big)_{bc}}{\X_{a}} \, \big( J_{12} \big)_{ad}
+ \fracp{\big( J_{22} (\X, \cb) \big)_{cd}}{\X_{a}} \, \big( J_{12} \big)_{ab} \\
+ \fracp{\big( J_{22} (\X, \cb) \big)_{db}}{\X_{a}} \, \big( J_{12} \big)_{ac}
\Bigg) = 0 ,
\end{multline}
for $1 \leq b, c, d \leq 3 N_{p}$, which corresponds to~\eqref{eq:poisson_condition} for $3 N_{p} < i, j, k \leq 6 N_{p}$.
Inserting the actual values for $J_{12}$ and $J_{22}$ and using that $\MM_{p}$ is diagonal, equation \eqref{eq:poisson_condition1} becomes
\begin{multline}\label{eq:poisson_condition1a}
  \fracp{\big( \MM_p^{-1} \MM_q \BB (\X,\cb) \, \MM_p^{-1} \big)_{bc}}{\X_{d}} \, \big( \MM_{p}^{-1} \big)_{dd}
+ \fracp{\big( \MM_p^{-1} \MM_q \BB (\X,\cb) \, \MM_p^{-1} \big)_{cd}}{\X_{b}} \, \big( \MM_{p}^{-1} \big)_{bb}
\\
+ \fracp{\big( \MM_p^{-1} \MM_q \BB (\X,\cb) \, \MM_p^{-1} \big)_{db}}{\X_{c}} \, \big( \MM_{p}^{-1} \big)_{cc}
= 0 .
\end{multline}
All outer indices of this expression belong to the inverse matrix $\MM_{p}^{-1}$. As this matrix is constant, symmetric and positive definite, we can contract the above expression with $\MM_{p}$ on all indices, to obtain
\begin{align}\label{eq:poisson_condition1b}
  \fracp{\big( \MM_q \BB (\X,\cb) \big)_{bc}}{\X_{d}}
+ \fracp{\big( \MM_q \BB (\X,\cb) \big)_{cd}}{\X_{b}}
+ \fracp{\big( \MM_q \BB (\X,\cb) \big)_{db}}{\X_{c}}
= 0 , \qquad 1 \leq b, c, d \leq 3 N_{p} .
\end{align}
If in the first term of~\eqref{eq:poisson_condition1b} one picks a particular index $k$, then this selects the $\sigma$ component of the position $\xa$ of some particle.
At the same time, in the second and third terms, this selects a block of~\eqref{eq:Bfield_twoform}, which is evaluated at the same particle position $\xa$. This means that the only non-vanishing contributions in~\eqref{eq:poisson_condition1b} will be those indexed by $b$ and $c$ with $\X_{b}$ and $\X_{c}$ corresponding to components $\mu$ and $\nu$ of the same particle position $\xa$.
Therefore, the condition~\eqref{eq:Bfield_twoform} reduces further to
\begin{align}
q_{a} \left( 
     \fracp{\widehat{B}_{\mu\nu} (\xa)}{\cx_{a}^{\si}}
   + \fracp{\widehat{B}_{\nu\si} (\xa)}{\cx_{a}^{\mu}}
   + \fracp{\widehat{B}_{\si\mu} (\xa)}{\cx_{a}^{\nu}}
\right) = 0 , \qquad 1 \leq a \leq N_{p} , \quad 1 \leq \mu, \nu, \si \leq 3 ,
\end{align}
where $\widehat{B}_{\mu\nu}$ denotes the components of the matrix in~\eqref{eq:Bfield_twoform}.
When all three indices are equal, this corresponds to diagonal terms of the matrix $\widehat{\B}_{h} (\xa, t) $ which vanish. When two of the three are equal, it cancels because of the skew-symmetry of the same matrix and for all three indices distinct, this condition corresponds to $\div \B_{h} = 0$.
Choosing initial conditions such that $\div \B_{h} (\xx,0) = 0$ and using a discrete deRham complex guarantees $\div \B_{h} (\xx,t) = 0$ for all times $t$.
Note that this was to be expected, because it is the discrete version of the $\div \B = 0$ condition for the continuous Poisson bracket~\eqref{eq:poisson_bracket_maxwell}~\cite{Morrison:1982}.
Further note that in the discrete case, just as in the continuous case, the Jacobi identity requires the magnetic field to be divergence-free, but it is not required to be the curl of some vector potential. In the language of differential forms this is to say that the magnetic field as a 2-form needs to be closed but does not need to be exact.

From the contributions involving $J_{22}$, $J_{23}$, $J_{32}$ and $J_{43}$ we have
\begin{multline}\label{eq:poisson_condition2}
\sum \limits_{a=1}^{3 N_{p}} \left(
    \fracp{\big( J_{23} (\X, \cb) \big)_{cI}}{\X_{a}} \, \big( J_{12} \big)_{ab}
  + \fracp{\big( J_{32} (\X, \cb) \big)_{Ib}}{\X_{a}} \, \big( J_{12} \big)_{ac}
\right) \\
+ \sum \limits_{J=1}^{N_{2}}
    \fracp{\big( J_{22} (\X, \cb) \big)_{bc}}{\cb_{J}} \, \big( J_{43} \big)_{JI}
= 0 ,
\end{multline}
for $1 \leq b, c \leq 3 N_{p}$ and $1 \leq I \leq 3 N_{1}$, which corresponds to~\eqref{eq:poisson_condition} for $3 N_{p} < i, j \leq 6 N_{p} < k \leq 6 N_{p} + 3 N_{1}$.
Writing out~\eqref{eq:poisson_condition2} and using that $\MM_{p}$ is diagonal gives
\begin{multline}\label{eq:poisson_condition2a}
    \fracp{\big( \MM_{1}^{-1} \LaB^1 (\X)^\top \MM_q \MM_p^{-1} \big)_{Ib}}{\X_{c}} \, \big( \MM_p^{-1} \big)_{cc} 
  - \fracp{\big( \MM_p^{-1} \MM_q \LaB^1(\X) \MM_{1}^{-1} \big)_{cI}}{\X_{b}} \, \big( \MM_p^{-1} \big)_{bb} =
  \\
= \sum \limits_{J=1}^{N_{2}} \fracp{\big( \MM_p^{-1} \MM_q \BB (\X,\cb) \, \MM_p^{-1} \big)_{bc}}{\cb_{J}} \, \big(  \C \MM_{1}^{-1} \big)_{JI} .
\end{multline}
Again, we can contract this with the matrix $\MM_{p}$ on the indices $b$ and $c$, in order to remove $\MM_{p}^{-1}$, and with $\MM_{1}$ on the index $I$, in order to remove $\MM_{1}^{-1}$. Similarly, $\MM_{q}$ can be removed by contracting with $\MM_{q}^{-1}$, noting that this matrix is constant and therefore commutes with the curl. This results in the simplified condition
\begin{align}\label{eq:poisson_condition2b}
    \fracp{\LaB^1_{bI} (\X)}{\X_{c}} 
  - \fracp{\LaB_{cI}^1(\X)}{\X_{b}}
= \sum \limits_{J=1}^{N_{2}} \fracp{\BB_{bc} (\X,\cb)}{\cb_{J}} \, \C_{JI} .
\end{align}
The $\cb_{J}$ derivative of $\BB$ results in the $3N_p \times 3N_p$ block diagonal matrix $\LaB^{2}_{J} (\X)$ with generic block
\begin{align}
\widehat{\Lab}^{2}_{J} (\xa) =
\begin{pmatrix}
\hphantom{-} 0 & \hphantom{-} \La^{2,3}_{J} (\xa) & - \La^{2,2}_{J} (\xa) \\
- \La^{2,3}_{J} (\xa) & \hphantom{-} 0 & \hphantom{-} \La^{2,1}_{J} (\xa) \\
\hphantom{-} \La^{2,2}_{J} (\xa) & - \La^{2,1}_{J} (\xa) & \hphantom{-} 0 \\
\end{pmatrix} ,
\end{align}
so that~\eqref{eq:poisson_condition2b} becomes
\begin{align}\label{eq:poisson_condition2c}
    \fracp{\big( \LaB^1 (\X)^\top \MM_q \big)_{Ib}}{\X_{c}} 
  - \fracp{\big( \MM_q \LaB^1(\X) \big)_{cI}}{\X_{b}}
= \sum \limits_{J=1}^{N_{2}} \C^{\top}_{IJ} \, \big( \MM_q \LaB^{2}_{J} (\X) \big)_{bc} .
\end{align}
This condition can be compactly written as
\begin{align}
\curl \LaB^1 (\X) = \LaB^{2} (\X) \C .
\end{align}
That the charge matrices $\MM_q$ can be eliminated can be seen as follows.
Similarly as before, picking a particular index $c$ in the first term selects the $\nu$ component of the position $\xa$ of some particle. At the same time, in the second term, this selects the $\nu$ component of $\Lab^{1}$, evaluated at the same particle position $\xa$. The only non-vanishing derivative of this term is therefore with respect to components of the same particle position, so that $\X_{b}$ denotes the $\mu$ component of $\xa$.
Hence, condition~\eqref{eq:poisson_condition2c} simplifies to
\begin{align}\label{eq:poisson_condition2d}
q_a \left( 
    \fracp{\Lab_{I}^{1,\mu} (\xa)}{\cx_a^{\nu}} 
  - \fracp{\Lab_{I}^{1,\nu} (\xa)}{\cx_a^{\mu}}
\right)
= q_a \sum \limits_{J=1}^{N_{2}} \C^{\top}_{IJ} \, \big( \widehat{\Lab}^{2}_{J} (\xa) \big)_{\mu\nu} , 
\end{align}
for $1 \leq a \leq N_{p}$, $1 \leq \mu, \nu \leq 3$ and $1 \leq I \leq 3 N_{1}$. The charge $q_{a}$ is the same on both sides and can therefore be removed.
This conditions states how the curl of the one-form basis, evaluated at some particle's position, is expressed in the two-form basis, evaluated at the same particle's position, using the curl matrix $\C$.
For spaces which build a deRham complex, this is always satisfied.
This concludes the verification of condition~\eqref{eq:poisson_condition} for the Jacobi identity to hold for the discrete bracket~\eqref{eq:poisson_bracket_maxwell_matrix}.

\subsection{Discrete Hamiltonian and Equations of Motion}

The Hamiltonian is discretized by inserting (\ref{eq:vlasov_maxwell_distribution_discrete}) and (\ref{eq:fields_discrete}) into \eqref{eq:hamiltonian_vlasov_maxwell},
\begin{align}
\dham (\V, \ce, \cb)
\nonumber
&\equiv \ham [\f_{h}, \E_{h}, \B_{h}] \\
\nonumber
&= \dfrac{1}{2} \int \abs{\vv}^{2} \, \sum \limits_{a=1}^{N_{p}} \, m_a w_a \, \de\big(\xx-\xa(t)\big) \, \de\big(\vv-\va(t)\big) \dx \dv \\
&+ \dfrac{1}{2} \int \bigg[ \Big\vert \sum \limits_{I=1}^{3N_{1}} \, \ce_I(t) \, \Lab^{1}_I (\xx) \Big\vert^{2}
 + \Big\vert \sum \limits_{J=1}^{3N_{2}} \, \cb_J(t) \, \Lab^{2}_J (\xx) \Big\vert^{2} \bigg] \dx ,
\end{align}
which in matrix notation becomes
\begin{align}
\dham
&= \tfrac{1}{2} \, \V^{\top}  \MM_{p} \V
 + \tfrac{1}{2} \, \ce^{\top} \MM_{1} \ce
 + \tfrac{1}{2} \, \cb^{\top} \MM_{2} \cb .
\label{eq:vlasov_maxwell_hamiltonian_discrete_matrix}
\end{align}
The semi-discrete equations of motion are given by
\begin{align}
\dot{\X}  &= \{ \X  , \dham \} , &
\dot{\V}  &= \{ \V  , \dham \} , &
\dot{\ce} &= \{ \ce , \dham \} , &
\dot{\cb} &= \{ \cb , \dham \} ,
\end{align}
which are equivalent to
\begin{align}
\dot{\uu} = \jac (\uu) \, \DD \dham (\uu) .
\end{align}
With $\DD \dham (\uu) = ( 0 , \, \MM_{p} \V , \, \MM_{1} \ce , \, \MM_{2} \cb)^{\top}$, we obtain
\begin{subequations}\label{eq:vlasov_maxwell_equations_of_motion}
\begin{align}
\label{eq:vlasov_maxwell_equations_of_motion_x}
\dot{\X}
&= \V ,
\\
\label{eq:vlasov_maxwell_equations_of_motion_v}
\dot{\V}
&= \MM_p^{-1} \MM_q \big( \LaB^1 (\X) \ce + \BB (\X,\cb) \V \big) ,
\\
\label{eq:vlasov_maxwell_equations_of_motion_e}
\dot{\ce}
&= \MM_{1}^{-1} \big( \C^\top \MM_{2} \cb (t) - \LaB^1 (\X)^\top \MM_q \V \big) ,
\\
\label{eq:vlasov_maxwell_equations_of_motion_b}
\dot{\cb}
&= - \C \ce (t) ,
\end{align}
\end{subequations}
where the first two equations describe the particle dynamics and the last two equations describe the evolution of the electromagnetic fields.
Note that $\MM_{p}$ and $\MM_{q}$ are diagonal matrices and that $\MM_p^{-1} \MM_q$ is nothing but the factor $q_a / m_a$ for the particle labelled by $a$. The purpose of introducing these matrices is solely to obtain a compact notation and to display the Poisson structure of the semi-discrete system. However, these matrices are neither explicitly constructed nor is there a need to invert them.

\subsection{Discrete Gauss' Law}

We will now show that the discrete bracket~\eqref{eq:poisson_bracket_maxwell_matrix} satisfies a discrete Gauss' law and a discrete continuity equation.
Multiplying \eqref{eq:vlasov_maxwell_equations_of_motion_e} by $\mathbb{G}^\top \MM_{1}$ on the left, we get
\begin{align}
\mathbb{G}^\top \MM_{1} \dot{\ce} (t) &= \mathbb{G}^\top \C^\top \MM_{2} \cb (t) - \mathbb{G}^\top \LaB^1 (\X (t))^\top \MM_q \V (t) .
\end{align}
As $\mathbb{C}\mathbb{G}=0$ from \eqref{eq:exseqmath_primal}, the first term on the right-hand side vanishes.
Observe that
\begin{align}\label{eq:rel_lambda}
\LaB^1 (\X) \G \psi = \grad \LaB^0 (\X) \psi \quad \forall \psi \in \R^{N_0} ,
\end{align}
and using $\dd \xa / \dd t = \va$ we find that
\begin{align}
\frac{\dd \Psi_h(\xa(t))}{\dd t} = \frac{\dd \xa(t)}{\dd t} \cdot \grad \Psi_h (\xa(t))
= \va (t) \cdot \grad \Psi_h(\xa(t)) ,
\end{align}
for any $\Psi_h \in V_0$ with $\grad \Psi_h \in V_1$, so that we obtain
\begin{align}\label{eq:rel_discrete_gauss_law}
\G^\top \MM_{1} \dot{\ce}
= - \G^\top \LaB^1(\X)^\top \MM_q \V
= - ( \grad \LaB^0 (\X) )^\top \MM_q \dfrac{\dd\X}{\dd t}
= - \frac{\dd \LaB^0 (\X)^\top}{\dd t} \MM_q \mathbb{1}_{N_p} ,
\end{align}
where $\mathbb{1}_{N_p}$ denotes the column vector with $N_p$ terms all being unity, needed for the sum over the particles when there is no velocity vector.
This shows that the discrete Gauss' law is conserved,
\begin{equation}\label{eq:mat_gauss}
\G^\top \MM_{1} \ce = -\LaB^0 (\X)^\top \MM_q \mathbb{1}_{N_p} .
\end{equation}
Moreover, note that Equation~\eqref{eq:rel_discrete_gauss_law} also contains a discrete version of the continuity equation~\eqref{eq:continuity},
\begin{align}\label{eq:mat_continuity}
\frac{\dd \crho}{\dd t} + \G^\top \cj = 0 ,
\end{align}
with the discrete charge and current density given by
\begin{align}
\crho = - \LaB^0 (\X)^\top \MM_q \mathbb{1}_{N_p} ,
~~~~
\cj = \LaB^1(\X)^\top \MM_q \V .
\end{align}
The conservation of this relation in the fully discrete system depends on the choice of the time stepping scheme.

\subsection{Discrete Casimir Invariants}

Let us now find the Casimir invariants of the semi-discrete Poisson structure. 
In order to obtain the discrete Casimir invariants, we assume that the discrete spaces in~\eqref{CGL} not only form a complex, so that $\Ima(\grad) \subseteq \Ker(\curl)$ and $\Ima(\curl) \subseteq \Ker(\div)$, but that they form an exact sequence, so that $\Ima(\grad) = \Ker(\curl)$ and $\Ima(\curl) = \Ker(\div)$, i.e., at each node in~\eqref{CGL}, the image of the previous operator is not only a subset of the kernel of the next operator but is exactly equal to the kernel of the next operator. We will then see that this requirement is not necessary for the identified functionals to be valid Casimir invariants. However, it is a useful assumption in their identification.

The Casimir invariants of the semi-discrete system are functions $\dcas (\X, \V, \ce, \cb)$ such that $\{ \dcas , \dff \}=0$ for any function $\dff$. In terms of the Poisson matrix $\jac$, this can be expressed as $\jac (\uu) \, \DD \dcas (\uu) = 0$. Upon writing this for each of the lines of $\jac$ of~\eqref{eq:poisson_matrix_discrete}, we find for the first line
\begin{align}
\MM_p^{-1} \DD_\V \dcas = 0.
\end{align}
This implies that $\dcas$ does not depend on $\V$, which we shall use in the sequel.
Then the third line simply becomes
\begin{align}
\MM_{1}^{-1} \C^\top \DD_\cb \dcas = 0 ~~~ \Rightarrow ~~~ \DD_\cb \dcas \in \Ker(\C^\top) .
\end{align}
Then, because of the exact sequence property, there exists $\bar{\cb} \in \R^{N_3}$ such that $\DD_\cb \dcas = \D^\top \bar{\cb}$. Hence all functions of the form
\begin{equation}\label{eq:casimir_divb}
\dcas (\cb)= \cb^\top \D^\top \bar{\cb} = \bar{\cb}^\top \D \cb,  ~~~~ \bar{\cb} \in \R^{N_3}
\end{equation}
are Casimirs, which means that $\D \cb$, the matrix form of $\div \B_{h}$, is conserved.

The fourth line, using that $\dcas$ does not depend on $\V$, becomes
\begin{align}
\C \MM_{1}^{-1} \DD_\ce \dcas = 0 ~~~ \Rightarrow ~~~ \MM_{1}^{-1} \DD_\ce \dcas \in \Ker(\C) .
\end{align}
Because of the exact sequence property there exists $\bar{\ce} \in \R^{N_1}$, such that $\DD_\ce \dcas = \MM_{1} \G \bar{\ce}$.
Finally, the second line couples $\ce$ and $\X$ and reads, upon multiplying by $M_p$,
\begin{align}
\DD_\X \dcas = \MM_q \LaB^1(\X) \MM_{1}^{-1} \DD_\ce \dcas
= \MM_q \LaB^1(\X) \G \bar{\ce} = \MM_q \grad \LaB_0(\X) \bar{\ce} ,
\end{align}
using the expression for $\DD_\ce \dcas$ derived previously and~\eqref{eq:rel_lambda}.
It follows that all functions of the form
\begin{equation}\label{eq:casimir_dive}
\dcas (\X, \ce) = \mathbb{1}_N^\top \MM_q \LaB^0(\X) \bar{\ce} + \ce^\top \MM_{1} \G \bar{\ce}
= \bar{\ce}^\top \LaB^0(\X)^\top \MM_q \mathbb{1}_N + \bar{\ce}^\top \G^\top \MM_{1} \ce , ~~~~
\bar{\ce} \in \R^{N_0} ,
\end{equation}
are Casimirs, so that $\G^\top \MM_{1} \ce + \LaB^0 (\X)^\top \MM_q \mathbb{1}_N$ is conserved. This is the matrix form of Gauss' law~\eqref{eq:mat_gauss}.

Having identified the discrete Casimir invariants~\eqref{eq:casimir_divb} and~\eqref{eq:casimir_dive}, it is easy to see by plugging them into the discrete Poisson bracket~\eqref{eq:poisson_bracket_maxwell_matrix} that they are valid Casimir invariants, even if the deRham complex is not exact, because all that is needed for $\jac (\uu) \, \DD \dcas (\uu) = 0$ is the complex property $\Ima \D^\top \subseteq \Ker \C^\top$ and $\grad \LaB^0 (\X) = \LaB^1(\X) \G$.

\section{Hamiltonian Splitting}\label{sec:splitting}

Following \cite{Crouseilles:2015, Qin:2015, He:2015}, we split the discrete Hamiltonian \eqref{eq:vlasov_maxwell_hamiltonian_discrete_matrix} into three parts,
\begin{align}
\dham = \dham_{p} + \dham_{E} + \dham_{B} ,
\end{align}
with
\begin{align}
\dham_{p} &= \tfrac{1}{2} \, \V^{\top}  \MM_{p} \V  , &
\dham_{E} &= \tfrac{1}{2} \, \ce^{\top} \MM_{1} \ce , &
\dham_{B} &= \tfrac{1}{2} \, \cb^{\top} \MM_{2} \cb .
\end{align}
Writing $\uu = (\X, \V, \ce, \cb)^{\top}$, we split the discrete Vlasov--Maxwell equations~\eqref{eq:vlasov_maxwell_equations_of_motion} into three subsystems,
\begin{align}
\dot{\uu} &= \{ \uu , \dham_{p} \} , &
\dot{\uu} &= \{ \uu , \dham_{E} \} , &
\dot{\uu} &= \{ \uu , \dham_{B} \} .
\end{align}
The exact solution to each of these subsystems will constitute a Poisson map.
Because a composition of Poisson maps is itself a Poisson map, we can construct Poisson structure-preserving integration methods for the Vlasov--Maxwell system by composition of the exact solutions of each of the subsystems.

\subsection{Solution of the Subsystems}

The discrete equations of motion for $\dham_{E}$ are
\begin{subequations}
\begin{align}
\dot{\X}  &= 0 , \\
\MM_{p} \dot{\V}  &= \MM_q \LaB^1 (\X) \ce , \\
\dot{\ce} &= 0 , \\
\dot{\cb} &= - \C \ce (t) .
\end{align}
\end{subequations}
For initial conditions $\big( \X(0), \V(0), \ce(0), \cb(0) \big)$ the exact solutions at time $\Dt$ are given by the map $\phy_{\Dt,E}$ defined as
\begin{subequations}
\begin{align}
\X (\Dt) &= \X (0) , \\
\MM_{p} \V (\Dt) &= \MM_{p} \V (0) + \Dt \, \MM_q \LaB^1 (\X (0)) \ce (0) , \\
\ce(\Dt) &= \ce(0) , \\
\cb(\Dt) &= \cb(0) - \Dt \, \C \ce (0) .
\end{align}
\end{subequations}
The discrete equations of motion for $\dham_{B}$ are
\begin{subequations}
\begin{align}
\dot{\X} &= 0 , \\
\dot{\V} &= 0 , \\
\MM_{1} \dot{\ce} &= \C^\top \MM_{2} \cb (t) , \\
\dot{\cb}&= 0 .
\end{align}
\end{subequations}
For initial conditions $\big( \X(0), \V(0), \ce(0), \cb(0) \big)$ the exact solutions at time $\Dt$ are given by the map $\phy_{\Dt,B}$ defined as
\begin{subequations}
\begin{align}
\X (\Dt) &= \X (0) , \\
\V (\Dt) &= \V (0) , \\
\MM_{1} \ce(\Dt) &= \MM_{1} \ce(0) + \Dt \, \C^\top \MM_{2} \cb (0) , \\
\cb(\Dt) &= \cb(0) .
\end{align}
\end{subequations}
The discrete equations of motion for $\dham_{p}$ are
\begin{subequations}
\begin{align}
\dot{\X}
&= \V ,
\\
\MM_p \dot{\V}
&= \MM_q \BB (\X,\cb) \V , \\
\MM_{1} \dot{\ce}
&= - \LaB^1 (\X)^\top \MM_q \V ,
\\
\dot{\cb}
&= 0 .
\end{align}
\end{subequations}
For general magnetic field coefficients $\cb$, this system cannot be exactly integrated~\cite{He:2015}.
Note that each component $\dot{\V}_{\mu}$ of the equation for $\dot{\V}$ does not depend on $\V_{\mu}$, where $\V_{\mu} = (\cv_{1,\mu}, \cv_{2,\mu}, \hdots, \cv_{N_{p},\mu} )^\top$, etc., with $1 \leq \mu \leq 3$. Therefore we can split this system once more into
\begin{align}
\dham_{p} = \dham_{p_1} + \dham_{p_2} + \dham_{p_3} ,
\end{align}
with
\begin{align}
\dham_{p_\mu} &= \tfrac{1}{2} \, \V_{\mu}^{\top} M_{p} \V_{\mu} 
\qquad \text{for} \qquad
1 \leq \mu \leq 3 .
\end{align}
For concise notation we introduce the $N_p \times N_1$ matrix $\LaB^{1}_{\mu} (\X)$ with generic term $\La^{1}_{i,\mu} (\xa)$, and the $N_p \times N_p$ diagonal matrix $\LaB^{2}_{\mu} (\cb, \X)$ with entries $\sum_{i=1}^{N_2} b_i(t) \La^{2}_{i,\mu} (\xa)$, where $1 \leq \mu \leq 3$, $1 \leq a \leq N_p$, $1 \leq i \leq N_1$, $1 \leq j \leq N_2$.
Then, for $\dham_{p_1}$ we have
\begin{subequations}
\begin{align}
\dot{\X}_{1}     &= \hphantom{-} \V_{1} (t) , \\
M_p \dot{\V}_{2} &=           -  M_q \LaB^{2}_{3} (\cb(t),\X(t))^\top \, \V_{1} (t) , \\
M_p \dot{\V}_{3} &= \hphantom{-} M_q \LaB^{2}_{2} (\cb(t),\X(t))^\top \, \V_{1} (t) , \\
\MM_{1} \dot{\ce}    &=           -  \LaB^{1}_{1} (\X(t))^\top M_q \V_{1} (t) ,
\end{align}
\end{subequations}
for $\dham_{p_2}$ we have
\begin{subequations}
\begin{align}
\dot{\X}_{2}     &= \hphantom{-} \V_{2} (t) , \\
M_p \dot{\V}_{1} &= \hphantom{-} M_q \LaB^{2}_{3} (\cb(t),\X(t))^\top \, \V_{2} (t) , \\
M_p \dot{\V}_{3} &=           -  M_q \LaB^{2}_{1} (\cb(t),\X(t))^\top \, \V_{2} (t) , \\
\MM_{1}\dot{\ce}     &=           -  \LaB^{1}_{2} (\X(t))^\top M_q \V_{2} (t) ,
\end{align}
\end{subequations}
and for $\dham_{p_3}$ we have
\begin{subequations}
\begin{align}
\dot{\X}_{3}     &= \hphantom{-} \V_{3} (t) , \\
M_p \dot{\V}_{1} &=           -  M_q \LaB^{2}_{2} (\cb(t),\X(t))^\top \, \V_{3} (t) , \\
M_p \dot{\V}_{2} &= \hphantom{-} M_q \LaB^{2}_{1} (\cb(t),\X(t))^\top \, \V_{3} (t) , \\
\MM_{1} \dot{\ce}    &=           -  \LaB^{1}_{3} (\X(t))^\top M_q \V_{3} (t) .
\end{align}
\end{subequations}
For $\dham_{p_1}$ and initial conditions $\big( \X(0), \V(0), \ce(0), \cb(0) \big)$ the exact solutions at time $\Dt$ are given by the map $\phy_{\Dt,p_1}$ defined as
\begin{subequations}
\begin{align}
\X_{1} (\Dt) &= \X_{1} (0) + \Dt \, \V_{1}(0) , \\
M_{p} \V_{2} (\Dt) &= M_{p} \V_{2} (0) - \int \limits_{0}^{\Dt} M_{q} \LaB^{2}_{3} (\cb(0),\X(t)) \, \V_{1} (0) \dt , \\
M_{p} \V_{3} (\Dt) &= M_{p} \V_{3} (0) + \int \limits_{0}^{\Dt} M_{q} \LaB^{2}_{2} (\cb(0),\X(t)) \, \V_{1} (0) \dt , \\
\MM_{1} \ce(\Dt)     &= \MM_{1} \ce    (0) - \int \limits_{0}^{\Dt} \LaB^{1}_{1} (\X(t))^\top M_{q} \V_{1} (0) \dt ,
\end{align}
\end{subequations}
for $\dham_{p_2}$ by the map $\phy_{\Dt,p_2}$ defined as
\begin{subequations}
\begin{align}
\X_{2} (\Dt) &= \X_{2} (0) + \Dt \, \V_{2}(0) , \\
M_{p} \V_{1} (\Dt) &= M_{p} \V_{1} (0) + \int \limits_{0}^{\Dt} M_{q} \LaB^{2}_{3} (\cb(0),\X(t)) \, \V_{2} (0) \dt , \\
M_{p} \V_{3} (\Dt) &= M_{p} \V_{3} (0) - \int \limits_{0}^{\Dt} M_{q} \LaB^{2}_{1} (\cb(0),\X(t)) \, \V_{2} (0) \dt , \\
\MM_{1} \ce (\Dt)    &= \MM_{1} \ce    (0) - \int \limits_{0}^{\Dt} \LaB^{1}_{2} (\X(t))^\top M_{q} \V_{2} (0) \dt ,
\end{align}
\end{subequations}
and for $\dham_{p_3}$ by the map $\phy_{\Dt,p_3}$ defined as
\begin{subequations}
\begin{align}
\X_{3} (\Dt) &= \X_{3} (0) + \Dt \, \V_{3}(0) , \\
M_{p} \V_{1} (\Dt) &= M_{p} \V_{1} (0) - \int \limits_{0}^{\Dt} M_{q} \LaB^{2}_{2} (\cb(0),\X(t)) \, \V_{3} (0) \dt , \\
M_{p} \V_{2} (\Dt) &= M_{p} \V_{2} (0) + \int \limits_{0}^{\Dt} M_{q} \LaB^{2}_{1} (\cb(0),\X(t)) \, \V_{3} (0) \dt , \\
\MM_{1} \ce (\Dt)    &= \MM_{1} \ce    (0) - \int \limits_{0}^{\Dt} \LaB^{1}_{3} (\X(t))^\top M_{q} \V_{3} (0) \dt ,
\end{align}
\end{subequations}
respectively, where all components not specified are constant.
The only challenge in solving these equations is the exact computation of line integrals along the trajectories~\cite{Campos.Pinto.Jund.Salmon.Sonnendrucker.2014.crm, Squire:2012, moon2015exact}. However, because only one component of the particle positions $\xx_{a}$ is changing in each step of the splitting and, moreover, the trajectory during one sub step of the splitting is approximated by a straight line, this is not very complicated.
Compared to standard particle-in-cell methods for the Vlasov--Maxwell system, the exact integration causes slightly increased computational costs. These are, however, comparable to existing charge-conserving algorithms like that of~\citet{Villasenor.Buneman.1992.cpc} or the Boris correction method~\cite{Boris:1970}.

\subsection{Splitting Methods}\label{sec:composition}

Given initial conditions $\uu (0) = \big( \X(0), \V(0), \ce(0), \cb(0) \big)^{\top}$, a numerical solution of the discrete Vlasov--Maxwell equations \eqref{eq:vlasov_maxwell_equations_of_motion_x}-\eqref{eq:vlasov_maxwell_equations_of_motion_b} at time $\Delta t$ can be obtained by composition of the exact solutions of all the subsystems.
A first order integrator can be obtained by the Lie--Trotter composition~\cite{Trotter:1959},
\begin{align}\label{eq:splitting_lie_trotter}
\phy_{\Dt,L}
         = \phy_{\Dt,p_3}
     \circ \phy_{\Dt,p_2}
     \circ \phy_{\Dt,p_1}
     \circ \phy_{\Dt,B}
     \circ \phy_{\Dt,E} .
\end{align}
A second order integrator can be obtained by the symmetric Strang composition~\cite{Strang:1968},
\begin{align}\label{eq:splitting_strang}
\phy_{\Dt,S2} = \phy_{\Dt/2,L} \circ \phy_{\Dt/2,L}^{*} ,
\end{align}
where $\phy_{\Dt,L}^{*}$ denotes the adjoint of $\phy_{\Dt,L}$, defined as
\begin{align}
\phy_{\Dt,L}^{*} = \phy_{-\Dt,L}^{-1} .
\end{align}
Explicitly, the Strang splitting can be written as
\begin{multline}
\phy_{\Dt,S2} = \phy_{\Dt/2,E}
          \circ \phy_{\Dt/2,B}
          \circ \phy_{\Dt/2,p_1}
          \circ \phy_{\Dt/2,p_2}
          \circ \phy_{\Dt/2,p_3} \\
          \circ \phy_{\Dt/2,p_3}
          \circ \phy_{\Dt/2,p_2}
          \circ \phy_{\Dt/2,p_1}
          \circ \phy_{\Dt/2,B}
          \circ \phy_{\Dt/2,E} .
\end{multline}
Let us note that the Lie splitting $\phy_{\Dt,L}$ and the Strang splitting $\phy_{\Dt,S2}$ are conjugate methods by the adjoint of $\phy_{\Dt,L}$ \cite{McLachlanQuispel:2002}, i.e.,
\begin{align}
\phy_{\Dt,S2}
= ( \phy_{\Dt/2,L}^{*} )^{-1} \circ \phy_{\Dt,L} \circ \phy_{\Dt/2,L}^{*} 
= \phy_{-\Dt/2,L} \circ \phy_{\Dt,L} \circ \phy_{\Dt/2,L}^{*}
= \phy_{\Dt/2,L} \circ \phy_{\Dt/2,L}^{*} .
\end{align}
The last equality holds by the group property of the flow, but is only valid when the exact solution of each subsystem is used in the composition (and not some general symplectic or Poisson integrator). This implies that the Lie splitting shares many properties with the Strang splitting which are not found in general first order methods.

Another second order integrator with a smaller error constant than $\phy_{\Dt,S2}$ can be obtained by the following composition,
\begin{align}\label{eq:splitting_4lie}
\phy_{\Dt,L2} = \phy_{\alpha \Dt,L} \circ \phy_{(1/2-\alpha) \Dt,L}^{*} \circ \phy_{(1/2-\alpha) \Dt,L} \circ \phy_{\alpha \Dt,L}^{*} .
\end{align}
Here, $\alpha$ is a free parameter which can be used to reduce the error constant. A particularly small error is obtained for $\alpha = 0.1932$~\cite{McLachlan:1995}.
Fourth order time integrators can easily be obtained from a second order integrator like $\phy_{\Dt,S}$ by the following composition~\cite{Suzuki:1990, Yoshida:1990},
\begin{align}\label{eq:splitting_3strang}
\phy_{\Dt,S4} = \phy_{\gamma_{1} \Dt,S2} \circ \phy_{\gamma_{2} \Dt,S2} \circ \phy_{\gamma_{1} \Dt,S2} ,
\end{align}
with
\begin{align*}
\gamma_{1} &= \dfrac{1}{2 - 2^{1/3}} , &
\gamma_{2} &= - \dfrac{2^{1/3}}{2 - 2^{1/3}} .
\end{align*}
Alternatively, we can compose the first order integrator $\phy_{\Dt,L}$ together with its adjoint $\phy_{\Dt,L}^{*}$ as follows~\cite{McLachlan:1995},
\begin{align}\label{eq:splitting_10lie}
\phy_{\Dt,L4} = \phy_{a_{5} \Dt, L}
        \circ \phy_{b_{5} \Dt, L}^{*}
		\circ \hdots
        \circ \phy_{a_{2} \Dt, L}
        \circ \phy_{b_{2} \Dt, L}^{*}
        \circ \phy_{a_{1} \Dt, L}
        \circ \phy_{b_{1} \Dt, L}^{*} ,
\end{align}
with
\begin{align*}
a_{1} &= b_{5} = \dfrac{146 + 5 \sqrt{19}}{540} , &
a_{2} &= b_{4} = \dfrac{-2 + 10 \sqrt{19}}{135} , &
a_{3} &= b_{3} = \dfrac{1}{5} , \\
a_{4} &= b_{2} = \dfrac{-23 - 20 \sqrt{19}}{270} , &
a_{5} &= b_{1} = \dfrac{14 - \sqrt{19}}{108} .
\end{align*}
For higher order composition methods see e.g. \cite{HairerLubichWanner:2006} and \cite{McLachlanQuispel:2002} and references therein.

\subsection{Backward Error Analysis}\label{sec:splitting_bae}

In the following, we want to compute the modified Hamiltonian for the Lie--Trotter composition~\eqref{eq:splitting_lie_trotter}.
For a splitting in two terms, $\dham = \dham_{A} + \dham_{B}$, the Lie--Trotter composition can be written as
\begin{align}
\phy_{\Dt} = \phy_{\Dt,B} \circ \phy_{\Dt,A}
= \exp ( \Dt \, X_{A} )
  \exp ( \Dt \, X_{B} )
= \exp ( \Dt \, \tilde{X} ) ,
\end{align}
where $X_{A}$ and $X_{B}$ are the Hamiltonian vector fields corresponding to $\dham_{A}$ and $\dham_{B}$, respectively, and $\tilde{X}$ is the modified vector field corresponding to the modified Hamiltonian $\tilde{\dham}$.
Following \citet[Chapter~IX]{HairerLubichWanner:2006}, the modified Hamiltonian $\tilde{\dham}$ is given by
\begin{align}
\tilde{\dham} = \dham + \Dt \tilde{\dham}_{1} + \Dt^{2} \tilde{\dham}_{2} + \order{\Dt^{3}} ,
\end{align}
with
\begin{align}
\tilde{\dham}_{1} &= \dfrac{1}{2} \{ \dham_{A} , \dham_{B} \} , \\
\tilde{\dham}_{2} &= \dfrac{1}{12} \big[ \{ \{ \dham_{A} , \dham_{B} \} , \dham_{B} \} + \{ \{ \dham_{B} , \dham_{A} \} , \dham_{A} \} \big] .
\end{align}
In order to compute the modified Hamiltonian for splittings with more than two terms, we have to apply the Baker--Campbell--Hausdorff formula recursively.
Denoting by $X_{E}$ the Hamiltonian vector field corresponding to $\dham_{E}$, that is $X_{E} = \jac ( \cdot , \DD \dham_{E})$, and similar for $X_{B}$ and $X_{p_i}$, the Lie--Trotter splitting~\eqref{eq:splitting_lie_trotter} can be expressed in terms of these Hamiltonian vector fields as
\begin{align}
\phy_{\Dt,L}
= \exp ( \Dt \, X_{E } )
  \exp ( \Dt \, X_{B } )
  \exp ( \Dt \, X_{p_1} )
  \exp ( \Dt \, X_{p_2} )
  \exp ( \Dt \, X_{p_3} ) .
\end{align}
Let us split $\phy_{\Delta t} = \exp ( \Dt \, \tilde{X} )$ into
\begin{align*}
\exp ( \Dt \, X_{E  } ) \exp ( \Dt \, X_{B  } ) &= \exp ( \Dt \, \tilde{X}_{EB     } ) , \\
\exp ( \Dt \, X_{p_1} ) \exp ( \Dt \, X_{p_2} ) &= \exp ( \Dt \, \tilde{X}_{p_{1,2}} ) , \\
\exp ( \Dt \, \tilde{X}_{p_{1,2}} ) \exp ( h X_{p_3} ) &= \exp ( \Dt \, \tilde{X}_{p_{1,2,3}} ) , \\
\exp ( \Dt \, \tilde{X}_{EB  } ) \exp ( \Dt \, \tilde{X}_{p_{1,2,3}} ) &= \exp ( \Dt \, \tilde{X} ) ,
\end{align*}
where the corresponding Hamiltonians are given by
\begin{align*}
\tilde{\dham}_{EB } &= \dham_{E  } + \dham_{B  } + \dfrac{\Dt}{2} \{ \dham_{E  } , \dham_{B  } \} + \order{\Dt^{2}} , \\
\tilde{\dham}_{p_{1,2}} &= \dham_{p_1} + \dham_{p_2} + \dfrac{\Dt}{2} \{ \dham_{p_1} , \dham_{p_2} \} + \order{\Dt^{2}} , \\
\nonumber
\tilde{\dham}_{p_{1,2,3}} &= \tilde{\dham}_{p_{1,2}} + \dham_{p_3} + \dfrac{\Dt}{2} \{ \tilde{\dham}_{p_{1,2}} , \dham_{p_3} \} + \order{\Dt^{2}} \\
 &= \dham_{p_1} + \dham_{p_2} + \dham_{p_3} + \dfrac{\Dt}{2} \{ \dham_{p_1} , \dham_{p_2} \} + \dfrac{\Dt}{2} \{ \dham_{p_1} + \dham_{p_2} , \dham_{p_3} \} + \order{\Dt^{2}} .
\end{align*}
The Hamiltonian $\tilde{\dham}$ corresponding to $\tilde{X}$ is given by
\begin{align}
\tilde{\dham} &= \dham_{E} + \dham_{B} + \dham_{p_1} + \dham_{p_2} + \dham_{p_3} + \Dt \, \tilde{\dham}_{1} + \order{\Dt^{2}} ,
\end{align}
with the first order correction $\tilde{\dham}_{1}$ obtained as
\begin{align}
\tilde{\dham}_{1} &= \dfrac{1}{2} \big(
    \{ \dham_{E  } , \dham_{B  } \}
  + \{ \dham_{p_1} , \dham_{p_2} \}
  + \{ \dham_{p_1} + \dham_{p_2} , \dham_{p_3} \}
  + \{ \dham_{E  } + \dham_{B  } , \dham_{p  } \}
\big) .
\end{align}
The various Poisson brackets are computed as follows,
\begin{align*}
\{ \dham_{E  } , \dham_{B}   \} &=   \ce^{\top} \C^{\top} \cb , \\
\{ \dham_{E  } , \dham_{p}   \} &= - \ce^\top \LaB^1 (\X)^\top \MM_q \V , \\
\{ \dham_{B  } , \dham_{p}   \} &= 0 , \\
\{ \dham_{p_1} , \dham_{p_2} \} &= \V_{1} M_q \BB_{3} (\cb,\X)^\top \, \V_{2} , \\
\{ \dham_{p_2} , \dham_{p_3} \} &= \V_{2} M_q \BB_{1} (\cb,\X)^\top \, \V_{3} , \\
\{ \dham_{p_3} , \dham_{p_1} \} &= \V_{3} M_q \BB_{2} (\cb,\X)^\top \, \V_{1} ,
\end{align*}
where $\BB_{\mu} (\X, \cb)$ denotes $N_p \times N_p$ diagonal matrix with elements $B_{h,\mu} (\xa)$.
The Lie--Trotter integrator~\eqref{eq:splitting_lie_trotter} preserves the modified energy $\tilde{\dham} = \dham + \Dt \, \tilde{\dham}_{1}$ to $\mcal{O} (\Dt^{2})$, while the original energy $\dham$ is preserved only to $\mcal{O} (\Dt)$.

\section{Example: Vlasov--Maxwell in 1D2V}\label{sec:example}

Starting from the Vlasov equation for a particle species $s$, charge $q_s$, and mass $m_s$ given in~\eqref{eq:vlasov}, we reduce by assuming a single species and 
\begin{align}\label{cordRed}
\xx = (x,0,0) \quad \text{and} \quad \vv = (v_1,v_2,0)
\end{align}
as well as
\begin{align}\label{fieldRed}
\E(x,t) &= (E_1,E_2,0) , &
\B(x,t) &= (0,0,B_3) , &
\f(\xx,\vv,t) &= \f(x,v_1,v_2,t) ,
\end{align}
so that the Vlasov equation takes the form
\begin{equation}\label{eq:vlasov_1d2v}
\fracp{\f(x,\vb,t)}{t} + v_1 \fracp{\f(x,\vb,t)}{x} + \frac{q_s}{m_s} \bigg[ \Eb(x,t) + B_3(x,t) \begin{pmatrix}
\hphantom{-} v_2 \\
 - v_1
\end{pmatrix} \bigg] \cdot \nabla_{\vb} f(x,v,t) = 0 ,
\end{equation}
while Maxwell's equations become
\begin{align}
\fracp{E_1(x,t)}{t} &= - J_1 (x) , \label{eq:maxwell_1d2v_Ex} \\
\fracp{E_2(x,t)}{t} &= - \fracp{B  (x,t)}{x} - J_2 (x) \label{eq:maxwell_1d2v_Ey} , \\
\fracp{B  (x,t)}{t}	&= - \fracp{E_2(x,t)}{x} \label{eq:maxwell_1d2v_Bz} , \\
\fracp{E_1(x,t)}{x} &= \rho + \rho_B \label{eq:maxwell_1d2v_last} ,
\end{align}
with sources given by
\begin{align}\label{eq:maxwell_1d2v_sources}
\rho &= q_s \int \f \dv , &
J_1 &=  q_s \int \f v_1 \dv , &
J_2 &=  q_s \int \f v_2 \dv .
\end{align}
Note that $\div \B = 0$ is manifest.

\subsection{Noncanonical Hamiltonian Structure}

Under the assumptions of \eqref{cordRed} and \eqref{fieldRed}, the bracket of~\eqref{eq:poisson_bracket_maxwell} reduces to 
\begin{align}\label{eq:poisson_bracket_maxwell_1d2v}
\{\ff,\fg\}  
\nonumber
&= \frac{1}{m} \int \f \left[ \frac{\de \ff}{\de \f} , \frac{\de \fg}{\de \f} \right] \dx[1] \dv[1] \dv[2] \\
\nonumber
&+ \frac{q}{m} \int \f \bigg(
        \fracp{}{v_1} \frac{\de \ff}{\de \f} \frac{\de \fg}{\de E_1}
      + \fracp{}{v_2} \frac{\de \ff}{\de \f} \frac{\de \fg}{\de E_2}
      - \fracp{}{v_1} \frac{\de \fg}{\de \f} \frac{\de \ff}{\de E_1}
      - \fracp{}{v_2} \frac{\de \fg}{\de \f} \frac{\de \ff}{\de E_2}
   \bigg) \dx[1] \dv[1] \dv[2] \\
\nonumber
&+ \frac{q}{m^2} \int \f \, B \bigg( 
        \fracp{}{v_1} \frac{\de \ff}{\de \f} \fracp{}{v_2} \frac{\de \fg}{\de \f} 
      - \fracp{}{v_2} \frac{\de \ff}{\de \f} \fracp{}{v_1} \frac{\de \fg}{\de \f} 
   \bigg) \dx[1] \dv[1] \dv[2] \\
&+ \int \bigg( 
        \frac{\de \fg}{\de E_2} \fracp{}{x} \frac{\de \ff}{\de B}
      - \frac{\de \ff}{\de E_2} \fracp{}{x} \frac{\de \fg}{\de B}
   \bigg) \dx[1] ,
\end{align}
where now 
\begin{align}
[f,g] := \fracp{f}{x} \fracp{g}{v_1} - \fracp{f}{v_1} \fracp{g}{x} .
\end{align}
The Hamiltonian~\eqref{eq:hamiltonian_vlasov_maxwell} becomes
\begin{align}
\ham &= \frac{m}{2} \int \f (x, v_1, v_2) \, (v_1^{2} + v_2^{2}) \dx[1] \dv[1] \dv[2] + \frac{1}{2} \int \big( E_1^{2} + E_2^{2} + B^{2} \big) \dx[1] .
\label{RNRham}
\end{align}
The bracket of~\eqref{eq:poisson_bracket_maxwell_1d2v} with Hamiltonian~\eqref{RNRham}, generates~\eqref{eq:maxwell_1d2v_Ex}, \eqref{eq:maxwell_1d2v_Ey} and~\eqref{eq:maxwell_1d2v_Bz}, with the $\J$ of \eqref{eq:maxwell_1d2v_sources}.

\subsection{Reduced Jacobi Identity}
\label{sec:jacobi_identity_1d2v}

The bracket of~\eqref{eq:poisson_bracket_maxwell_1d2v} can be shown to satisfy the Jacobi identity by direct calculations. However since~\eqref{eq:poisson_bracket_maxwell} satisfies the Jacobi identity for all functionals, it must also satisfy it for all reduced functionals.  Here we  have closure, i.e., if $\ff$ and $\fg$ are reduced functionals, then $\{\ff,\fg\}$ is a reduced functional, where the bracket is the full bracket of~\eqref{eq:poisson_bracket_maxwell}.

More specifically, to understand this closure, observe that the  reduced functionals have the form
\begin{align}\label{eq:reduced_functional}
\ff [E_1,E_2, B, f]
\nonumber
&= \int \ff(\xx, \E, \B, \f, \DD \E, \DD \B, \DD \f, \dots) \dx \dv \\
\nonumber
&= \int \ff(x, E_1, E_2, B, \f, \DD E_1, \DD E_2, \DD B, \DD \f, \dots) \dx \dv \\
&= \int \bar{\ff}(x, E_1, E_2, B, \f, \DD E_1, \DD E_2, \DD B, \DD \f, \dots) \dx[1] \dv[1] \dv[2] ,
\end{align}
where in \eqref{eq:reduced_functional} we assumed \eqref{cordRed} and \eqref{fieldRed}. In the second equality of \eqref{eq:reduced_functional} the integrations over $x_2$, $x_3$ and $v_3$ are easily performed because the integrand is independent of these variables or it has been performed with an explicit dependence on $x_2$, $x_3$ and $v_3$ that makes the integrals converge. Any constant factors resulting from the integrations are absorbed into the definition of $\bar{\ff}$. The closure condition on $\{\ff,\fg\}$ amounts to the statement that given  any two functionals $\ff$, $\fg$ of the form of~\eqref{eq:reduced_functional}, their bracket  is again a reduced functional of this form. This follows from the fact that the bracket of two such functionals reduces~\eqref{eq:poisson_bracket_maxwell} to~\eqref{eq:poisson_bracket_maxwell_1d2v}, which of course is reduced.

That not all reductions of functionals have  closure can be seen by considering ones of the form  $\ff[\E,\f]$, i.e., ones for which dependence on $\B$ is absent. The bracket of two functions of this form gives the bracket of~\eqref{eq:poisson_bracket_maxwell} with the absence of the last term. Clearly this bracket depends on $\B$ and thus there is no closure. A consequence of this is that the bracket~\eqref{eq:poisson_bracket_maxwell} with the absence of the last term does not satisfy the Jacobi identity. We note, however, that adding a projector can remedy this, as was shown in~\cite{Chandre:2013}.

\subsection{Discrete deRham Complex in one Dimension}

Here, we consider the components of the electromagnetic fields separately and we have that $E_{1}$ is a one-form, $E_{2}$ is a zero-form and $B_{3}$ is again a one-form. We denote the semi-discrete fields by $D_{h}$, $E_{h}$ and $B_{h}$ respectively, and write
\begin{equation}\begin{aligned}
D_{h} (x, t) &= \sum \limits_{i=1}^{N_{1}} d_{i} (t) \, \Lambda^{1}_{i} (x) , \\
E_{h} (x, t) &= \sum \limits_{i=1}^{N_{0}} e_{i} (t) \, \Lambda^{0}_{i} (x) , \\
B_{h} (x, t) &= \sum \limits_{i=1}^{N_{1}} b_{i} (t) \, \Lambda^{1}_{i} (x) .
\end{aligned}\end{equation}
Next we introduce an equidistant grid in $x$ and denote the spline of degree $p$ with support starting at $x_i$ by $N_i^p$. We can express the derivative of $N_i^p$ as follows
\begin{equation}\label{eq:spline_derivative_p_pm1}
\frac{\dd}{\dd x} N_i^p(x) = \frac{1}{\Delta x} \left( N_i^{p-1}(x) - N_{i+1}^{p-1}(x) \right).
\end{equation}
In the Finite Element field solver, we represent $E_{h}$ by an expansion in splines of degree $p$ and $D_{h}$, $B_{h}$ by an expansion in splines of degree $p-1$, that is
\begin{equation}\label{eq:spline_expansion}\begin{aligned}
D_{h} (x, t) &= \sum_{i=1}^{N_{1}} d_i (t) \, N_i^{p-1}(x) , \\
E_{h} (x, t) &= \sum_{i=1}^{N_{0}} e_i (t) \, N_i^{p  }(x) , \\
B_{h} (x, t) &= \sum_{i=1}^{N_{1}} b_i (t) \, N_i^{p-1}(x) . \\
\end{aligned}\end{equation}

\subsection{Discrete Poisson Bracket}

The discrete Poisson bracket for this reduced system reads
\begingroup
\allowdisplaybreaks
\begin{align}\label{eq:eq:poisson_bracket_maxwell_matrix_1d2v}
\nonumber
\{\dff,\dfg\} & [\X, \V_{1}, \V_{2}, \cd, \ce, \cb] = \\
\nonumber
&= \fracp{\dff}{\X_{1}} \MM_{p}^{-1} \fracp{\dfg}{\V_{1}}
 - \fracp{\dfg}{\X_{1}} \MM_{p}^{-1} \fracp{\dff}{\V_{1}}
\\
\nonumber
&+ \bigg( \fracp{\dff}{\V_1}\bigg)^{\top} \MM_{p}^{-1} \MM_{q} \LaB^1(\X)^{\top} \MM_{1}^{-1} \bigg( \fracp{\dfg}{\cd} \bigg)
 - \bigg( \fracp{\dff}{\cd} \bigg)^{\top} \MM_{1}^{-1} \LaB^1(\X)        \MM_{q} \MM_{p}^{-1} \bigg( \fracp{\dfg}{\V_1}\bigg)
\\
\nonumber
&+ \bigg( \fracp{\dff}{\V_2}\bigg)^{\top} \MM_{p}^{-1} \MM_{q} \LaB^0(\X)^{\top} \MM_{0}^{-1} \bigg( \fracp{\dfg}{\ce} \bigg)
 - \bigg( \fracp{\dff}{\ce} \bigg)^{\top} \MM_{0}^{-1} \LaB^0(\X)        \MM_{q} \MM_{p}^{-1} \bigg( \fracp{\dfg}{\V_2}\bigg)
\\
\nonumber
&+ \bigg( \fracp{\dff}{\V_1}  \bigg)^{\top} \MM_{p}^{-1} \MM_{q} \BB(\X, \cb) \MM_{p}^{-1} \bigg( \fracp{\dfg}{\V_2} \bigg)
 - \bigg( \fracp{\dff}{\V_2}  \bigg)^{\top} \MM_{p}^{-1} \BB(\X, \cb) \MM_{q} \MM_{p}^{-1} \bigg( \fracp{\dfg}{\V_1} \bigg)
\\
&+ \bigg( \fracp{\dff}{\ce} \bigg)^{\top} \MM_{0}^{-1} \C^{\top} \bigg( \fracp{\dfg}{\cb} \bigg)
 - \bigg( \fracp{\dff}{\cb} \bigg)^{\top} \C        \MM_{0}^{-1} \bigg( \fracp{\dfg}{\ce} \bigg)
.
\end{align}
\endgroup
Here, we denote by $\MM_{p} = M_p$ and $\MM_{q} = M_q$ the $N_{p} \times N_{p}$ diagonal matrices holding the particle masses and charges, respectively. We denote by $\LaB^0(\X)$ the $N_p \times N_0$ matrix with generic term $\La^0_i(\cx_a)$, where $1 \leq a \leq N_p$ and $1 \leq i \leq N_0$, and by $\LaB^1(\X)$ the $N_p \times N_1$ matrix with generic term $\La^1_i(\cx_a)$, where $1 \leq a \leq N_p$ and $1 \leq i \leq N_1$. Further, $\BB(\X, \cb)$ denotes the $N_p \times N_p$ diagonal matrix with entries
\begin{align}
B_{h} (\cx_a, t) = \sum_{i=1}^{N_1} b_i(t) \, \La^{1}_{i} (\cx_a) .
\end{align}
The reduced bracket can be shown to satisfy the Jacobi identity by direct proof in full analogy to the proof for the full bracket shown in Section~\ref{sec:jacobi_identity}.
However, one can also follow along the lines of Section~\ref{sec:jacobi_identity_1d2v} in order to arrive at the same result.

\subsection{Discrete Hamiltonian and Equations of Motion}

The discrete Hamiltonian is given in terms of the reduced set of degrees of freedom $\uu = (\X, \V_{1}, \V_{2}, \cd, \ce, \cb)$ by
\begin{align}
\dham
&= \tfrac{1}{2} \, \V_{1}^{\top} \MM_{p} \V_{1}
 + \tfrac{1}{2} \, \V_{2}^{\top} \MM_{p} \V_{2}
 + \tfrac{1}{2} \, \cd^{\top} \MM_{1} \cd
 + \tfrac{1}{2} \, \ce^{\top} \MM_{0} \ce
 + \tfrac{1}{2} \, \cb^{\top} \MM_{1} \cb .
\end{align}
and the equations of motion are obtained as
\begin{equation}\label{eq:vlasov_maxwell_equations_of_motion_1d2v}
\begin{aligned}
\dot{\X}
&= \V_{1} ,
\\
\dot{\V}_{1}
&= \MM_p^{-1} \MM_q \big( \LaB^1 (\X) \cd + \BB (\X,\cb) \V_{2} \big) ,
\\
\dot{\V}_{2}
&= \MM_p^{-1} \MM_q \big( \LaB^0 (\X) \ce - \BB (\X,\cb) \V_{1} \big) ,
\\
\dot{\cd}
&= - \MM_{1}^{-1} \LaB^1 (\X)^\top \MM_q \V_{1} ,
\\
\dot{\ce}
&= \MM_{0}^{-1} \big( \C^\top \MM_{1} \cb (t) - \LaB^0 (\X)^\top \MM_q \V_{2} \big) ,
\\
\dot{\cb}
&= - \C \ce (t) ,
\end{aligned}
\end{equation}
which is seen to be in direct correspondence with~\eqref{eq:vlasov_1d2v}-\eqref{eq:maxwell_1d2v_last}.

\subsection{Hamiltonian Splitting}

The solution of the discrete equations of motion for $\dham_{D} + \dham_{E}$ at time $\Dt$ is
\begin{equation}\begin{aligned}
\MM_{p} \V_1 (\Dt) &= \MM_{p} \V_1 (0) + \Dt \, \MM_{q} \LaB^{1} (\X(0)) \, \cd (0) , \\
\MM_{p} \V_2 (\Dt) &= \MM_{p} \V_2 (0) + \Dt \, \MM_{q} \LaB^{0} (\X(0)) \, \ce (0) , \\
\cb  (\Dt) &= \cb (0) - \Dt \, \C \ce (0) .
\end{aligned}\end{equation}
The solution of the discrete equations of motion for $\dham_{B}$ is
\begin{equation}\begin{aligned}
\MM_{0} \ce (\Dt) &= \MM_{0} \ce(0) + \Dt \, \C^{\top} \MM_{1} \cb (0) .
\end{aligned}\end{equation}
The solution of the discrete equations of motion for $\dham_{p_1}$ is
\begin{equation}\begin{aligned}
\X   (\Dt) &= \X_1 (0) + \Dt \, \V_1 (0) , \\
\MM_{p} \V_2 (\Dt) &= \MM_{p} \V_2 (0) - \int \limits_{0}^{\Dt} \MM_{q} \BB (\X (t), \cb(0)) \V_{1} (0) \dt , \\
\MM_{1} \cd  (\Dt) &= \MM_{1} \cd  (0) - \int \limits_{0}^{\Dt} \LaB^1 (\X (t))^\top \MM_q \V_{1} (0) \dt ,
\end{aligned}\end{equation}
and for $\dham_{p_2}$ it is
\begin{equation}\begin{aligned}
\MM_{p} \V_1 (\Dt) &= \MM_{p} \V_1 (0) + \int \limits_{0}^{\Dt} \MM_{q} \BB (\X (0), \cb(0)) \V_{2} (0) \dt , \\
\MM_{0} \ce (\Dt) &= \MM_{0} \ce (0) - \int \limits_{0}^{\Dt} \LaB^1 (\X (0))^\top \MM_q \V_{2} (0) \dt ,
\end{aligned}\end{equation}
respectively.

\section{Numerical Experiments}\label{sec:numerics}

We have implemented the Hamiltonian splitting scheme as well as the Boris--Yee scheme from Appendix~\ref{app:boris-yee} as part of the SeLaLib library~\cite{selalib}. In this section, we present results for various test cases in 1d2v, comparing the conservation properties of the total energy and of the Casimirs for the two schemes. We simulate the electron distribution function in a neutralizing ion background. In all experiments, we have used splines of order three for the 0-forms. The particle loading was done using Sobol numbers and antithetic sampling (symmetric around the middle of the domain in $x$ and around the mean value of the Gaussian from which we are sampling in each velocity dimension). We sample uniformly in $x$ and from the Gaussians of the initial distribution in each velocity dimension.

\subsection{Weibel Instability}\label{sec:weibel}

We consider the Weibel instability studied in Weibel \cite{Weibel:1959} in the form simulated in \cite{Crouseilles:2015}. We study a reduced 1d2v model with a perturbation along $x_1$, a magnetic field along $x_3$ and electric fields along the $x_1$ and $x_2$ directions. Moreover, we assume that the distribution function is independent of $v_3$. The initial distribution and fields are of the form
\begin{align}
\f(x, \vb,t=0) &= \frac{1}{2\pi \sigma_1 \sigma_2} \exp \left(- \frac{1}{2}\left( \frac{v_1^2}{\sigma_1^2} + \frac{v_2^2}{\sigma_2^2} \right) \right) \left( 1 + \alpha \cos( kx)\right), \quad x \in [0,2\pi/k),\\
B_3(x,t=0) &= \beta \cos(kx),\\
E_2(x,t=0) &= 0,
\end{align}
and $E_1(x,t=0)$ is computed from Poisson's equation. 
In our simulations, we use the following choice of parameters, $\sigma_1 = 0.02/\sqrt{2}$, $\sigma_2 = \sqrt{12} \sigma_1$, $k=1.25$, $\alpha = 0$ and $\beta = -10^{-4}$. Note that these are the same parameters as in \cite{Crouseilles:2015} except for the fact that we sample from the Maxwellian without perturbation in $x_1$.

\begin{figure}
\centering
\includegraphics[width=.7\textwidth]{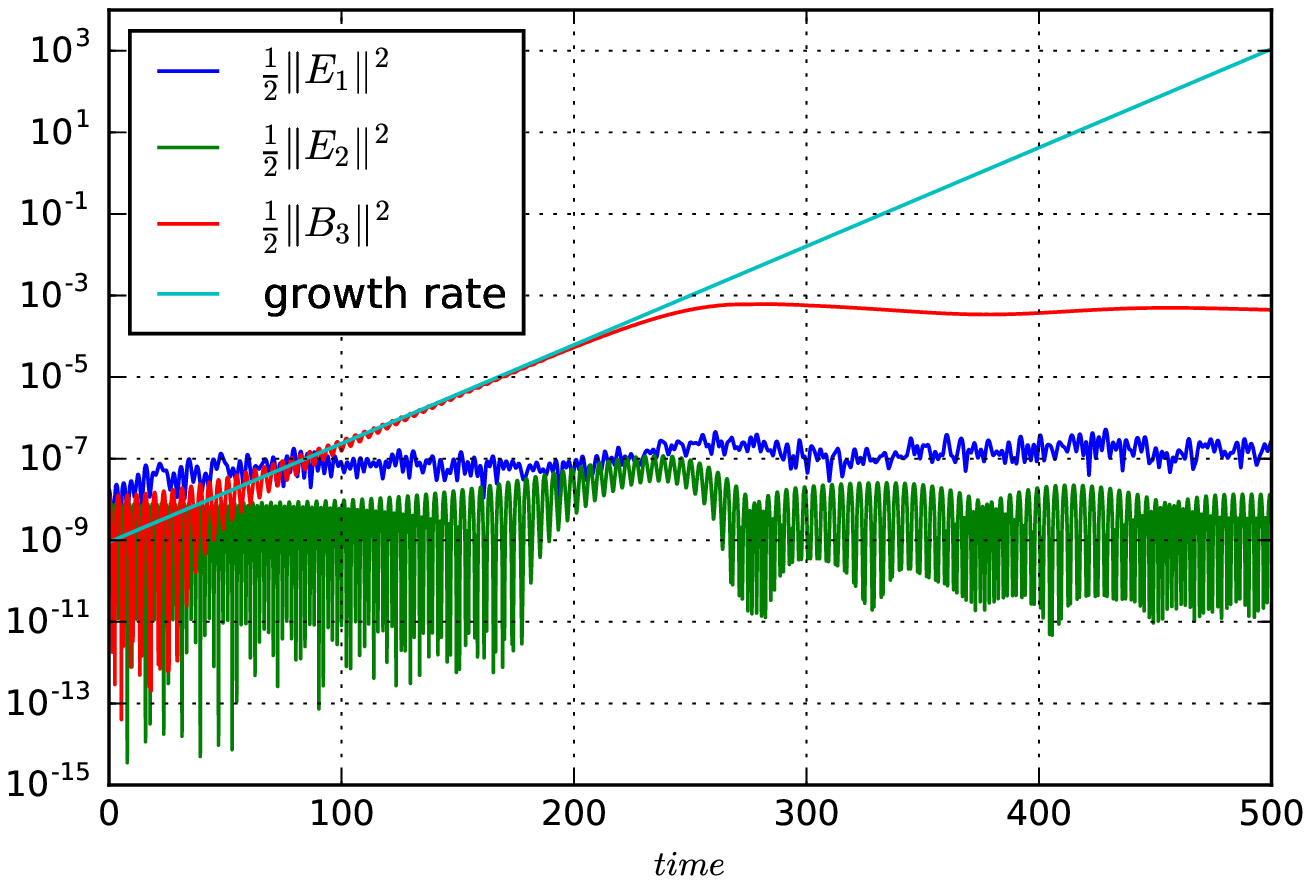}
\caption{Weibel instability: The two electric and the magnetic energies together with the analytic growth rate.}\label{fig:weibel_growth_rate}
\end{figure}

The dispersion relation from \cite{Weibel:1959} applied to our model reads
\begin{align}
D(\omega, k) = \omega^2 - k^2 + \left(\frac{\sigma_2}{\sigma_1}\right)^2-1 -  \left(\frac{\sigma_2}{\sigma_1}\right)^2 \phi\left( \frac{\omega}{\sigma_1 k}\right) \frac{\omega}{\sigma_1 k},
\end{align}
where $\phi(z) = \exp\left( -\frac{1}{2} z^2\right) \int_{-\mathrm{i} \infty}^z \exp\left(\frac{1}{2} \xi^2\right) \, \mathrm{d} \xi$. For our parameter choice, this gives a growth rate of 0.02784. In Fig.~\ref{fig:weibel_growth_rate}, we show the electric and magnetic energies together with the analytic growth rate. We see that the growth rate is verified in the numerical solution. This simulation was performed with 100,000 particles, 32 grid points, splines of degree 3 and 2 and $\Dt = 0.05$. Note that we have chosen a very large number of particles in order to obtain a solution of very high quality. In practice, the Weibel instability can also be simulated with much   fewer particles (cf.~Sec.~\ref{sec:numerics_momentum}).

In Table \ref{tab:weibel_conservation}, we show the conservation properties of our splitting with various orders of the splitting (cf.~Sec.~\ref{sec:composition}) and compare them also to the Boris--Yee scheme. The other numerical parameters are kept as before. 

We can see that Gauss' law is satisfied in each time step for the Hamiltonian splitting. This is a Casimir (cf. Sec.~\ref{sec:casimirs}) and therefore naturally conserved by the Hamiltonian splitting. On the other hand, this is not the case for the Boris--Yee scheme. Note that the numbers presented for Gauss' law in the table give the maximum difference of the coefficients for the electric field computed by Ampere's and Poisson's equation, respectively. 

\begin{figure}
\centering
\includegraphics[width=.7\textwidth]{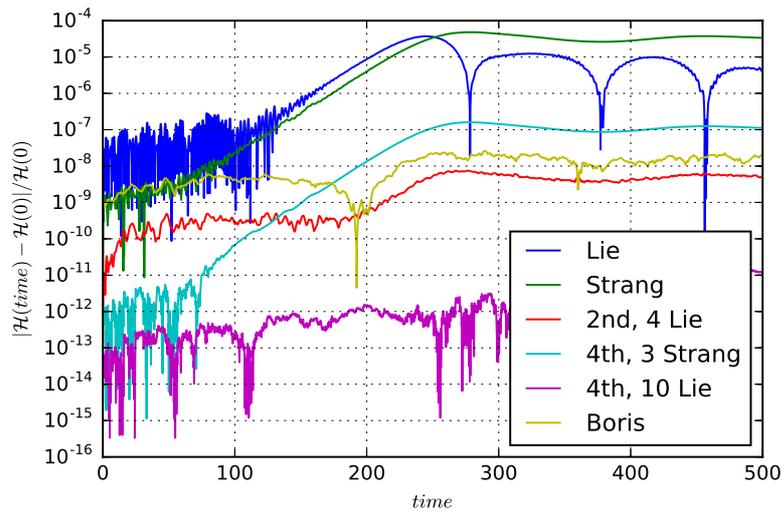}
\caption{Weibel instability: Difference of total energy and its initial value as a function of time for various integrators.}
\label{fig:weibel_energy}
\end{figure}

We can also see that the energy error improves with the order of the splitting; however, the Hamiltonian splitting method as well as the Boris--Yee scheme are not energy conserving. The time evolution of the total energy error is depicted in Fig.~\ref{fig:weibel_energy} for the various methods.

\begin{table}
\begin{center}
\begin{tabular}{|l|c|c|}
\hline 
Propagator & Total Energy & Gauss' law\\
\hline
Lie & 4.9E-7  & 8.7E-15  \\
Strang & 6.3E-7 & 1.5E-14 \\
2nd, 4 Lie & 9.8E-7 & 1.6E-14  \\
4th, 3 Strang & 2.1E-9 & 2.2E-14 \\
4th, 10 Lie & 2.1E-13 & 3.9E-14  \\
\hline
Boris & 3.4E-10 & 1.0E-4  \\
\hline
\end{tabular}
\caption{Weibel instability: Maximum error in the total energy and  Gauss' law until time 500 for simulation with various integrators: Lie--Trotter splitting from \eqref{eq:splitting_lie_trotter} (Lie), Strang splitting from \eqref{eq:splitting_strang} (Strang), second order splitting with 4 Lie parts defined in \eqref{eq:splitting_4lie} (2nd, 4 Lie), fourth order splitting with 3 Strang parts defined in \eqref{eq:splitting_3strang} (4th, 3 Strang) and 10 Lie parts defined in \eqref{eq:splitting_10lie} (4th, 10 Lie). }
\label{tab:weibel_conservation}
\end{center}
\end{table}

\subsection{Streaming Weibel Instability}\label{sec:streaming_weibel}

As a second test case, we consider the streaming Weibel instability. We study the same reduced model as in the previous section, but following~\cite{Califano:1997,Cheng:2014a} the initial distribution and fields are prescribed as
\begin{align}
\nonumber
\f(x, \vb,t=0) &= \frac{1}{2\pi \sigma^2} \exp \bigg(- \frac{ v_1^2 }{2\sigma^2} \bigg)  \bigg( \delta \exp\bigg(- \frac{(v_2-v_{0,1})^2}{2\sigma^2} \bigg)
\label{eq:streaming_weibel_dist} \\
&\hspace{12em} + (1-\delta) \exp\bigg(- \frac{(v_2-v_{0,2})^2}{2\sigma^2}\bigg)\bigg) , \\
B_3(x,t=0) &= \beta \sin(kx) , \\
E_2(x,t=0) &= 0 ,
\end{align}
and $E_1(x,t=0)$ is computed from Poisson's equation.  

\begin{figure}
\centering
\includegraphics[width=.7\textwidth]{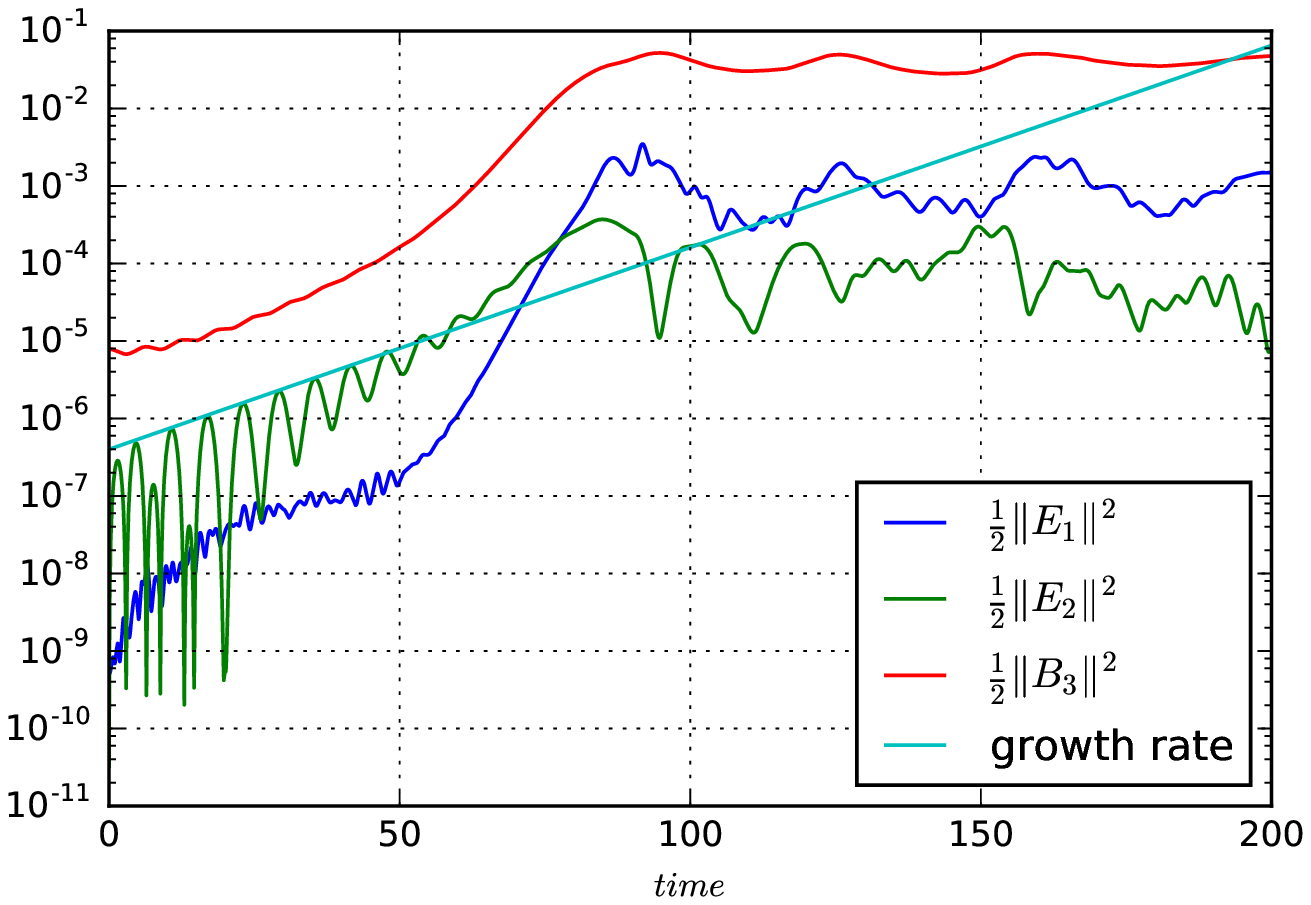}
\caption{Streaming Weibel instability: The two electric and the magnetic energies together with the analytic growth rate.}
\label{fig:sweibel2_growth_rate}
\end{figure}

We set the parameters to the following values $\sigma = 0.1/\sqrt{2}$, $k= 0.2$, $\beta = -10^{-3}$, $v_{0,1} = 0.5$, $v_{0,2} = -0.1$ and $\delta = 1/6$. The parameters are chosen as in the case 2 of \cite{Cheng:2014a}. The growth rate of energy of the second component of the electric field was determined to be 0.03 in \cite{Califano:1997}. In Fig.~\ref{fig:sweibel2_growth_rate}, we show the electric and magnetic energies together with the analytic growth rate. We see that the growth rate is verified in the numerical solution. This simulation was performed on the domain $x \in [0,2\pi/k)$ with 20,000,000 particles, 128 grid points, splines of degree 3 and 2 and $\Dt = 0.01$.
Observe that the energy of the~$E_{1}$ component of the electric field starts to increase at times earlier than in~\cite{Califano:1997}, which is caused by particle noise.

\begin{figure}
\centering
\includegraphics[width=.7\textwidth]{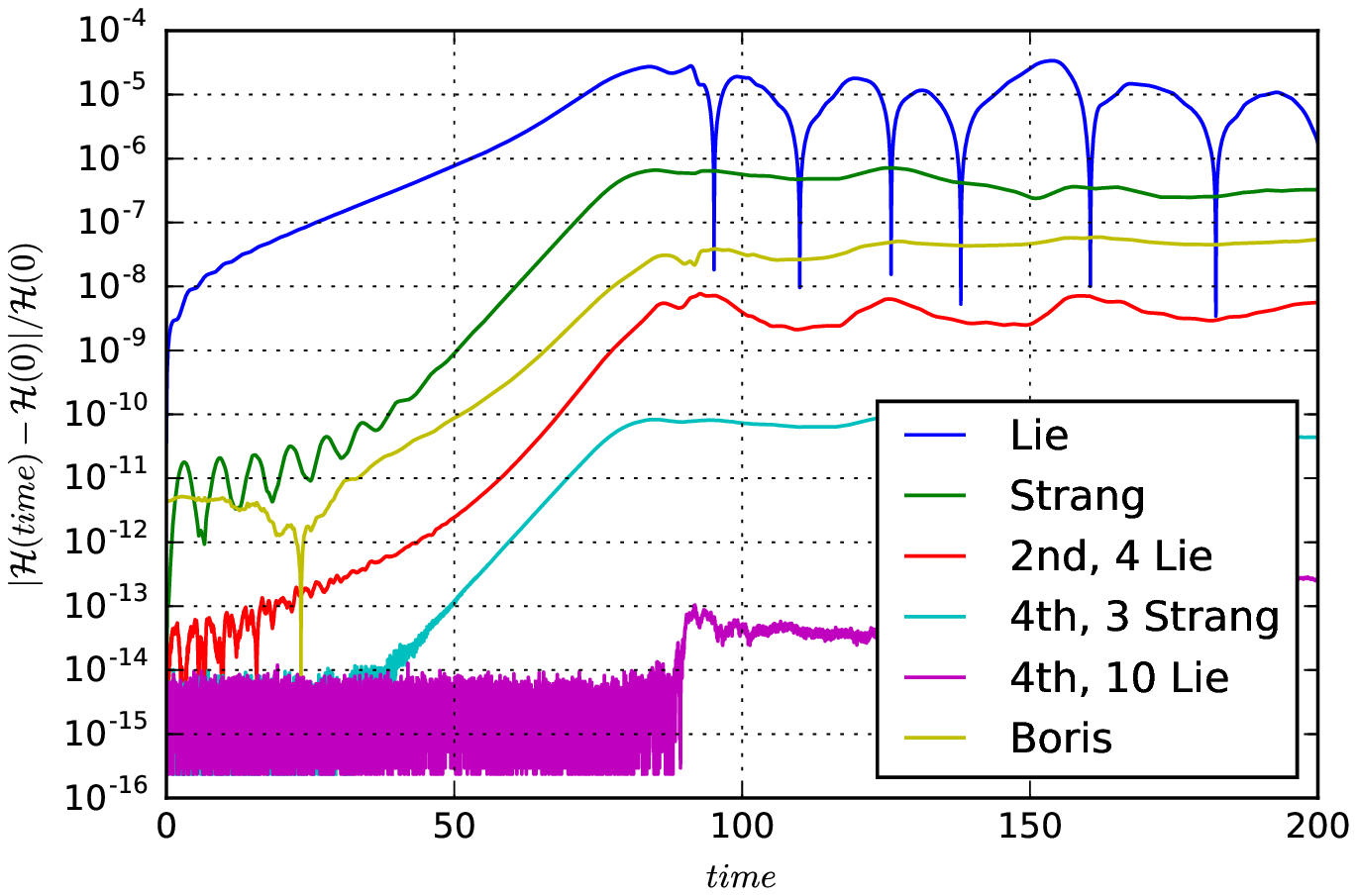}
\caption{Streaming Weibel instability: Difference of total energy and its initial value as a function of time for various integrators.}
\label{fig:sweibel2_energy}
\end{figure}

As for the Weibel instability, we compare the conservation properties in Table \ref{tab:sweibel_conservation} for various integrators. Again we see that the Hamiltonian splitting conserves Gauss' law as opposed to the Boris--Yee scheme. The energy conservation properties of the various schemes show approximately the same behaviour as in the previous test case (see also Fig.~\ref{fig:sweibel2_energy} for the time evolution of the energy error).

\begin{table}
\begin{center}
\begin{tabular}{|l|c|c|}
\hline 
Propagator & Total Energy & Gauss' law \\
\hline
Lie & 6.4E-5  & 3.1E-16 \\
Strang & 1.4E-6 & 3.7E-16 \\
2nd, 4 Lie & 1.5E-8 & 4.4E-16 \\
4th, 3 Strang & 1.7E-10 & 4.7E-16 \\
4th, 10 Lie & 5.7E-13 & 5.4E-16 \\
\hline
Boris & 1.1E-7 & 5.1E-4  \\
\hline
\end{tabular}
\caption{Streaming Weibel instability: Maximum error in the total energy and  Gauss' law until time 200 for simulation with various integrators.}
\label{tab:sweibel_conservation}
\end{center}
\end{table}

\subsection{Strong Landau Damping}

\begin{figure}[ht]
\centering
\includegraphics[width=.7\textwidth]{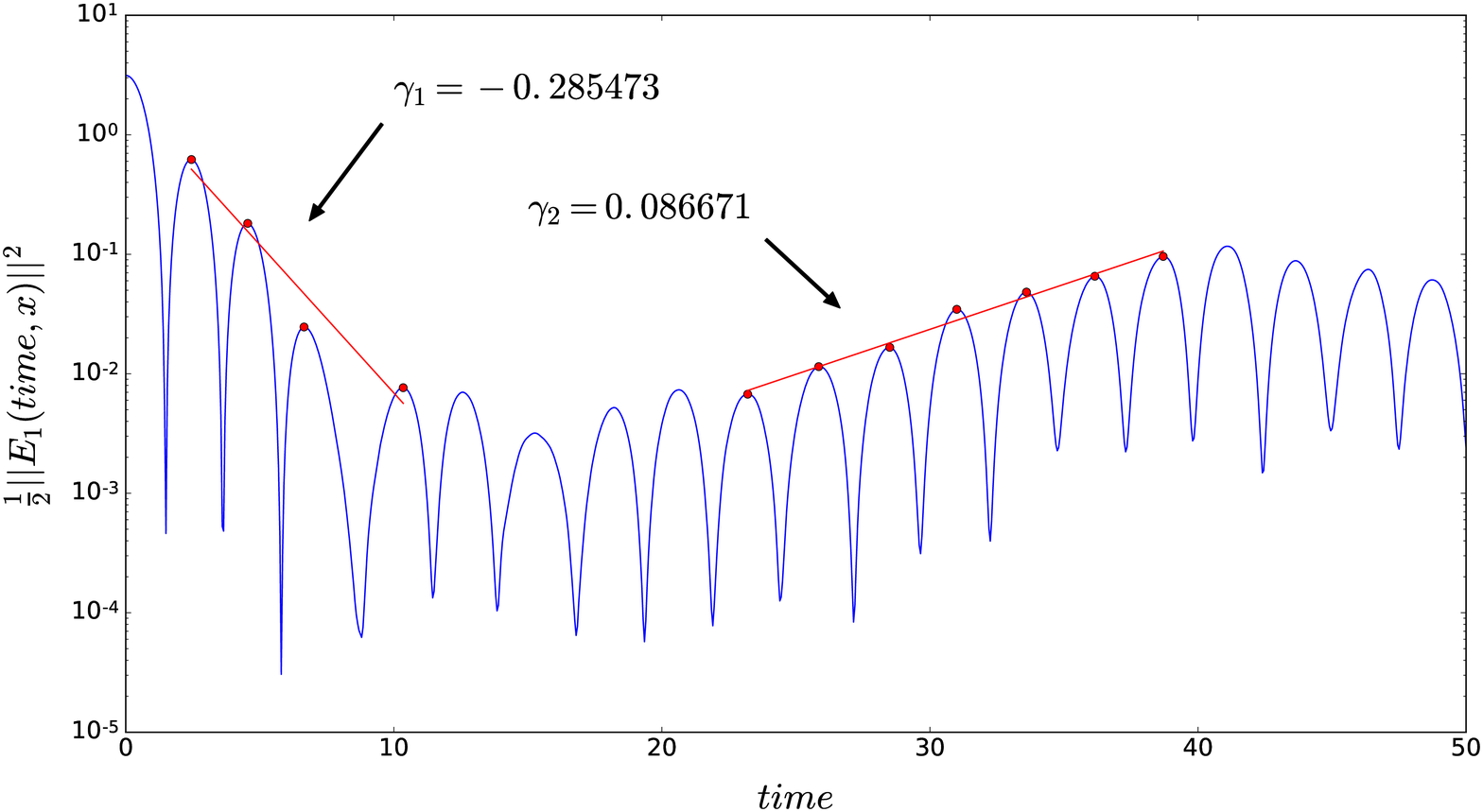}
\caption{Landau damping: Electric energy with fitted damping and growth rates.}
\label{fig:landau_eenergy}
\end{figure}

\begin{table}[ht]
\begin{center}
\vspace{1em}
\begin{tabular}{|l|c|c|c|}
\hline
Integrator												& $\gamma_{1}$	& $\gamma_{2}$ \\
\hline
GEMPIC													&  $-0.286$		& $+0.087$ \\
viVlasov1D \cite{KrausMajSonnendruecker:2015}			&  $-0.286$		& $+0.085$ \\
\hline
Cheng \& Knorr \cite{ChengKnorr:1976}					& $-0.281$		& $+0.084$ \\
Nakamura \& Yabe \cite{NakamuraYabe:1999}				& $-0.280$		& $+0.085$ \\
Ayuso \& Hajian \cite{deDios:2012b}						& $-0.292$		& $+0.086$ \\
Heath, Gamba, Morrison, Michler 	\cite{HeathGamba:2012}	& $-0.287$		& $+0.075$ \\
Cheng, Gamba, Morrison \cite{ChengGamba:2013}			& $-0.291$		& $+0.086$ \\
\hline
\end{tabular}
\caption{Damping and growth rates for strong Landau damping.}
\label{tab:landau_damping}
\end{center}
\end{table}

\begin{figure}[p]
\centering
\includegraphics[width=.7\textwidth]{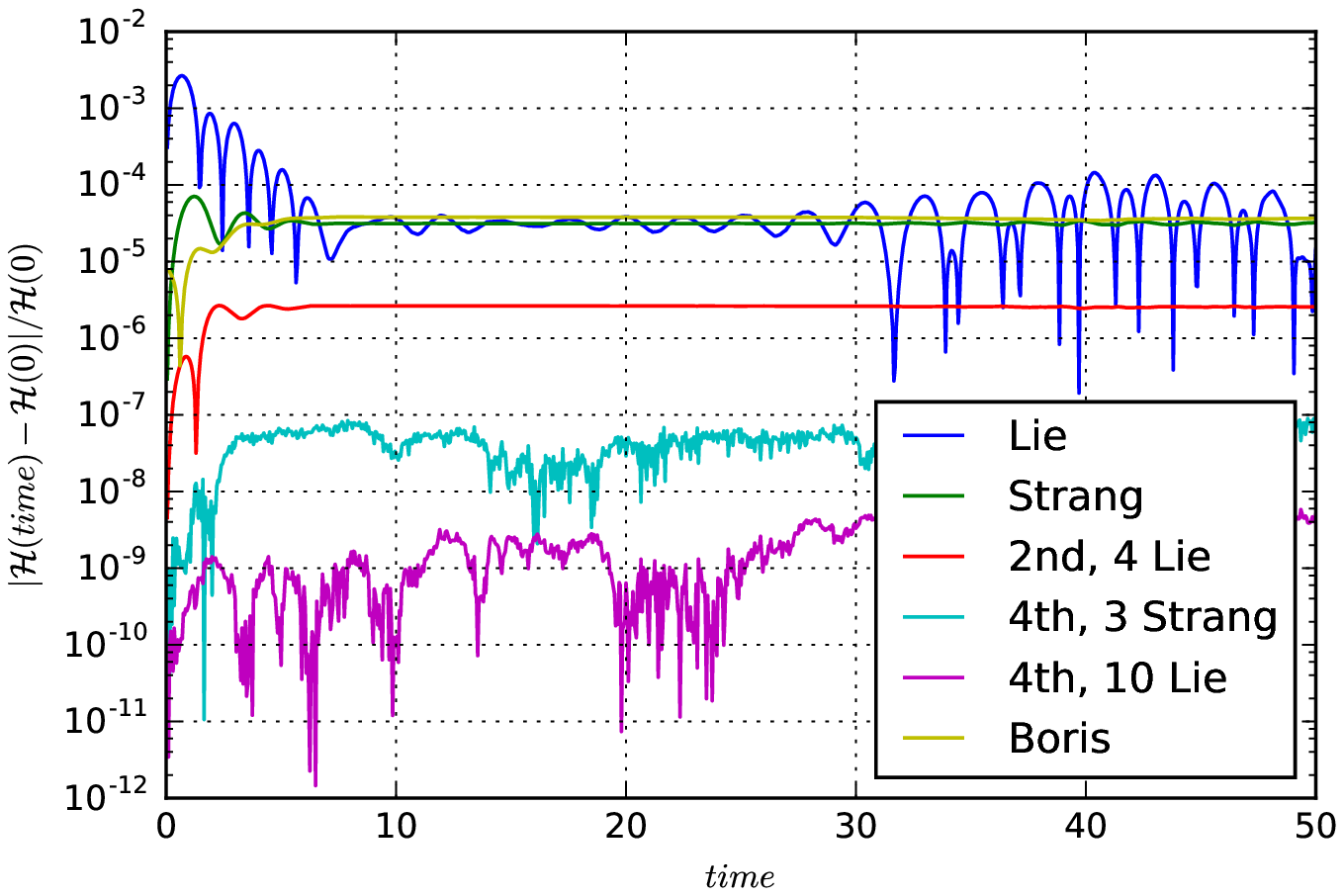}
\caption{Landau damping: Total energy error.}
\label{fig:landau_energy}
\end{figure}

Finally, we also study the electrostatic example of strong Landau damping with initial distribution
\begin{align}
\f (x,\vb) &= \exp \frac{1}{2 \pi \sigma^2}\left( - \tfrac{v_{1}^{2} + v_{2}^{2}}{2\sigma^{2}}\right) \, \big( 1 + \alpha \, \cos (kx) \big), \quad x \in [0, 2\pi / k ), \vb \in \R^{2}.
\end{align}
The physical parameters are chosen as $\sigma = 1$, $k = 0.5$, $\alpha = 0.5$ and the numerical parameters as $\De t = 0.05$, $n_{x} = 32$ and 100,000 particles.
The fields $B_3$ and $E_2$ are initialized to zero and remain zero over time. In this example, we essentially solve the Vlasov--Amp{\`e}re equation with results equivalent to the Vlasov--Poisson equations. In Fig.~\ref{fig:landau_eenergy} we show the time evolution of the electric energy associated with $E_1$. We have also fitted a damping and growth rate (using the marked local maxima in the plot). These are in good agreement with other codes (see Table~\ref{tab:landau_damping}). Again the energy conservation for the various method is visualized as a function of time in Fig.~\ref{fig:landau_energy}. And again we see that the fourth order methods give excellent energy conservation. 

\subsection{Backward Error Analysis}

\begin{figure}
\centering
\includegraphics[width=.7\textwidth]{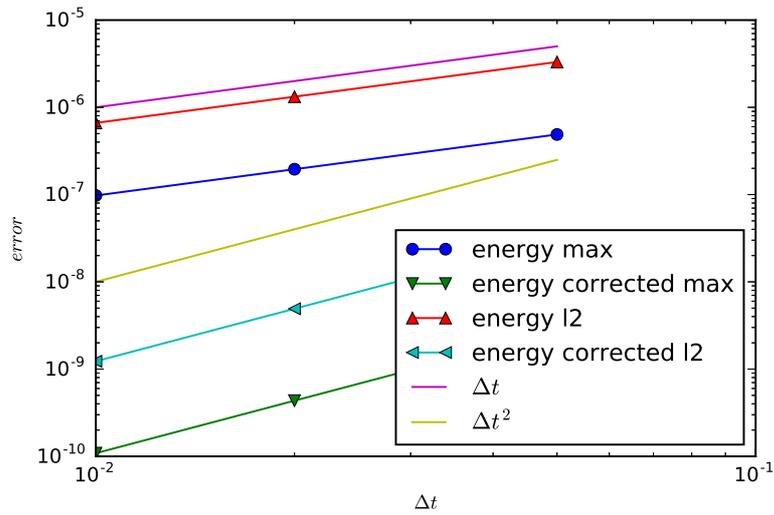}
\caption{Weibel instability: Error (maximum norm and $\ell_2$ norm) in the total energy for simulations with $\Dt = 0.01, 0.02, 0.05$.}
\label{fig:weibel_bea_convergence}
\end{figure}

\begin{figure}
\centering
\includegraphics[width=.7\textwidth]{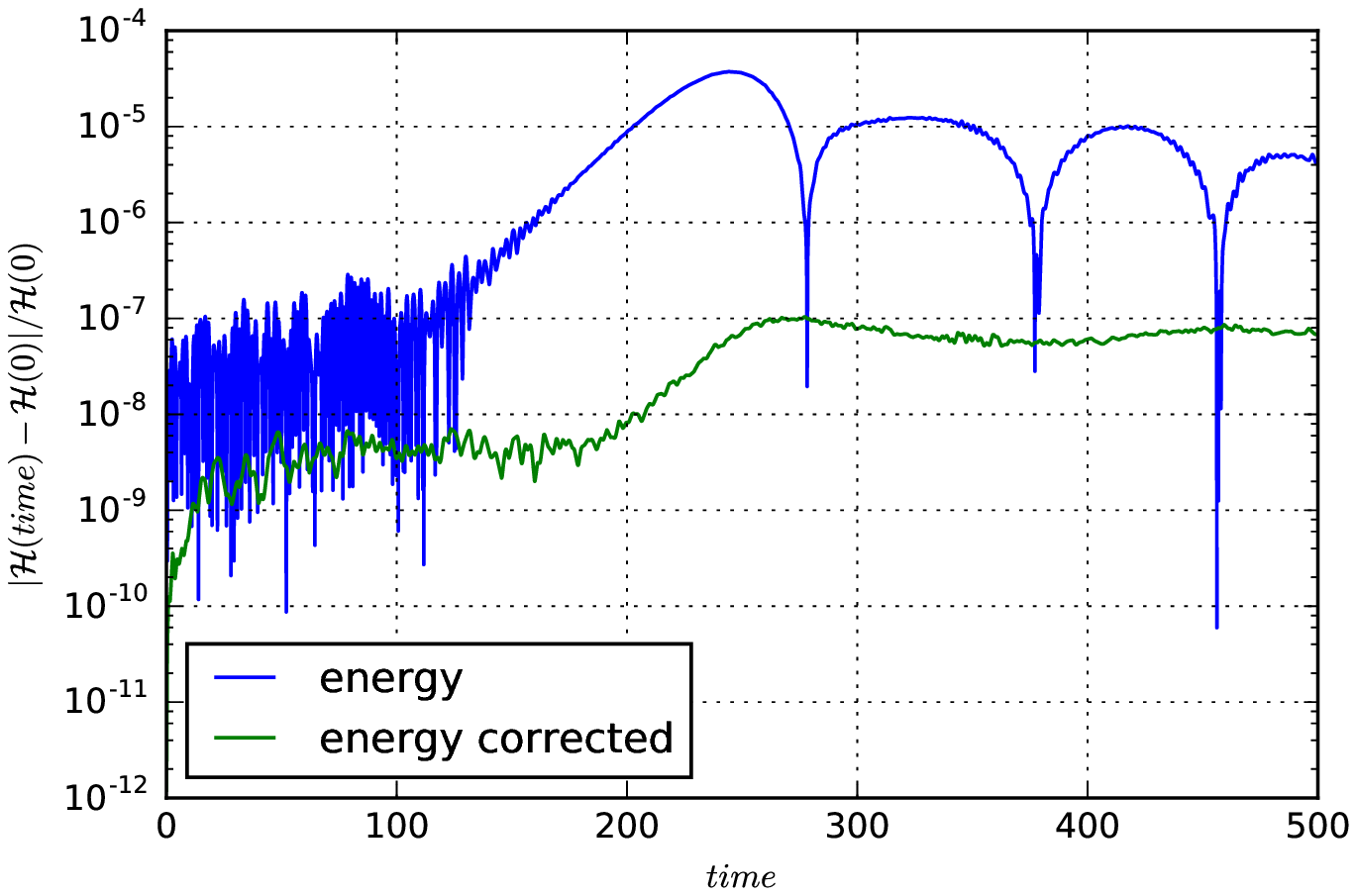}
\caption{Weibel instability: Energy and first-order corrected energy for simulation with $\Dt = 0.05$.}
\label{fig:weibel_bea_energy}
\end{figure}

For the Lie--Trotter splitting, the error in the Hamiltonian $\dham$ is of order $\Dt$. However, using backward error analysis (cf. Section~\ref{sec:splitting_bae}), modified Hamiltonians can be computed, which are preserved to higher order.
Accounting for first order corrections $\tilde{\dham}_{1}$, the error in the modified Hamiltonian,
\begin{align}
\tilde{\dham} = \dham + \Dt \tilde{\dham}_{1} + \order{\Dt^{2}} ,
\end{align}
is of order $(\Dt)^2$.
For the 1d2v example, this correction is obtained as
\begin{align}
\tilde{\dham}_{1} &= \dfrac{1}{2} \left[
    \{ \dham_{D  } + \dham_{E  } , \dham_{B  } \}
  + \{ \dham_{p_1} , \dham_{p_2} \}
  + \{ \dham_{D  } + \dham_{E  } + \dham_{B  } , \dham_{p_1} + \dham_{p_2} \}
\right] ,
\end{align}
with the various Poisson brackets computed as
\begin{align*}
\{ \dham_{D  } , \dham_{B}    \} &= 0 , \\
\{ \dham_{E  } , \dham_{B}    \} &=   \ce^{\top} \C^{\top} \cb , \\
\{ \dham_{p_1} , \dham_{p_2} \} &= \V_{1} M_q \BB (\cb,\X)^\top \, \V_{2} , \\
\{ \dham_{D  } , \dham_{p_1} \} &= - \cd^\top \LaB^1 (\X)^\top M_q \V_{1} , \\
\{ \dham_{D  } , \dham_{p_2} \} &= 0 , \\
\{ \dham_{E  } , \dham_{p_1} \} &= 0 , \\
\{ \dham_{E  } , \dham_{p_2} \} &= - \ce^\top \LaB^0 (\X)^\top M_q \V_{2} , \\
\{ \dham_{B  } , \dham_{p_1} \} &= 0 , \\
\{ \dham_{B  } , \dham_{p_2} \} &= 0 .
\end{align*}
In Fig.~\ref{fig:weibel_bea_convergence}, we show the maximum and $\ell_2$ error of the energy and the corrected energy for the Weibel instability test case with the parameters in Sec.~\ref{sec:weibel}. The simulations were performed with 100,000 particles, 32 grid points, splines of degree 3 and 2 and $\Dt = 0.01, 0.02, 0.05$.  We can see that the theoretical convergence rates are verified in the numerical experiments. 
Figure~\ref{fig:weibel_bea_energy} shows the energy as well as the modified energy for the Weibel instability test case simulated with a time step of $0.05$. 

\subsection{Momentum Conservation}\label{sec:numerics_momentum}

Finally, let us discuss conservation of momentum for our one species 1d2v equations. In this case, Eq.~\eqref{eq:dotPe} becomes
\begin{align}\label{eq:dPdt}
\frac{\dd P_{e,1,2}}{\dt} &= - \int E_{1,2}(x,t) \,\dd x,
\end{align}
For $P_{e,1}(t)$ from \eqref{eq:maxwell_1d2v_Ex} we get 
\begin{align}
\frac{\dd }{\dd t} \int E_1(x,t) \,\dd x = - \int j_1(x,t) \,\dd x.
\end{align}
Since in general $\int j_1(x,t) \,\dd x \neq 0$, momentum is not conserved. As a means of testing momentum conservation, Crouseilles et al. \cite{Crouseilles:2015} replaced \eqref{eq:maxwell_1d2v_Ex} by
\begin{align}
\fracp{E_1(x,t)}{t} &= - j_1 (x,t) + \int j_1(x,t) \,\dd x.  
\end{align}
Thus, for this artificial dynamics, momentum will be conserved if  $\int E_1(x,0) \,\dd x = 0$.
Instead, we do not modify the equations but check the validity of~\eqref{eq:dPdt}. We define a discrete version of~\eqref{eq:dPdt}, integrated over time, in the following way:
\begin{align}
\tilde{P}_{e,1}^n &= P_{e,1}^0 - \frac{\Dt}{2} \sum_{m=1}^n \left( \int D_h^{m-1}(x)\, \dd x + \int D_h^{m}(x)\, \dd x  \right),\\
\tilde{P}_{e,2}^n &= P_{e,2}^0 - \frac{\Dt}{2} \sum_{m=1}^n \left( \int E_h^{m-1}(x)\, \dd x + \int E_h^{m}(x)\, \dd x  \right)\label{eq:momentum2_discrete}.
\end{align} 

\begin{figure}
\centering
\includegraphics[width=.7\textwidth]{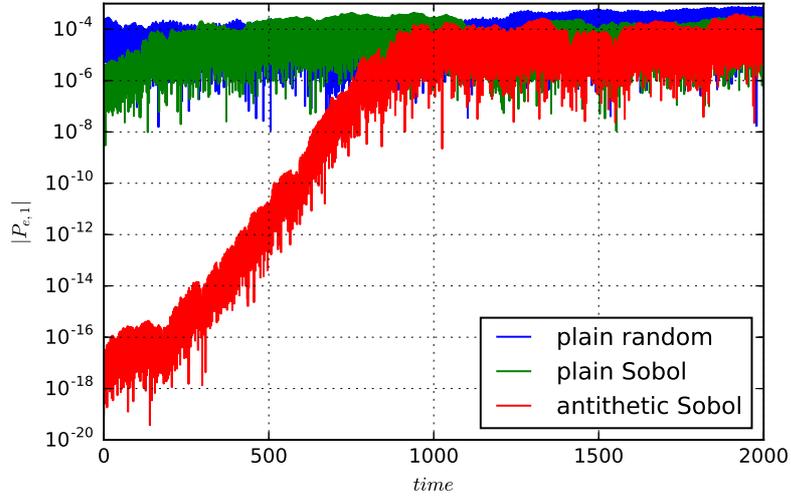}
\caption{Weibel instability: Error in the first component of the momentum for plain and antithetic Sobol sampling.}
\label{fig:weibel_momentum}
\end{figure}

\begin{figure}
\centering
\includegraphics[width=.7\textwidth]{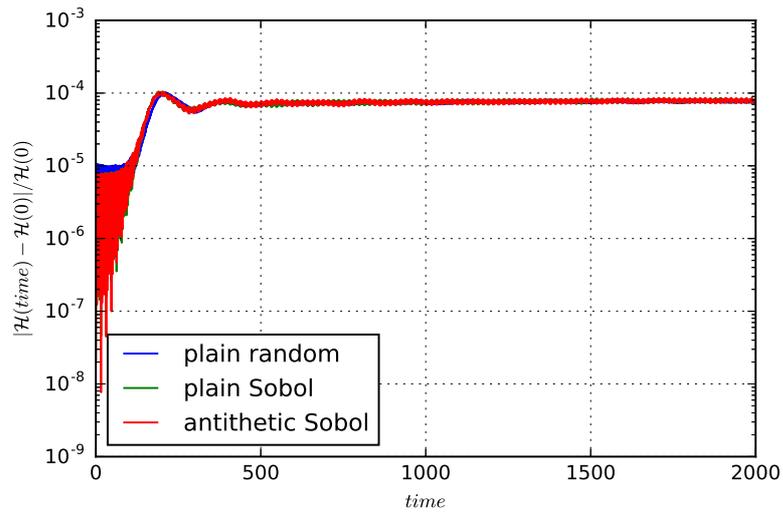}
\caption{Weibel instability: Total energy error for plain and antithetic Sobol sampling.}
\label{fig:weibel_energy_sampling}
\end{figure}

Our numerical scheme does not conserve momentum exactly. However, the error in momentum can be kept rather small during the linear phase of the simulation. Note that in all our examples, $\int j_{1,2}(x,t) \,\dd x = 0$. For a Gaussian initial distribution, the antithetic sampling ensures that  $\int j_{1,2}(x,t) \,\dd x = 0$ holds in the discrete sense. Figure \ref{fig:weibel_momentum} shows the momentum error in a simulation of the Weibel instability  as considered in Section \ref{sec:weibel}, but up to time 2,000 and with 25,600 particles sampled from pseudo-random numbers and Sobol numbers, for both plain and antithetical sampling. For the plain sampling, momentum error is a bit smaller for Sobol numbers compared to pseudo-random numbers during the linear phase. For the antithetic sampling, we can see that the momentum error is very small until time 200 (linear phase). However, when nonlinear effects start to dominate, the momentum error slowly increases until it has reached the same level as the momentum error for plain Sobol number sampling. The level depends on the number of particles (cf.~Table \ref{tab:weibel_momentum}). Note that the sampling does not seem to have an influence on the energy conservation as can be seen in Figure \ref{fig:weibel_energy_sampling} that compares  the energy error for the various sampling techniques. The curves show that the energy error is related to the increase in potential energy during the linear phase but does not further grow during the nonlinear phase.

\begin{table}[t]
\begin{center}
\begin{tabular}{|c|c|c|}
\hline
$N_p$ & $P_{e,1}$ & $P_{e,2}$ \\
\hline
320 & 1.92E-2 & 1.44E-2 \\
3,200 & 7.17E-3 & 8.96E-3 \\
25,600 & 4.45E-4 & 1.49E-3 \\
100,000 & 1.72E-4 & 9.02E-4 \\
\hline
\end{tabular}
\caption{Weibel instability: Maximum error in both components of the momentum for simulations until time 2,000 with various numbers of particles and 32 grid cells.}
\label{tab:weibel_momentum}
\end{center}
\end{table}

For the streaming Weibel instability, on the other hand, we have a sum of two Gaussians in the second component of the velocity.
Since in our sampling method we draw the particles based on Sobol quasi-random numbers, the fractions drawn from each of the Gaussians are not exactly given by $\delta$ and $1-\delta$ as in~\eqref{eq:streaming_weibel_dist}. Hence $\int j_{2}(x,t) \,\dd x$ is small but non-zero. Comparing the discrete momentum defined by Eq.~\eqref{eq:momentum2_discrete} with the discretization of its definition \eqref{eq:momentum},
\begin{align}
P_{e,2}^n = \sum_a v_{a,2}^n m_a w_a - \int D_h^n(x) B_h^n(x) \, \dd x,
\end{align}
we see a maximum deviation over time of $1.6 \times 10^{-16}$ for a simulation with  Strang splitting, $20,000,000$ particles, 128 grid points and $\Dt=0.01$. This shows that our discretization with a Strang splitting conserves \eqref{eq:dPdt} in the sense of \eqref{eq:momentum2_discrete} to very high accuracy.

\section{Summary}

In this work, a general framework for geometric Finite Element Particle-in-Cell methods for the Vlasov--Maxwell system was presented. The discretization proceeded in two steps. First, a semi-discretization of the noncanonical Poisson bracket was obtained, which preserves the Jacobi identity and important Casimir invariants, so that the resulting finite dimensional system is still Hamiltonian. Then, the system was discretized in time by Hamiltonian splitting methods, still retaining exact conservation of Casimirs, which in practice means exact conservation of Gauss' law and $\div B = 0$. Therefore the resulting method corresponds to one of the rare instances of a genuine Poisson integrator. Energy is not preserved exactly, but backward error analysis showed that the energy error does not depend on the degrees of freedom, the number of particles or the number of time steps. The favourable properties of the method were verified in various numerical experiments in 1d2v using splines as basis functions for the electromagnetic fields. One of the advantages of our approach is that conservation laws such as those for   energy and  charge  are not manufactured into the scheme ``by hand'' but follow automatically from preserving the underlying geometric structure of the equations.

The basic structure and implementation strategy of the code is very similar to existing Finite Element Particle-in-Cell methods for the Vlasov--Maxwell system. The main difference is the use of basis functions of mixed polynomial degree for the electromagnetic fields. The particle pusher is very similar to usual schemes, such as the Boris scheme, the only additional complexity being the exact computation of some line integrals.
The cost of the method is comparable to existing charge-conserving algorithms like the method of~\citeauthor{Villasenor.Buneman.1992.cpc} or the Boris-correction method. It is somewhat more expensive than non-charge-conserving methods, but such schemes are known to be prone to spurious instabilities that can lead to unphysical simulation results.
Even though only examples in 1d2v were shown, there are no conceptional differences or difficulties when going from one to two or to  three spatial dimensions. The building blocks of the code are identical in all three cases due to the tensor product structure of the Eulerian grid and the splitting in time. Details on the implementation of a three-dimensional version of the code as well as a comparison with existing methods will be disseminated in a subsequent publication.

The generality of the framework opens up several new paths for subsequent research. Instead of splines, other Finite Element spaces that form a deRham complex could be used, e.g., mimetic spectral elements or N\'ed\'elec elements for one-forms and Raviart--Thomas elements for two-forms. 
Further, it also should be possible to apply this approach to other systems like the gyrokinetic Vlasov--Maxwell system~\cite{Burby:2015, Burby:2017}, although in this case the necessity for new splitting schemes or other time integration strategies might arise.
Energy-preserving time stepping methods might provide an alternative to Hamiltonian splitting algorithms, where a suitable splitting cannot be easily found. This is a topic currently under investigation.
So is the treatment of the relativistic Vlasov--Maxwell system, which follows very closely along the lines of the non-relativistic system. This is a problem of interest in its own right, featuring an even larger set of invariants one should attempt to preserve in the discretization.
Another extension of the GEMPIC framework under development is the inclusion of non-ideal, non-Hamiltonian effects, most importantly collisions. An appropriate geometric description for such effects can be provided either microscopically by stochastic Hamiltonian processes or macroscopically by metriplectic brackets~\cite{Morrison:1986}.

\section*{Acknowledgments}

\noindent Useful discussions with Omar Maj and Marco Restelli are gratefully appreciated.
MK has received funding from the European Union's Horizon 2020 research and innovation programme under the Marie Sklodowska-Curie grant agreement No 708124.
PJM received support from the U.S. Dept.\ of Energy Contract \# DE-FG05-80ET-53088 and a Forschungspreis from the Humboldt Foundation. 
This work has been carried out within the framework of the EUROfusion Consortium and has received funding from the Euratom research and training programme 2014-2018 under grant agreement No 633053. The views and opinions expressed herein do not necessarily reflect those of the European Commission.
A part of this work was carried out using the HELIOS supercomputer system at Computational Simulation Centre of International Fusion Energy Research Centre (IFERC-CSC), Aomori, Japan, under the Broader Approach collaboration between Euratom and Japan, implemented by Fusion for Energy and QST.

\appendix

\section{The Boris--Yee scheme}
\label{app:boris-yee}

As an alternative discretization scheme, we consider a Boris--Yee scheme \cite{Boris:1970, Yee:1966} with our conforming Finite Elements. The scheme uses a time staggering working with the variables $\X^{n+1/2} = \X (t^n +\Dt/2)$, $\cd^{n+1/2} = \cd(t^n + \Dt/2)$, $\ce^{n+1/2} = \eb(t^n + \Dt /2)$, $\V^{n} = \V (t^n)$, and $\bb^{n} = \bb(t^n)$ in the $n$th time step $t^n = t^0 + n \Dt$.
The Hamiltonian at time $t^n$ is defined as
\begin{multline}
\dham
 = \tfrac{1}{2} \, (\V_1^{n})^{\top} M_{p} \V_1^{n}
 + \tfrac{1}{2} \, (\V_2^{n})^{\top} M_{p} \V_2^{n}
 + \tfrac{1}{2} \, (\cd^{n-1/2})^{\top} \MM_{1} \cd^{n+1/2} \\
 + \tfrac{1}{2} \, (\ce^{n-1/2})^{\top} \MM_{0} \ce^{n+1/2}
 + \tfrac{1}{2} \, (\cb^n)^{\top} \MM_{1} \cb^n .
\end{multline}
Given $\X^{n-1/2}$, $\db^{n-1/2}$, $\eb^{n-1/2}$, $\V^{n-1}$, $\bb^{n-1}$ the Vlasov--Maxwell system is propagated by the following time step:
\begin{enumerate}
\item Compute $\bb^{n}$ according to 
\begin{equation}
b_i^{n} = b_i^{n-1} - \frac{\Dt}{\Delta x} \left( e_i^{n-1/2} - e_{i-1}^{n-1/2} \right).
\end{equation}
 and $\bb^{n-1/2} = (\bb^{n-1}+\bb^{n})/2$.

\item Propagate $\vb^{n-1} \rightarrow \vb^{n}$ by equation 
\begin{align}
\label{eq:boris_ve1}
\vb_a^{-} &= \vb_a^{n-1} + \frac{\Dt}{2} \frac{q_s}{m_s} \Eb^{n-1/2}(\cx_a^{n-1/2}) ,
\\
\label{eq:boris_partB1}
\cv_{a,1}^{+} &= \frac{1-\alpha^2}{1+\alpha^2} \cv_{a,1}^{-} + \frac{2\alpha}{1+\alpha^2} \cv_{a,2}^{-},
\\
\label{eq:boris_partB2}
\cv_{a,2}^{+} &= \frac{-2\alpha}{1+\alpha^2} \cv_{a,1}^{-} +\frac{1-\alpha^2}{1+\alpha^2} \cv_{a,2}^{-},
\\
\label{eq:boris_ve2}
\vb_a^{n} &= \vb_a^{+} + \frac{\Dt}{2} \frac{q_s}{m_s} \Eb^{n-1/2}(\cx_a^{n-1/2}) ,
\end{align}
where $\alpha = \frac{q_a}{m_a} \frac{\Dt}{2} B^{n-1/2}(x^{n-1/2})$.

\item Propagate $\cx^{n-1/2} \rightarrow \cx^{n+1/2}$ by 
\begin{equation}\label{eq:boris_x}
\cx_a^{n+1/2} = \cx_a^{n-1/2} + \Dt \, \cv_{a,1}^{n}
\end{equation}
and accumulate $\jb_1^n$, $\jb_2^n$ by
\begin{align}
\jb_1^n &= \sum_{a=1}^{N_{p}} w_a \cv_{a,1}^{n} \Lambda^1 \big( ( \cx_a^{n-1/2} + \cx_a^{n+1/2} )/2 \big) , \\
\jb_2^n &= \sum_{a=1}^{N_{p}} w_a \cv_{a,2}^{n} \Lambda^0 \big( ( \cx_a^{n-1/2} + \cx_a^{n+1/2} )/2 \big) .
\end{align}

\item Compute $\db^{n+1/2}$ according to 
\begin{equation}\label{eq:boris_ex}
\MM_{1} \db^{n+1/2} = \MM_{1} \db^{n-1/2} - \Dt \, \jb_1^n ,
\end{equation}
 and  $\eb^{n+1/2}$ according to 
\begin{equation}\label{eq:boris_ey}
 \MM_{0} \eb^{n+1/2} = \MM_{0} \eb^{n-1/2} + \frac{\Dt}{\Delta x} \C^{T} \bb^{n} - \Dt \, \jb_2^n .
\end{equation}
\end{enumerate} 
For the initialization, we sample $\X^0$ and $\V^0$ from the initial sampling distribution, set $\ce^0$, $\cb^0$ from the given initial fields and solve Poisson's equation for $\cd^0$. Then, we compute $\X^{1/2}$, $\cd^{1/2}$ and $\ce^{1/2}$ from the corresponding equations of the Boris--Yee scheme for a half time step, using $\cb^0$, $\V^0$ instead of the unknown values at time $\Dt/4$.
Note that the error in this step is of order $(\Dt)^2$. But since we only introduce this error in the first time step, the overall scheme is still of order two.

\pagebreak

\bibliographystyle{plainnat}
\bibliography{gempic}

\begin{thebibliography}{102}
\providecommand{\natexlab}[1]{#1}
\providecommand{\url}[1]{\texttt{#1}}
\expandafter\ifx\csname urlstyle\endcsname\relax
  \providecommand{\doi}[1]{doi: #1}\else
  \providecommand{\doi}{doi: \begingroup \urlstyle{rm}\Url}\fi

\bibitem[sel()]{selalib}
{SeLaLib}.
\newblock \url{http://selalib.gforge.inria.fr/}.

\bibitem[Arnold et~al.(2006)Arnold, Falk, and
  Winther]{Arnold.Falk.Winther.2006.anum}
Douglas~N. Arnold, Richard~S. Falk, and Ragnar Winther.
\newblock Finite element exterior calculus, homological techniques, and
  applications.
\newblock \emph{Acta Numerica}, 15:\penalty0 1--155, 2006.
\newblock \doi{10.1017/S0962492906210018}.

\bibitem[Arnold et~al.(2010)Arnold, Falk, and
  Winther]{Arnold.Falk.Winther.2010.bams}
Douglas~N. Arnold, Richard~S. Falk, and Ragnar Winther.
\newblock Finite element exterior calculus: From hodge theory to numerical
  stability.
\newblock \emph{Bulletin of the American Mathematical Society}, 47:\penalty0
  281--354, 2010.
\newblock \doi{10.1090/S0273-0979-10-01278-4}.

\bibitem[Back and Sonnendr{\"u}cker(2014)]{Back:2014}
Aurore Back and Eric Sonnendr{\"u}cker.
\newblock {Finite Element Hodge for spline discrete differential forms.
  Application to the Vlasov--Poisson system}.
\newblock \emph{Applied Numerical Mathematics}, 79:\penalty0 124--136, 2014.
\newblock \doi{10.1016/j.apnum.2014.01.002}.

\bibitem[Baez and Muniain(1994)]{BaezMuniain:1994}
John Baez and Javier~P. Muniain.
\newblock \emph{Gauge Fields, Knots and Gravity}.
\newblock World Scientific, 1994.

\bibitem[Barthelm{\'e} and Parzani(2005)]{barthelme2005}
R{\'e}gine Barthelm{\'e} and C{\'e}line Parzani.
\newblock Numerical charge conservation in particle-in-cell codes.
\newblock In \emph{Numerical Methods for Hyperbolic and Kinetic Problems},
  pages 7--28. European Mathematical Society, 2005.
\newblock \doi{10.4171/012-1/1}.

\bibitem[Boffi(2006)]{Boffi.2006.csd}
Daniele Boffi.
\newblock {Compatible Discretizations for Eigenvalue Problems}.
\newblock In \emph{Compatible Spatial Discretizations}, pages 121--142.
  Springer New York, 2006.
\newblock \doi{10.1007/0-387-38034-5_6}.

\bibitem[Boffi(2010)]{boffi2010finite}
Daniele Boffi.
\newblock Finite element approximation of eigenvalue problems.
\newblock \emph{Acta Numerica}, 19:\penalty0 1--120, 2010.
\newblock \doi{10.1017/S0962492910000012}.

\bibitem[Boris(1970)]{Boris:1970}
Jay~P. Boris.
\newblock Relativistic plasma-simulation -- optimization of a hybrid code.
\newblock In \emph{Proceedings of the Fourth Conference on Numerical Simulation
  of Plasmas}, pages 3--67. Naval Research Laboratory, Washington DC, 1970.

\bibitem[Bossavit(1990)]{Bossavit:1990}
Alain Bossavit.
\newblock Differential geometry for the student of numerical methods in
  electromagnetism, 1990.
\newblock Lecture notes. Available online on
  \url{http://butler.cc.tut.fi/~bossavit/Books/DGSNME/DGSNME.pdf}.

\bibitem[Bossavit(2006)]{Bossavit:2005}
Alain Bossavit.
\newblock Applied differential geometry, 2006.
\newblock Lecture notes. Available online on
  \url{http://butler.cc.tut.fi/~bossavit/BackupICM/Compendium.html}.

\bibitem[Buffa and Perugia(2006)]{Buffa.Perugia.2006.sinum}
Annalisa Buffa and Ilaria Perugia.
\newblock {Discontinuous Galerkin approximation of the Maxwell eigenproblem}.
\newblock \emph{SIAM Journal on Numerical Analysis}, 44:\penalty0 2198--2226,
  2006.
\newblock \doi{10.1137/050636887}.

\bibitem[Buffa et~al.(2010)Buffa, Sangalli, and
  V{\'a}zquez]{buffa2010isogeometric}
Annalisa Buffa, Giancarlo Sangalli, and Rafael V{\'a}zquez.
\newblock Isogeometric analysis in electromagnetics: B-splines approximation.
\newblock \emph{Computer Methods in Applied Mechanics and Engineering},
  199\penalty0 (17):\penalty0 1143--1152, 2010.
\newblock \doi{10.1016/j.cma.2009.12.002}.

\bibitem[Buffa et~al.(2011)Buffa, Rivas, Sangalli, and V\'{a}zquez]{Buffa:2011}
Annalisa Buffa, Judith Rivas, Giancarlo Sangalli, and Rafael V\'{a}zquez.
\newblock Isogeometric discrete differential forms in three dimensions.
\newblock \emph{SIAM Journal on Numerical Analysis}, 49:\penalty0 818--844,
  2011.
\newblock \doi{10.1137/100786708}.

\bibitem[Burby(2017)]{Burby:2017}
Joshua~W. Burby.
\newblock Finite-dimensional collisionless kinetic theory.
\newblock \emph{Physics of Plasmas}, 24\penalty0 (3):\penalty0 032101, 2017.
\newblock \doi{10.1063/1.4976849}.

\bibitem[Burby et~al.(2015)Burby, Brizard, Morrison, and Qin]{Burby:2015}
Joshua~W. Burby, Alain~J. Brizard, Philip~J. Morrison, and Hong Qin.
\newblock {Hamiltonian Gyrokinetic Vlasov--Maxwell system}.
\newblock \emph{Physics Letters A}, 379\penalty0 (36):\penalty0 2073--2077,
  2015.
\newblock \doi{10.1016/j.physleta.2015.06.051}.

\bibitem[Califano et~al.(1997)Califano, Pegoraro, and Bulanov]{Califano:1997}
Francesco Califano, Francesco Pegoraro, and Sergei~V. Bulanov.
\newblock Spatial structure and time evolution of the weibel instability in
  collisionless inhomogeneous plasmas.
\newblock \emph{Phys. Rev. E}, 56:\penalty0 963--969, 1997.
\newblock \doi{10.1103/PhysRevE.56.963}.

\bibitem[Campos~Pinto et~al.(2014)Campos~Pinto, Jund, Salmon, and
  Sonnendr{\"u}cker]{Campos.Pinto.Jund.Salmon.Sonnendrucker.2014.crm}
Martin Campos~Pinto, S{\'e}bastien Jund, St{\'e}phanie Salmon, and Eric
  Sonnendr{\"u}cker.
\newblock {Charge conserving FEM-PIC schemes on general grids}.
\newblock \emph{Comptes Rendus Mecanique}, 342\penalty0 (10-11):\penalty0
  570--582, 2014.
\newblock \doi{10.1016/j.crme.2014.06.011}.

\bibitem[Caorsi et~al.(2000)Caorsi, Fernandes, and
  Raffetto]{Caorsi.Fernandes.Raffetto.2000.sinum}
Salvatore Caorsi, Paolo Fernandes, and Mirco Raffetto.
\newblock {On the convergence of Galerkin Finite Element Approximations of
  Electromagnetic Eigenproblems}.
\newblock \emph{SIAM Journal on Numerical Analysis}, 38:\penalty0 580--607,
  2000.
\newblock \doi{10.1137/S0036142999357506}.

\bibitem[Cendra et~al.(1998)Cendra, Holm, Hoyle, and Marsden]{Cendra:1998}
Hern{\'a}n Cendra, Darryl~D. Holm, Mark J.~W. Hoyle, and Jerrold~E. Marsden.
\newblock {The Maxwell--Vlasov Equations in Euler--Poincar{\'e} Form}.
\newblock \emph{Journal of Mathematical Physics}, 39:\penalty0 3138--3157,
  1998.
\newblock \doi{10.1063/1.532244}.

\bibitem[Chandre et~al.(2013)Chandre, Guillebon, Back, Tassi, and
  Morrison]{Chandre:2013}
Cristel Chandre, Lo{\"i}c Guillebon, Aurore Back, Emanuele Tassi, and Philip~J.
  Morrison.
\newblock On the use of projectors for {Hamiltonian} systems and their
  relationship with {Dirac} brackets.
\newblock \emph{Journal of Physics A: Mathematical and Theoretical},
  46:\penalty0 125203, 2013.
\newblock \doi{10.1088/1751-8113/46/12/125203}.

\bibitem[Channell and Scovel(1990)]{Channell:1990}
Paul~J. Channell and Clint Scovel.
\newblock {Symplectic integration of Hamiltonian systems}.
\newblock \emph{Nonlinearity}, 3\penalty0 (2):\penalty0 231 -- 259, 1990.
\newblock \doi{10.1088/0951-7715/3/2/001}.

\bibitem[Cheng and Knorr(1976)]{ChengKnorr:1976}
Chio-Zong Cheng and Georg Knorr.
\newblock The integration of the {Vlasov} equation in configuration space.
\newblock \emph{Journal of Computational Physics}, 22:\penalty0 330--351, 1976.
\newblock \doi{10.1016/0021-9991(76)90053-X}.

\bibitem[Cheng et~al.(2013)Cheng, Gamba, and Morrison]{ChengGamba:2013}
Yingda Cheng, Irene~M. Gamba, and Philip~J. Morrison.
\newblock Study of conservation and recurrence of {Runge}--{Kutta}
  discontinuous {Galerkin} schemes for {Vlasov}--{Poisson} systems.
\newblock \emph{Journal of Scientific Computing}, 56:\penalty0 319--349, 2013.
\newblock \doi{10.1007/s10915-012-9680-x}.

\bibitem[Cheng et~al.(2014{\natexlab{a}})Cheng, Christlieb, and
  Zhong]{Cheng:2014b}
Yingda Cheng, Andrew~J. Christlieb, and Xinghui Zhong.
\newblock Energy-conserving discontinuous {Galerkin} methods for the
  {Vlasov}--{Amp{\`e}re} system.
\newblock \emph{Journal of Computational Physics}, 256:\penalty0 630--655,
  2014{\natexlab{a}}.
\newblock \doi{10.1016/j.jcp.2013.09.013}.

\bibitem[Cheng et~al.(2014{\natexlab{b}})Cheng, Christlieb, and
  Zhong]{Cheng:2014c}
Yingda Cheng, Andrew~J. Christlieb, and Xinghui Zhong.
\newblock Energy-conserving discontinuous {Galerkin} methods for the
  {Vlasov}--{Maxwell} system.
\newblock \emph{Journal of Computational Physics}, 279:\penalty0 145--173,
  2014{\natexlab{b}}.
\newblock \doi{10.1016/j.jcp.2014.08.041}.

\bibitem[Cheng et~al.(2014{\natexlab{c}})Cheng, Gamba, Li, and
  Morrison]{Cheng:2014a}
Yingda Cheng, Irene~M. Gamba, Fengyan Li, and Philip~J. Morrison.
\newblock Discontinuous {Galerkin} methods for the {Vlasov}--{Maxwell}
  equations.
\newblock \emph{SIAM Journal on Numerical Analysis}, 52\penalty0 (2):\penalty0
  1017--1049, 2014{\natexlab{c}}.
\newblock \doi{10.1137/130915091}.

\bibitem[Christiansen et~al.(2011)Christiansen, Munthe-Kaas, and
  Owren]{Christiansen:2011}
Snorre~H. Christiansen, Hans~Z. Munthe-Kaas, and Brynjulf Owren.
\newblock Topics in structure-preserving discretization.
\newblock \emph{Acta Numerica}, 20:\penalty0 1--119, 2011.
\newblock \doi{10.1017/S096249291100002X}.

\bibitem[Crouseilles et~al.(2014)Crouseilles, Navaro, and
  Sonnendr\"ucker]{crouseilles2014}
Nicolas Crouseilles, Pierre Navaro, and Eric Sonnendr\"ucker.
\newblock Charge conserving grid based methods for the {Vlasov}--{Maxwell}
  equations.
\newblock \emph{Comptes Rendus Mecanique}, 342\penalty0 (10-11):\penalty0
  636--646, 2014.
\newblock \doi{10.1016/j.crme.2014.06.012}.

\bibitem[Crouseilles et~al.(2015)Crouseilles, Einkemmer, and
  Faou]{Crouseilles:2015}
Nicolas Crouseilles, Lukas Einkemmer, and Erwan Faou.
\newblock {Hamiltonian} splitting for the {Vlasov--Maxwell} equations.
\newblock \emph{Journal of Computational Physics}, 283:\penalty0 224--240,
  2015.
\newblock \doi{10.1016/j.jcp.2014.11.029}.

\bibitem[Darling(1994)]{Darling:1994}
Richard W.~R. Darling.
\newblock \emph{Differential Forms and Connections}.
\newblock Cambridge University Press, 1994.

\bibitem[de~Dios and Hajian(2012)]{deDios:2012b}
Blanca~Ayuso de~Dios and Soheil Hajian.
\newblock High order and energy preserving discontinuous {Galerkin} methods for
  the {Vlasov}--{Poisson} system.
\newblock \href{https://arxiv.org/abs/1209.4025}{arXiv:1209.4025}, 2012.

\bibitem[de~Dios et~al.(2011)de~Dios, Carrillo, and Shu]{deDios:2011}
Blanca~Ayuso de~Dios, Jos{\'e} Carrillo, and Chi-Wang Shu.
\newblock {Discontinuous Galerkin methods for the one--dimensional
  Vlasov--Poisson system}.
\newblock \emph{Kinetic and Related Models}, 4:\penalty0 955--989, 2011.
\newblock \doi{10.3934/krm.2011.4.955}.

\bibitem[de~Dios et~al.(2012)de~Dios, Carrillo, and Shu]{deDios:2012a}
Blanca~Ayuso de~Dios, Jos{\'e}~A Carrillo, and Chi-Wang Shu.
\newblock {Discontinuous Galerkin Methods For The Multi--Dimensional
  Vlasov--Poisson Problem}.
\newblock \emph{Mathematical Models and Methods in Applied Sciences},
  22:\penalty0 1250042, 2012.
\newblock \doi{10.1142/S021820251250042X}.

\bibitem[Desbrun et~al.(2008)Desbrun, Kanso, and Tong]{Desbrun:2008}
Mathieu Desbrun, Eva Kanso, and Yiying Tong.
\newblock \emph{Discrete Differential Geometry}, chapter Discrete Differential
  Forms for Computational Modeling, pages 287--324.
\newblock Birkh{\"a}user Basel, 2008.
\newblock ISBN 978-3-7643-8621-4.
\newblock \doi{10.1007/978-3-7643-8621-4_16}.

\bibitem[Dray(2014)]{Dray:2014}
Tevian Dray.
\newblock \emph{Differential Forms and the Geometry of General Relativity}.
\newblock CRC Press, 2014.

\bibitem[Eastwood(1991)]{eastwood1991virtual}
James~W. Eastwood.
\newblock The virtual particle electromagnetic particle-mesh method.
\newblock \emph{Computer Physics Communications}, 64\penalty0 (2):\penalty0
  252--266, 1991.
\newblock \doi{10.1016/0010-4655(91)90036-K}.

\bibitem[Esirkepov(2001)]{esirkepov2001exact}
Timur~Zh. Esirkepov.
\newblock Exact charge conservation scheme for particle-in-cell simulation with
  an arbitrary form-factor.
\newblock \emph{Computer Physics Communications}, 135\penalty0 (2):\penalty0
  144--153, 2001.
\newblock \doi{10.1016/S0010-4655(00)00228-9}.

\bibitem[Evstatiev and Shadwick(2013)]{Evstatiev:2013}
Evstati~G. Evstatiev and Bradley~A. Shadwick.
\newblock {Variational formulation of particle algorithms for kinetic plasma
  simulations}.
\newblock \emph{Journal of Computational Physics}, 245:\penalty0 376--398,
  2013.
\newblock \doi{10.1016/j.jcp.2013.03.006}.

\bibitem[Gerritsma(2012)]{Gerritsma:2012}
Marc Gerritsma.
\newblock {An Introduction to a Compatible Spectral Discretization Method}.
\newblock \emph{Mechanics of Advanced Materials and Structures}, 19:\penalty0
  48--67, 2012.
\newblock \doi{10.1080/15376494.2011.572237}.

\bibitem[Hairer et~al.(2006)Hairer, Lubich, and
  Wanner]{HairerLubichWanner:2006}
Ernst Hairer, Christian Lubich, and Gerhard Wanner.
\newblock \emph{Geometric Numerical Integration}.
\newblock Springer, 2006.

\bibitem[He et~al.(2015)He, Qin, Sun, Xiao, Zhang, and Liu]{He:2015}
Yang He, Hong Qin, Yajuan Sun, Jianyuan Xiao, Ruili Zhang, and Jian Liu.
\newblock {Hamiltonian integration methods for Vlasov--Maxwell equations}.
\newblock \emph{Physics of Plasmas}, 22:\penalty0 124503, 2015.
\newblock \doi{10.1063/1.4938034}.

\bibitem[He et~al.(2016)He, Sun, Qin, and Liu]{he2016}
Yang He, Yajuan Sun, Hong Qin, and Jian Liu.
\newblock {Hamiltonian} particle-in-cell methods for {Vlasov}--{Maxwell}
  equations.
\newblock \emph{Physics of Plasmas}, 23\penalty0 (9):\penalty0 092108, 2016.
\newblock \doi{10.1063/1.4962573}.

\bibitem[Heath et~al.(2012)Heath, Gamba, Morrison, and
  Michler]{HeathGamba:2012}
Ross~E. Heath, Irene~M. Gamba, Philip~J. Morrison, and Christian Michler.
\newblock A discontinuous {Galerkin} method for the {Vlasov}--{Poisson} system.
\newblock \emph{Journal of Computational Physics}, 231:\penalty0 1140--1174,
  2012.
\newblock \doi{10.1016/j.jcp.2011.09.020}.

\bibitem[Hesthaven and Warburton(2004)]{Hesthaven.Warburton.2004.phil}
Jan~S. Hesthaven and Tim Warburton.
\newblock {High-order nodal discontinuous Galerkin methods for the Maxwell
  eigenvalue problem}.
\newblock \emph{Philosophical Transactions of the Royal Society A:
  Mathematical, Physical and Engineering Sciences}, 362\penalty0
  (1816):\penalty0 493--524, 2004.
\newblock \doi{10.1098/rsta.2003.1332}.

\bibitem[Hirani(2003)]{Hirani:2003}
Anil~N. Hirani.
\newblock \emph{Discrete Exterior Calculus}.
\newblock PhD thesis, California Institute of Technology, 2003.
\newblock URL \url{http://resolver.caltech.edu/CaltechETD:etd-05202003-095403}.

\bibitem[Holloway(1996)]{Holloway:1996}
James~Paul Holloway.
\newblock {On Numerical Methods for Hamiltonian PDEs and a Collocation Method
  for the Vlasov--Maxwell Equations}.
\newblock \emph{Journal of Computational Physics}, 129\penalty0 (1):\penalty0
  121--133, 1996.
\newblock \doi{10.1006/jcph.1996.0238}.

\bibitem[Kraus(2013)]{Kraus:2013:thesis}
Michael Kraus.
\newblock \emph{Variational Integrators in Plasma Physics}.
\newblock PhD thesis, Technische Universit{\"a}t M{\"u}nchen, 2013.
\newblock \href{https://arxiv.org/abs/1307.5665}{arXiv:1307.5665}.

\bibitem[Kraus et~al.()Kraus, Maj, and
  Sonnendr\"ucker]{KrausMajSonnendruecker:2015}
Michael Kraus, Omar Maj, and Eric Sonnendr\"ucker.
\newblock {Variational Integrators for the Vlasov--Poisson System}.
\newblock In preparation.

\bibitem[Kreeft et~al.(2011)Kreeft, Palha, and Gerritsma]{kreeft2011}
Jasper Kreeft, Artur Palha, and Marc Gerritsma.
\newblock Mimetic framework on curvilinear quadrilaterals of arbitrary order.
\newblock \href{https://arxiv.org/abs/1111.4304}{arXiv:1111.4304}, 2011.

\bibitem[Langdon(1992)]{langdon1992}
A.~Bruce Langdon.
\newblock On enforcing {G}auss' law in electromagnetic particle-in-cell codes.
\newblock \emph{Computer Physics Communications}, 70:\penalty0 447--450, 1992.
\newblock \doi{10.1016/0010-4655(92)90105-8}.

\bibitem[Lewis(1970)]{Lewis:1970}
H.~Ralph Lewis.
\newblock Energy-conserving numerical approximations for {Vlasov} plasmas.
\newblock \emph{Journal of Computational Physics}, 6\penalty0 (1):\penalty0
  136--141, 1970.
\newblock \doi{10.1016/0021-9991(70)90012-4}.

\bibitem[Lewis(1972)]{Lewis:1972}
H.~Ralph Lewis.
\newblock {Variational algorithms for numerical simulation of collisionless
  plasma with point particles including electromagnetic interactions}.
\newblock \emph{Journal of Computational Physics}, 10\penalty0 (3):\penalty0
  400--419, 1972.
\newblock \doi{10.1016/0021-9991(72)90044-7}.

\bibitem[Low(1958)]{Low:1958}
Francis~E. Low.
\newblock {A Lagrangian Formulation of the Boltzmann-Vlasov Equation for
  Plasmas}.
\newblock In \emph{Proceedings of the Royal Society of London. Series A,
  Mathematical and Physical Sciences}, volume 248, pages 282--287, 1958.
\newblock \doi{10.1098/rspa.1958.0244}.

\bibitem[Madaule et~al.(2014)Madaule, Restelli, and
  Sonnendr{\"u}cker]{Madaule:2014}
{\'E}ric Madaule, Marco Restelli, and Eric Sonnendr{\"u}cker.
\newblock Energy conserving discontinuous {Galerkin} spectral element method
  for the {Vlasov}--{Poisson} system.
\newblock \emph{Journal of Computational Physics}, 279:\penalty0 261--288,
  2014.
\newblock \doi{10.1016/j.jcp.2014.09.010}.

\bibitem[Marder(1987)]{marder1987}
Barry Marder.
\newblock A method for incorporating {G}auss's law into electromagnetic {PIC}
  codes.
\newblock \emph{Journal of Computational Physics}, 68:\penalty0 48--55, 1987.
\newblock \doi{10.1016/0021-9991(87)90043-X}.

\bibitem[Marsden and Weinstein(1982)]{Marsden:1982}
Jerrold~E. Marsden and Alan Weinstein.
\newblock {The Hamiltonian structure of the Maxwell--Vlasov equations}.
\newblock \emph{Physica D: Nonlinear Phenomena}, 4\penalty0 (3):\penalty0
  394--406, 1982.
\newblock \doi{10.1016/0167-2789(82)90043-4}.

\bibitem[McLachlan(1995)]{McLachlan:1995}
Robert~I. McLachlan.
\newblock {On the Numerical Integration of Ordinary Differential Equations by
  Symmetric Composition Methods}.
\newblock \emph{SIAM Journal on Scientific Computing}, 16\penalty0
  (1):\penalty0 151--168, 1995.
\newblock \doi{10.1137/0916010}.

\bibitem[McLachlan and Quispel(2002)]{McLachlanQuispel:2002}
Robert~I. McLachlan and G.~Reinout~W. Quispel.
\newblock {Splitting methods}.
\newblock \emph{Acta Numerica}, 11:\penalty0 341--434, 2002.
\newblock \doi{10.1017/S0962492902000053}.

\bibitem[Monk(2003)]{Monk:2003}
Peter Monk.
\newblock \emph{Finite Element Methods for Maxwell's Equations}.
\newblock Oxford University Press, 2003.

\bibitem[Moon et~al.(2015)Moon, Teixeira, and Omelchenko]{moon2015exact}
Haksu Moon, Fernando~L. Teixeira, and Yuri~A. Omelchenko.
\newblock Exact charge-conserving scatter-gather algorithm for particle-in-cell
  simulations on unstructured grids: A geometric perspective.
\newblock \emph{Computer Physics Communications}, 194:\penalty0 43--53, 2015.
\newblock \doi{10.1016/j.cpc.2015.04.014}.

\bibitem[Moore and Reich(2003)]{Moore:2003}
Brian~E. Moore and Sebastian Reich.
\newblock {Multi-symplectic integration methods for Hamiltonian PDEs}.
\newblock \emph{Future Generation Computer Systems}, 19\penalty0 (3):\penalty0
  395--402, 2003.
\newblock \doi{10.1016/S0167-739X(02)00166-8}.

\bibitem[Morita(2001)]{Morita:2001}
Shigeyuki Morita.
\newblock \emph{Geometry of Differential Forms}.
\newblock American Mathematical Society, 2001.

\bibitem[Morrison(1980)]{Morrison:1980}
Philip~J. Morrison.
\newblock The {Maxwell}--{Vlasov} equations as a continuous {Hamiltonian}
  system.
\newblock \emph{Physics Letters A}, 80\penalty0 (5--6):\penalty0 383--386,
  1980.
\newblock \doi{10.1016/0375-9601(80)90776-8}.

\bibitem[Morrison(1981{\natexlab{a}})]{Morrison:1981}
Philip~J. Morrison.
\newblock {Hamiltonian field description of the one-dimensional Poisson--Vlasov
  equations}.
\newblock Technical Report PPPL--1788, Princeton Plasma Physics Laboratory,
  1981{\natexlab{a}}.

\bibitem[Morrison(1981{\natexlab{b}})]{Morrison:1981b}
Philip~J. Morrison.
\newblock {Hamiltonian Field Description of the Two-Dimensional Vortex Fluids
  and Guiding Center Plasmas}.
\newblock Technical Report PPPL--1783, Princeton Plasma Physics Laboratory,
  1981{\natexlab{b}}.

\bibitem[Morrison(1982)]{Morrison:1982}
Philip~J. Morrison.
\newblock Poisson brackets for fluids and plasmas.
\newblock In \emph{AIP Conference Proceedings}, volume~88, pages 13--46. AIP,
  1982.
\newblock \doi{10.1063/1.33633}.

\bibitem[Morrison(1986)]{Morrison:1986}
Philip~J. Morrison.
\newblock {A paradigm for joined Hamiltonian and dissipative systems}.
\newblock \emph{Physica D: Nonlinear Phenomena}, 18:\penalty0 410--419, 1986.
\newblock \doi{10.1016/0167-2789(86)90209-5}.

\bibitem[Morrison(1987)]{Morrison:1987}
Philip~J. Morrison.
\newblock {Variational Principle and Stability of Nonmonotonic Vlasov--Poisson
  Equilibria}.
\newblock \emph{Zeitschrift für Naturforschung A}, 42:\penalty0 1115--1123,
  1987.
\newblock \doi{10.1515/zna-1987-1009}.

\bibitem[Morrison(1998)]{Morrison:1998}
Philip~J. Morrison.
\newblock Hamiltonian description of the ideal fluid.
\newblock \emph{Reviews of Modern Physics}, 70:\penalty0 467--521, Apr 1998.
\newblock \doi{10.1103/RevModPhys.70.467}.

\bibitem[Morrison(2013)]{Morrison:2013}
Philip~J. Morrison.
\newblock A general theory for gauge-free lifting.
\newblock \emph{Physics of Plasmas}, 20\penalty0 (1):\penalty0 012104, 2013.
\newblock \doi{10.1063/1.4774063}.

\bibitem[Morrison and Pfirsch(1989)]{MorrisonPfirsch:1989}
Philip~J. Morrison and Dieter Pfirsch.
\newblock {Free-energy expressions for Vlasov equilibria}.
\newblock \emph{Physical Review A}, 40:\penalty0 3898--3910, 1989.
\newblock \doi{10.1103/PhysRevA.40.3898}.

\bibitem[Munz et~al.(1999)Munz, Schneider, Sonnendr{\"u}cker, and
  Vo{\ss}]{munz1999}
Claus-Dieter Munz, Rudolf Schneider, Eric Sonnendr{\"u}cker, and Ursula
  Vo{\ss}.
\newblock Maxwell's equations when the charge conservation is not satisfied.
\newblock \emph{Comptes Rendus de l'Acad{\'e}mie des Sciences -- Series I --
  Mathematics}, 328\penalty0 (5):\penalty0 431--436, 1999.
\newblock \doi{10.1016/S0764-4442(99)80185-2}.

\bibitem[Munz et~al.(2000)Munz, Omnes, Schneider, Sonnendr{\"u}cker, and
  Voss]{munz2000divergence}
Claus-Dieter Munz, Pascal Omnes, Rudolf Schneider, Eric Sonnendr{\"u}cker, and
  Ursula Voss.
\newblock Divergence correction techniques for {M}axwell solvers based on a
  hyperbolic model.
\newblock \emph{Journal of Computational Physics}, 161:\penalty0 484--511,
  2000.
\newblock \doi{10.1006/jcph.2000.6507}.

\bibitem[Nakamura and Yabe(1999)]{NakamuraYabe:1999}
Takashi Nakamura and Takashi Yabe.
\newblock Cubic interpolated propagation scheme for solving the
  hyper--dimensional {Vlasov}--{Poisson} equation in phase space.
\newblock \emph{Journal of Computational Physics}, 120:\penalty0 122--154,
  1999.
\newblock \doi{10.1016/S0010-4655(99)00247-7}.

\bibitem[Palha et~al.(2014)Palha, Rebelo, Hiemstra, Kreeft, and
  Gerritsma]{Palha:2014}
Artur Palha, Pedro~Pinto Rebelo, Ren{\'e} Hiemstra, Jasper Kreeft, and Marc
  Gerritsma.
\newblock {Physics-compatible discretization techniques on single and dual
  grids, with application to the Poisson equation of volume forms}.
\newblock \emph{Journal of Computational Physics}, 257:\penalty0 1394--1422,
  2014.
\newblock \doi{10.1016/j.jcp.2013.08.005}.

\bibitem[Qin et~al.(2015)Qin, He, Zhang, Liu, Xiao, and Wang]{Qin:2015}
Hong Qin, Yang He, Ruili Zhang, Jian Liu, Jianyuan Xiao, and Yulei Wang.
\newblock {Comment on ``Hamiltonian splitting for the Vlasov--Maxwell
  equations''}.
\newblock \emph{Journal of Computational Physics}, 297:\penalty0 721--723,
  2015.
\newblock \doi{10.1016/j.jcp.2015.04.056}.

\bibitem[Qin et~al.(2016)Qin, Liu, Xiao, Zhang, He, Wang, Sun, Burby, Ellison,
  and Zhou]{Qin:2016}
Hong Qin, Jian Liu, Jianyuan Xiao, Ruili Zhang, Yang He, Yulei Wang, Yajuan
  Sun, Joshua~W. Burby, Leland Ellison, and Yao Zhou.
\newblock Canonical symplectic particle-in-cell method for long-term
  large-scale simulations of the {Vlasov}--{Maxwell} equations.
\newblock \emph{Nuclear Fusion}, 56\penalty0 (1):\penalty0 014001, 2016.
\newblock \doi{10.1088/0029-5515/56/1/014001}.

\bibitem[Ratnani and Sonnendr\"ucker(2012)]{Ratnani:2012}
Ahmed Ratnani and Eric Sonnendr\"ucker.
\newblock An arbitrary high-order spline finite element solver for the time
  domain {Maxwell} equations.
\newblock \emph{Journal of Scientific Computing}, 51\penalty0 (1):\penalty0
  87--106, 2012.
\newblock \doi{10.1007/s10915-011-9500-8}.

\bibitem[Reich(2000)]{Reich:2000}
Sebastian Reich.
\newblock {Multi-Symplectic Runge--Kutta Collocation Methods for Hamiltonian
  Wave Equations}.
\newblock \emph{Journal of Computational Physics}, 157\penalty0 (2):\penalty0
  473--499, 2000.
\newblock \doi{10.1006/jcph.1999.6372}.

\bibitem[Sanz-Serna and Calvo(1993)]{SanzSerna:1993}
Jes{\'u}s~Maria Sanz-Serna and Mari-Paz Calvo.
\newblock {Symplectic Numerical Methods for Hamiltonian Problems}.
\newblock \emph{International Journal of Modern Physics C}, 04\penalty0
  (02):\penalty0 385--392, 1993.
\newblock \doi{10.1142/S0129183193000410}.

\bibitem[Scovel and Weinstein(1994)]{Scovel:1994}
Clint Scovel and Alan Weinstein.
\newblock {Finite dimensional Lie-Poisson approximations to Vlasov--Poisson
  equations}.
\newblock \emph{Communications on Pure and Applied Mathematics}, 47\penalty0
  (5):\penalty0 683--709, 1994.
\newblock \doi{10.1002/cpa.3160470505}.

\bibitem[Shadwick et~al.(2014)Shadwick, Stamm, and Evstatiev]{Shadwick:2014}
Bradley~A. Shadwick, Alexander~B. Stamm, and Evstati~G. Evstatiev.
\newblock {Variational formulation of macro-particle plasma simulation
  algorithms}.
\newblock \emph{Physics of Plasmas}, 21\penalty0 (5):\penalty0 055708, 2014.
\newblock \doi{10.1063/1.4874338}.

\bibitem[Sircombe and Arber(2009)]{sircombe2009split}
Nathan~J. Sircombe and Tony~D. Arber.
\newblock {VALIS: A split-conservative scheme for the relativistic 2D
  Vlasov--Maxwell system}.
\newblock \emph{Journal of Computational Physics}, 228:\penalty0 4773--4788,
  2009.
\newblock \doi{10.1016/j.jcp.2009.03.029}.

\bibitem[Sonnendr{\"u}cker(2017)]{Sonnendruecker:2017}
Eric Sonnendr{\"u}cker.
\newblock \emph{{Numerical Methods for the Vlasov--Maxwell Equations}}.
\newblock Springer, 2017.

\bibitem[Squire et~al.(2012)Squire, Qin, and Tang]{Squire:2012}
Jonathan Squire, Hong Qin, and William~M. Tang.
\newblock Geometric integration of the {Vlasov}--{Maxwell} system with a
  variational particle-in-cell scheme.
\newblock \emph{Physics of Plasmas}, 19:\penalty0 084501, 2012.
\newblock \doi{10.1063/1.4742985}.

\bibitem[Stamm and Shadwick(2014)]{Stamm:2014}
Alexander~B. Stamm and Bradley~A. Shadwick.
\newblock {Variational formulation of macroparticle models for electromagnetic
  plasma simulations}.
\newblock \emph{IEEE Transactions on Plasma Science}, 42\penalty0 (6):\penalty0
  1747--1758, 2014.
\newblock \doi{10.1109/TPS.2014.2320461}.

\bibitem[Stern et~al.(2014)Stern, Tong, Desbrun, and Marsden]{Stern:2014}
Ari Stern, Yiying Tong, Mathieu Desbrun, and Jerrold~E. Marsden.
\newblock Geometric computational electrodynamics with variational integrators
  and discrete differential forms.
\newblock In \emph{Geometry, mechanics, and dynamics: the legacy of {Jerry
  Marsden}}, Fields Institute Communications. Springer, 2014.
\newblock \doi{10.1007/978-1-4939-2441-7_19}.

\bibitem[Strang(1968)]{Strang:1968}
Gilbert Strang.
\newblock {On the Construction and Comparison of Difference Schemes}.
\newblock \emph{SIAM Journal on Numerical Analysis}, 5\penalty0 (3):\penalty0
  506--517, 1968.
\newblock \doi{10.1137/0705041}.

\bibitem[Suzuki(1990)]{Suzuki:1990}
Masuo Suzuki.
\newblock {Fractal decomposition of exponential operators with applications to
  many-body theories and Monte Carlo simulations}.
\newblock \emph{Physics Letters A}, 146\penalty0 (6):\penalty0 319--323, 1990.
\newblock \doi{10.1016/0375-9601(90)90962-N}.

\bibitem[Trotter(1959)]{Trotter:1959}
Hale~F. Trotter.
\newblock {On the product of semi-groups of operators}.
\newblock \emph{Proceedings of the American Mathematical Society}, 10\penalty0
  (4):\penalty0 545--551, 1959.
\newblock \doi{10.1090/S0002-9939-1959-0108732-6}.

\bibitem[Umeda et~al.(2003)Umeda, Omura, Tominaga, and Matsumoto]{umeda2003}
Takayuki Umeda, Yoshiharu Omura, T.~Tominaga, and Hiroshi Matsumoto.
\newblock A new charge conservation method in electromagnetic particle-in-cell
  simulations.
\newblock \emph{Computer Physics Communications}, 156\penalty0 (1):\penalty0
  73--85, 2003.
\newblock \doi{10.1016/S0010-4655(03)00437-5}.

\bibitem[Villasenor and Buneman(1992)]{Villasenor.Buneman.1992.cpc}
John Villasenor and Oscar Buneman.
\newblock Rigorous charge conservation for local electromagnetic field solvers.
\newblock \emph{Computer Physics Communications}, 69\penalty0 (2):\penalty0
  306--316, 1992.
\newblock \doi{10.1016/0010-4655(92)90169-Y}.

\bibitem[Warnick and Russer(2006)]{Warnick:2006}
Karl~F. Warnick and Peter~H. Russer.
\newblock {Two, Three and Four-Dimensional Electromagnetics Using Differential
  Forms}.
\newblock \emph{Turkish Journal Of Electrical Engineering {\&} Computer
  Sciences}, 14\penalty0 (1):\penalty0 153--172, 2006.

\bibitem[Warnick et~al.(1998)Warnick, Selfridge, and Arnold]{Warnick:1997}
Karl~F. Warnick, Richard~H. Selfridge, and David~V. Arnold.
\newblock {Teaching Electromagnetic Field Theory Using Differential Forms}.
\newblock \emph{IEEE Transactions on Education}, 1:\penalty0 53--68, 1998.
\newblock \doi{10.1109/13.554670}.

\bibitem[Weibel(1959)]{Weibel:1959}
Erich~S. Weibel.
\newblock Spontaneously growing transverse waves in a plasma due to an
  anisotropic velocity distribution.
\newblock \emph{Physical Review Letters}, 2:\penalty0 83--84, 1959.
\newblock \doi{10.1103/PhysRevLett.2.83}.

\bibitem[Weinstein and Morrison(1981)]{Weinstein:1981}
Alan Weinstein and Philip~J. Morrison.
\newblock {Comments on: The Maxwell--Vlasov equations as a continuous
  Hamiltonian system}.
\newblock \emph{Physics Letters A}, 86\penalty0 (4):\penalty0 235--236, 1981.
\newblock \doi{10.1016/0375-9601(81)90496-5}.

\bibitem[Xiao et~al.(2015)Xiao, Qin, Liu, He, Zhang, and Sun]{Xiao:2015}
Jianyuan Xiao, Hong Qin, Jian Liu, Yang He, Ruili Zhang, and Yajuan Sun.
\newblock Explicit high-order non-canonical symplectic particle-in-cell
  algorithms for {Vlasov}--{Maxwell} systems.
\newblock \emph{Physics of Plasmas}, 22:\penalty0 112504, 2015.
\newblock \doi{10.1063/1.4935904}.

\bibitem[Ye and Morrison(1992)]{Ye:1992}
Huanchun Ye and Philip~J. Morrison.
\newblock Action principles for the {Vlasov} equation.
\newblock \emph{Physics of Fluids B}, 4\penalty0 (4):\penalty0 771--777, 1992.
\newblock \doi{10.1063/1.860231}.

\bibitem[Yee(1966)]{Yee:1966}
Kane~S. Yee.
\newblock {Numerical solution of initial boundary value problems involving
  {Maxwell}'s equations in isotropic media}.
\newblock \emph{IEEE Trans. Antennas Propag.}, 14:\penalty0 302--307, 1966.
\newblock \doi{10.1109/TAP.1966.1138693}.

\bibitem[Yoshida(1990)]{Yoshida:1990}
Haruo Yoshida.
\newblock {Construction of higher order symplectic integrators}.
\newblock \emph{Physics Letters A}, 150\penalty0 (5-7):\penalty0 262--268,
  1990.
\newblock \doi{10.1016/0375-9601(90)90092-3}.

\bibitem[Yu et~al.(2013)Yu, Jin, Zhou, Li, and Gu]{yu2013high}
Jinqing Yu, Xiaolin Jin, Weimin Zhou, Bin Li, and Yuqiu Gu.
\newblock High-order interpolation algorithms for charge conservation in
  particle-in-cell simulations.
\newblock \emph{Communications in Computational Physics}, 13:\penalty0
  1134--1150, 2013.
\newblock \doi{10.4208/cicp.290811.050412a}.

\end{thebibliography}

\end{document}